\documentclass[twoside,11pt]{article}

\usepackage{amsthm}
\usepackage
{jmlr2e}

\usepackage{hyperref}
\usepackage{amsmath, amssymb, amsfonts, mathtools}

\usepackage{graphicx}
\usepackage{subfigure}
\usepackage{booktabs}
\usepackage[textsize=small]{todonotes}
\usepackage{multirow}
\usepackage{bbm}

\DeclareMathOperator*{\argmax}{arg\,max}
\DeclareMathOperator*{\argmin}{arg\,min}
\DeclareMathOperator*{\supp}{\operatorname{supp}}

\theoremstyle{plain}
\newtheorem{theor}{Theorem}[section]
\newtheorem{prop}[theor]{Proposition}
\newtheorem{lem}[theor]{Lemma}
\newtheorem{corol}[theor]{Corollary}

\newcommand{\lav}[1]{\left\langle#1\right\rangle_x}
\newcommand{\sgn}[1]{\operatorname{sgn}\left(#1\right)}

\usepackage{lastpage}

\jmlrheading{25}{2024}{1-\pageref{LastPage}}{5/23; Revised
12/23}{3/24}{23-0698}{Maksim Velikanov and Dmitry Yarotsky}

\ShortHeadings{Tight Convergence Rate Bounds for Optimization}{Velikanov and Yarotsky}
\firstpageno{1}

\begin{document}

\title{Tight Convergence Rate Bounds for Optimization Under Power Law Spectral Conditions}

\author{\name Maksim Velikanov \email maksim.velikanov@tii.ae \\
\addr   Technology Innovation Institute, Abu Dhabi, UAE; \\
        CMAP, Ecole Polytechnique, Paris, France
        \AND
        \name Dmitry Yarotsky \email d.yarotsky@skoltech.ru \\
       \addr Center for Artificial Intelligence Technology\\
       Skolkovo Institute of Science and Technology\\
       Moscow, Russia}

\editor{Mehryar Mohri}
\maketitle

\begin{abstract}

Performance of optimization on quadratic problems sensitively depends on the low-lying part of the spectrum. For large (effectively infinite-dimensional) problems, this part of the spectrum can often be naturally represented or approximated by power law distributions, resulting in power law convergence rates for iterative solutions of these problems by gradient-based algorithms. In this paper, we propose a new spectral condition providing tighter upper bounds for problems with power law optimization trajectories. We use this condition to build a complete picture of upper and lower bounds for a wide range of optimization algorithms -- Gradient Descent, Steepest Descent, Heavy Ball, and Conjugate Gradients -- with an emphasis on the underlying schedules of learning rate and momentum. In particular, we demonstrate how an optimally accelerated method, its schedule, and convergence upper bound can be obtained in a unified manner for a given shape of the spectrum. Also, we provide first proofs of tight lower bounds for convergence rates of Steepest Descent and Conjugate Gradients under spectral power laws with general exponents. Our experiments show that the obtained convergence bounds and acceleration strategies are not only relevant for exactly quadratic optimization problems, but also fairly accurate when applied to the training of neural networks.  

\end{abstract}

\begin{keywords}
Gradient Descent, Steepest Descent, Heavy Ball,  Conjugate Gradients, power-law spectrum, convergence rate, tight bounds, non-strongly-convex least squares, acceleration, neural networks 
\end{keywords}

\newpage

\renewcommand{\baselinestretch}{0.75}\normalsize
\tableofcontents
\renewcommand{\baselinestretch}{1.0}\normalsize

\section{Introduction}\label{sec:intro}
Modern large-scale optimization problems, such as training of neural networks, are typically solved by some variants of Gradient Descent (GD) or its accelerated versions. Examples of such methods include Stochastic Gradient Descent (SGD), GD with momentum \citep{polyak1964some, qian1999momentum},  Nesterov's accelerated gradient \citep{nesterov1983method}, Conjugate Gradients (CG, \citet{Hestenes&Stiefel:1952}), ADAM \citep{kingma2014adam}. These first order methods strike a good balance between universality, efficiency and complexity, which is crucial for high-dimensional applications (where, for example, higher-order Hessian-based methods would be prohibitively expensive).

While real world optimization problems can be characterized by a multitude of different aspects, the key features of first order methods are well captured by examining optimization of  quadratic loss functions $L(\mathbf w)=\tfrac{1}{2}\mathbf w^TA\mathbf w-\mathbf w^T\mathbf b$, that typically serve as reasonable approximation to the actual objective functions near local or global minima. The main challenge in optimizing such quadratic losses is their ill-conditioning, i.e. some of the eigenvalues of $A$ being much smaller than the others. The convergence rate of the optimization is determined by the condition number of $A$ and notably degrades as this number tends to infinity. The extreme case is when the condition number is effectively infinite, i.e. the eigenvalues of $A$ can be arbitrarily small. In this case, there are well-known classical bounds (see e.g. section 6.1 of \cite{Polyak87}) for the convergence rate in terms of the initial error $\|\mathbf w_0-\mathbf w_*\|$, where $\mathbf w_0$ is the starting point and $\mathbf w_*$ is the minimizer. Specifically, for the vanilla GD $\mathbf w_{n+1}=\mathbf w_{n}-\alpha\nabla L(\mathbf w_n)$ with learning rate $\alpha<2/\lambda_{\max}$, where $\lambda_{\max}$ denotes the largest eigenvalue of $A$,  we have 
\begin{equation}\label{eq:gd_bound_polyak}L(\mathbf w_n)-L(\mathbf w_*)\le \frac{\|\mathbf w_0-\mathbf w_*\|^2}{4\alpha e n}.
\end{equation}  
For optimization by Conjugate Gradients (CG), we have
\begin{equation}\label{eq:cgd_bound_polyak}
L(\mathbf w_n)-L(\mathbf w_*)\le \frac{\lambda_{\max}\|\mathbf w_0-\mathbf w_*\|^2}{2(2n+1)^2}.
\end{equation} 
These bounds suggest, in particular, that the convergence is $O(n^{-1})$ for GD, and $O(n^{-2})$ for CG.

\begin{figure*}[t]
    \centering
    \includegraphics[width=0.99\textwidth,clip,trim= 4mm 97mm 4mm 4mm ]{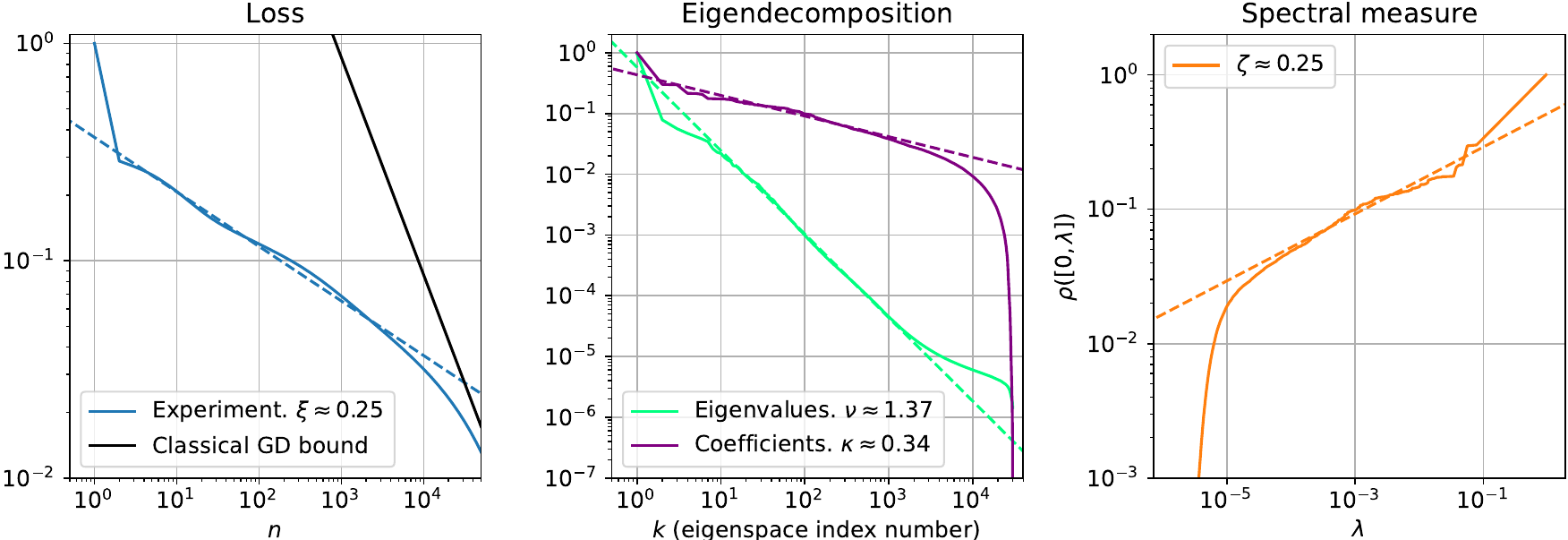}
    
    \vspace{-1mm}
    \caption{Spectral properties and GD loss of a MNIST classifier learned in a kernel regime. \textbf{Left:} The experimental trajectory $L(\mathbf w_n)$ of GD loss, the fitted power law \eqref{eq:lpropto}, and the classical $O(n^{-1})$ bound \eqref{eq:gd_bound_polyak}. \textbf{Center:} The eigenvalues $\lambda_k\propto k^{-\nu},\nu\approx1.37,$ and the cumulative distribution of target expansion coefficients, $\sum_{s=k}^{k_{\max}} \lambda_s c_s^2\propto k^{-\kappa}, \kappa\approx 0.34$. The target expansion coefficients $c_k$ are defined by the spectral expansion $\mathbf w_*=\sum_{k}c_k\mathbf e_k$ of the minimizer vector $\mathbf w_*$ over the eigenvectors $\mathbf e_k$. \textbf{Right:} The spectral measure $\rho([0,\lambda])=\sum_{k:\lambda_k<\lambda}\lambda_kc_k^2\propto \lambda^{\zeta}, \xi=\zeta=\tfrac{\kappa}{\nu}\approx 0.25$. See Section \ref{sec:probdef} for a general definition.} 
  \label{fig:MNIST_loss_and_spectral_distributions}
\end{figure*}

However, bounds \eqref{eq:gd_bound_polyak}, \eqref{eq:cgd_bound_polyak} are crude in that they do not use any information about the distribution of eigenvalues in the segment $[0,\lambda_{\max}]$ and about the expansion coefficients of the initial displacement $\mathbf w_0-\mathbf w_*$ over the eigenbasis of $A$. As a results, actual convergence rates in practical problems can be drastically different from  the above $O(n^{-1})$ or $O(n^{-2})$. In fact, the experimentally observed convergence can even be slower than $O(n^{-1})$, seemingly contradicting the theory. In Figure \ref{fig:MNIST_loss_and_spectral_distributions} (left) we show the loss trajectory of a neural network in a very basic example -- learning the standard MNIST digit classifier \citep{lecun2010mnist} by basic GD in a kernel regime (see Section \ref{sec:exp_details} for details). We see that up to very late iterations, loss evolves  as 
\begin{equation}\label{eq:lpropto}
    L(\mathbf w_n)\propto n^{-\xi},\quad \xi\approx 0.25.
\end{equation} This power law  can be explained theoretically by observing that  both the eigenvalue distribution and the cumulative distribution of target expansion coefficients in this problem are also approximate power laws, with exponents $\kappa\approx 0.34$ and $\nu\approx 1.35$, respectively (see Figure \ref{fig:MNIST_loss_and_spectral_distributions} (center) and later sections for details). The power law \eqref{eq:lpropto} for the loss can then be derived from these spectral laws with the exponent given by $\xi=\tfrac{\kappa}{\nu}\approx 0.25$.

The apparent contradiction between the $O(n^{-1})$  theoretical bound \eqref{eq:gd_bound_polyak} and the much slower experimental convergence \eqref{eq:lpropto} is explained by the heavy tail of the  eigendecomposition of the fitted function. If we attempt to view the problem as effectively infinite-dimensional (which is convenient for abstract theory involving spectral power laws), then under condition $\tfrac{\kappa}{\nu}\le 1$ the minimizer $\mathbf w_*$ does not exist as a finite-norm vector (is ``unattainable''). Accordingly, the norm $\|\mathbf w_0-\mathbf w_*\|$ appearing in \eqref{eq:lpropto} becomes infinite and bound \eqref{eq:lpropto} becomes vacuous. If instead we treat the space as high- but finite-dimensional, then the norm $\|\mathbf w_0-\mathbf w_*\|$ is finite but very large, so that bound \eqref{eq:gd_bound_polyak} is still too crude to reflect the actual convergence. Note, at the same time, that in the context of predictive modeling we are primarily interested in the loss $L(\mathbf w_n)$ rather than the norm $\|\mathbf w_n-\mathbf w_*\|$, since the former directly reflects the performance of the models while the latter only characterizes convergence in terms of the internal structure of the model and depends on the model parameterization, the choice of the norm, etc. In the case of MNIST, despite large values of the norm $\|\mathbf w_n-\mathbf w_*\|$, the model trains well and achieves high accuracy even on the test set (see Section \ref{sec:exp}). This suggests that a theory describing training realistic machine learning models even as simple as a MNIST classifier need not in general assume existence of a finite-norm solution $\mathbf w_*$.

A power-law structure of the spectrum is a common property of many large-scale optimization problems, in particular in machine learning: see e.g. recent works \cite{cui2021generalization, bahri2021explaining,lee2020finite,canatar2021spectral,kopitkov2020neural,dou2021training,atanasov2021neural,bordelon2021learning,basri2020frequency,bietti2021approximation}. One particularly interesting modern scenario is optimization of neural networks in the ``infinitely wide'' NTK regime \citep{jacot2018neural} (or some other ``lazy training'' setting where the learning problem is linearized \citep{chizat2019lazy}). Ill-conditioning here results naturally from overparameterization. In the NTK regime, the neural network effectively becomes a linear model with an explicit kernel \citep{NEURIPS2019_0d1a9651}. This can be used to derive explicit power laws for the corresponding spectral distributions and GD convergence rates. For example, when fitting a $d$-variate indicator function by a ReLU network using the continuous-time GD, the leading term in the loss evolution can be found as $Cn^{-1/(d+1)}$ with some explicit constant $C$ \citep{velikanov2021explicit}. A number of recent works experimentally verify and exploit power law asymptotics of the kernel eigenvalues, e.g. for the analysis of generalization \citep{bahri2021explaining, canatar2021spectral, lee2020finite, jin2021learning}.

This shows that power law spectral conditions are natural assumptions for abstract optimization theory. The power-law structure of the spectrum is commonly described by ``source condition'' and ``capacity condition'' \citep{caponnetto2007optimal}. To the best of our knowledge, the first comprehensive study of several fundamental algorithms such as GD, CG and Heavy Ball (HB) in this setting was performed by Nemirovsky and Polyak who established a number of upper and lower bounds for convergence rates \citep{nemirovskiy1984iterative1, nemirovsky1984iterative2} under the source condition. Their work was later extended in various directions by multiple authors. In particular, \citet{brakhage1987ill} introduced HB with a special schedule based on Jacobi polynomials, providing improved convergence bounds. \citet{hanke1991accelerated, hanke1996asymptotics} pointed out several important connections between CG and the theory of orthogonal polynomial and proved tight lower bounds for convergence of CG in some special cases. \citet{gilyazov2013regularization} established upper bounds for convergence of the method of Steepest Descent (SD). In recent years, capacity and source conditions have been used in the context of kernel methods and Stochastic GD (SGD) to obtain power law convergence rate bounds $O(n^{-\xi})$ with different exponents $\xi$  \citep{berthier2020tight, zou2021benign, nitanda2021optimal, varre2021last}.  

\paragraph{Our contribution.} The present work is a comprehensive study of the fundamental first order optimization algorithms  GD, CG, HB and SD in problems with a power-law type of the spectrum. On the one hand, our aim is to paint a complete rigorous picture of attainable convergence rates. We consider separately the scenarios with constant, non-constant predefined, and adaptive learning rates. For each algorithm we prove a power-law upper bound and a matching lower bound showing that the upper bound is tight. On the other hand, we introduce a new type of assumption to describe problems with a power-law type of the spectrum. We show that our new assumption provides a more accurate description of convergence rates, and develop a methodology of working with it. We highlight now some particular contributions of our work.

\begin{enumerate}
    \item We give first general proofs of tight lower bounds for SD and CG, which were previously missing in the literature. This completes the full picture of upper and lower bounds for all considered algorithms.
    \item For optimization problems with a power-law loss asymptotic $L(\mathbf{w}_n)\sim C n^{-\zeta}$, we show that the optimal upper bound under the classical source condition acquires an additional logarithmic factor $O(n^{-\zeta}\log n)$, while the upper bound under our spectral assumption recovers the correct rate $O(n^{-\zeta})$.  
    \item Our spectral assumption naturally treats attainable and unattainable problems in a unified way, in particular covering practical scenarios in which loss converges as a power law with an exponent close to 0, like in the above MNIST example.   
    \item We show that our spectral assumption simplifies the logic of derivation of optimally accelerated gradient descent methods. As a byproduct, we give a new simple expression for an optimal HB schedule. 
    \item Our experiments show that the considered accelerated gradient descent methods may achieve their theoretically expected convergence rate $O(n^{-2\zeta})$ for practical quadratic problems as well as for non-linear optimization of neural networks. 
\end{enumerate}

\paragraph{Paper organization.} We describe our assumptions and optimization algorithms in Section \ref{sec:setting}.  In Section \ref{sec:overview} we summarize and briefly discuss our results. Detailed statements of the theoretical upper and lower convergence bounds are presented in Section \ref{sec:detail}. In Section \ref{sec:tighterub} we accurately compare upper bounds obtained for the same problem using either our spectral condition or the classical source condition. 
In Section \ref{sec:exp} we present experiments with all our optimization algorithms, including applications to neural network training. An additional literature review and proof details are deferred to the appendix.

\section{The setting}\label{sec:setting}
\subsection{Problem definition and spectral assumptions}\label{sec:probdef}
\paragraph{The assumptions.} We assume that the optimized quadratic loss function $L$ is defined on a Hilbert space $\mathcal H$ by
\begin{equation}\label{eq:lwj}
    L(\mathbf w)= \frac{1}{2}\|J\mathbf w-\mathbf f_*\|^2=\frac{1}{2}\langle \mathbf w, A \mathbf w\rangle-\langle\mathbf w,\mathbf b\rangle+\frac{1}{2}\|\mathbf f_*\|^2,
\end{equation}
where $J:\mathcal H\to\widetilde{\mathcal H}$ is a bounded linear operator mapping $\mathcal H$ to another Hilbert space $\widetilde{\mathcal H}$, $\mathbf f_*\in \widetilde{\mathcal H}$, and
\begin{align}
    A = J^\dagger J: \mathcal H\to\mathcal H,\quad
    \mathbf b = J^\dagger \mathbf f_*\in\mathcal H
\end{align}
($J^\dagger$ denotes the adjoint operator). The spaces $\mathcal H$ and $\widetilde{\mathcal H}$ are, in general, infinite-dimensional.
The form \eqref{eq:lwj} of the quadratic function appears naturally in the setting where $J$ represents a linearized model fitting a target vector $\mathbf f_*$ (that, e.g., represents a large number of scalar measurements). If $J$ is written as a matrix, its columns correspond to different ``features'' used to predict the target. 

In the sequel, it is convenient to assume  that $\ker(J)=\{0\}$ and that the range $\operatorname{Ran}(J)$ is dense in $\widetilde{\mathcal H}$.\footnote{The extension to the general case is obtained easily by projecting or restricting all the vectors and operators to $\mathcal H\ominus\ker(J)$ in the space $\mathcal H$ and to $\operatorname{Ran}(J)$ in the space $\widetilde{\mathcal H}$.} This implies, in particular, that $\inf_{\mathbf w\in\mathcal H}L(\mathbf w)=0$.  

Along with $A$, consider the unitarily equivalent positive definite operator
\begin{equation}
    \widetilde A=JJ^\dagger: \widetilde{\mathcal H}\to\widetilde{\mathcal H}.
\end{equation}
Given $\mathbf f_*\in\widetilde{\mathcal H},$ there is a (unique, scalar-valued) associated spectral measure $\rho=\rho_{\widetilde A,\mathbf f_*}$  such that
\begin{equation}\label{eq:rhodef}
    \langle p(\widetilde A)\mathbf f_*,\mathbf f_*\rangle 
    = \int_{\mathbb R} p(\lambda)\rho(d\lambda)
\end{equation}
for any polynomial $p$; this relation can then be extended to general Borel functions (see e.g. \cite{Birman1987SpectralTO}). In particular, if $\widetilde{\mathcal H}$ is finite-dimensional or $\widetilde A$ is compact, then $\rho=\sum_{k=1}^{\dim \widetilde{\mathcal H}}c_k^2\delta_{\lambda_k},$ where $\lambda_k$ are the eigenvalues of $\widetilde A$, and $c_k$ are the respective coefficients in the expansion of $\mathbf f_*$ over the orthonormal eigenvectors of $\widetilde A$. The measure $\rho$ is finite ($\rho(\mathbb R)=\|\mathbf f_*\|^2<\infty$) and supported on the finite interval $[0,\lambda_{\max}]$, where $\lambda_{\max}=\|\widetilde A\|=\| A\|.$ 

Our \textbf{main (``target expansion'') spectral condition} is a growth condition on the cumulative distribution function of $\rho$: 
\begin{equation}\label{eq:rhozeta}
    \rho((0,\lambda]) \le Q \lambda^{\zeta},\quad \lambda\in[0,\lambda_{\max}],
\end{equation}
where  $Q$ and $\zeta$ are some positive constants. Note that this condition does not require $\widetilde A$ to have a discrete spectrum.  
It is sometimes convenient  to fix $\lambda_{\max}=1$ and $Q=1$ for brevity:  
\begin{equation}\label{eq:rhozeta1}
    \rho((0,\lambda]) \le \lambda^{\zeta},\quad \lambda\in[0,1].
\end{equation}
Results for general $\lambda_{\max}$ and $Q$ can be recovered by rescaling  $J\mapsto \lambda_{\max}^{1/2}J$ and $\mathbf f_*\mapsto Q^{1/2}\lambda_{\max}^{\zeta/2}\mathbf f_*$; in particular, the loss $L(\mathbf w_n)$ is simply multiplied by $Q\lambda_{\max}^{\zeta}$ (see Section \ref{sec:attain}). 
In the rest of the paper we will always assume that $\lambda_{\max}=1$, but occasionally keep the coefficient $Q$ (e.g., this will be convenient for comparison with the classical source condition).

Our \textbf{secondary (``eigenvalue decay'') spectral condition} assumes that the operator $\widetilde A$ is compact so that its spectrum is discrete, and that the sorted eigenvalues $\lambda_1\ge \lambda_2\ge\ldots>0$ obey 
\begin{equation}\label{eq:lkcknu}
 \lambda_k\le \Lambda k^{-\nu}
\end{equation}
with some positive constants $\Lambda,\nu$\footnote{Bounding eigenvalues from above may seem counter-intuitive as faster eigenvalue decay for fixed coefficients $c_k$ leads to slower loss convergence. However, our main spectral condition \eqref{eq:rhozeta1} forces the coefficients $c_k$ to decrease if we decrease the eigenvalues $\lambda_k$. Moreover, for all considered algorithms except CG the eigenvalue decay condition \eqref{eq:lkcknu} will not actually matter given condition \eqref{eq:rhozeta1}. For CG, faster eigenvalue decay leads to faster loss convergence, which justifies the $\leq$ sign in \eqref{eq:lkcknu}}. We will see that this condition will only matter for the algorithm CG, but not the other algorithms we consider (GD, SD and HB). 

We say that \textbf{the solution is attainable} if there exists $\mathbf w_*\in \mathcal H$ such that $L(\mathbf w_*)=0$. In terms of the measure $\rho$, attainabilitity means that $\|\mathbf w_*\|^2=\|J^{-1}\mathbf f_*\|^2=\int \lambda^{-1}\rho(d\lambda)<\infty$. This holds if $\zeta>1$ in Eq. \eqref{eq:rhozeta}, and generally does not hold for $\zeta\le 1$ (see Section \ref{sec:attain}). We don't require attainability: in our setting $\zeta$ can be any positive number. 

\paragraph{Comparison with a standard ``source condition''.} Our ``target expansion'' condition \eqref{eq:rhozeta} is closely related to so-called ``source condition'' \citep{nemirovskiy1984iterative1, caponnetto2007optimal, berthier2020tight, varre2021last}, which is traditionally used to describe the problems with a power-law type of the spectral distributions. It is convenient to write this latter condition in the form 
\begin{equation}\label{eq:sourcecond}
    \| A^{-(\zeta'-1)/2}\mathbf w_*\|^2\le Q'
\end{equation} 
with some  parameters $\zeta', Q'$. This inequality can be written as an integral inequality w.r.t. the spectral measure $\rho$ given by Eq. \eqref{eq:rhodef}: using the identity $\| A^{a}\mathbf w_*\|^2=\int_0^\infty \lambda^{2a-1}\rho(d\lambda)$, 
\begin{equation}\label{eq:sourcecond1}
    \int_0^{\lambda_{\max}} \lambda^{-\zeta'}\rho(d\lambda) \leq Q'.
\end{equation}
Accordingly, the difference between our \eqref{eq:rhozeta} and classical \eqref{eq:sourcecond1} conditions is akin to the difference between $L^\infty$- and $L^1$-norm bounds. 

There is an approximate correspondence between the two conditions under which our exponent $\zeta$ matches the exponent $\zeta'$ of the classical condition. More precisely, let, as agreed, $\lambda_{\max}=1$. Denote by $\mathrm{P}(\zeta,Q)$ the set of all spectral measures $\rho$ on $[0,1]$ satisfying condition \eqref{eq:rhozeta}, and analogously denote by $\mathrm{P}'(\zeta',Q')$ the set of spectral measures $\rho$ on $[0,1]$ satisfying the classical source condition \eqref{eq:sourcecond1}. Then we prove (see section \ref{sec:attain})
\begin{lem} \label{lem:conditions_comp} Assuming $\zeta,\zeta',Q,Q'>0,$
\begin{align}\label{eq:ourcond_to_sourcecond}
    \hspace{7mm} \mathrm{P}(\zeta,Q) \subseteq \mathrm{P}'(\zeta',Q') \; &\iff \; \begin{cases}
        \zeta' < \zeta\\
        Q' \geq Q \frac{\zeta}{\zeta-\zeta'}
    \end{cases}\\
\label{eq:sourcecond_to_ourcond}
    \mathrm{P}'(\zeta',Q') \subseteq \mathrm{P}(\zeta,Q) \; &\iff \; \begin{cases}
        \zeta \leq \zeta'\\
        Q\geq Q'
    \end{cases}
\end{align}
\end{lem}
This lemma shows that our condition with parameters $\zeta, Q$ is slightly weaker than the classical source condition with the same parameters. In particular, while the classical condition with some exponent $\zeta'$ always implies our condition with the same exponent, the converse is not true: our condition with some $\zeta$ implies the classical condition only for $\zeta'<\zeta$, and the allowed constant $Q'\propto \tfrac{1}{\zeta-\zeta'}$ deteriorates as $\zeta'\nearrow \zeta$. Nevertheless, we will see that our weaker condition implies loss upper bounds analogous to those available with the classical condition.

\subsection{Optimization algorithms}\label{sec:algo}

We consider several classical iterative optimization algorithms \citep{Polyak87}. All of them are first-order in the sense that they use only the values of the loss function and its gradients from current and previous iterations. It will be convenient to assume that the starting point of these algorithms is $\mathbf w_0=0.$

\paragraph{Gradient Descent (GD)} is given by 
\begin{align}
    \mathbf w_{n+1}
    ={}&\mathbf w_n-\alpha_n\nabla L(\mathbf w_n)\\ 
    ={}&\mathbf w_n-\alpha_n(A \mathbf w_n-\mathbf b).\label{eq:discrete_gd}
\end{align}
We consider two scenarios for GD: the learning rate $\alpha_n$ either does not depend on $n$, or may depend on $n$, but with a schedule predefined prior to optimization and depending only on the exponent $\zeta$ from the main spectral condition \eqref{eq:rhozeta1}. 

\paragraph{Steepest Descent (SD)} is a modification of GD in which learning rate $\alpha_n$ is adaptively chosen at each iteration to optimize the loss: 
\begin{equation}
    \alpha_n = \argmin_\alpha L(\mathbf w_n-\alpha\nabla L(\mathbf w_n)).
\end{equation}
In our quadratic setting $\alpha_n$ can be explicitly written as
\begin{equation}\label{eq:SD_quadratic_choice}
    \alpha_n = \frac{\|\nabla L(\mathbf w_n)\|^2}{\langle A \nabla L(\mathbf w_n), \nabla L(\mathbf w_n)\rangle}.
\end{equation}

\paragraph{Heavy Ball (HB)} is a basic multi-step method (a.k.a. ``GD with momentum'') given by
\begin{align}
    \mathbf w_{n+1}
    ={}&\mathbf w_n-\alpha_n\nabla L(\mathbf w_n)+\beta_n(\mathbf w_{n}-\mathbf w_{n-1})\label{eq:mgd1}\\ 
    ={}&\mathbf w_n-\alpha_n(A \mathbf w_n-\mathbf b)+\beta_n(\mathbf w_{n}-\mathbf w_{n-1}) \label{eq:mgd2}
\end{align}
(for $n=0$ the term $\beta_n(\mathbf w_n-\mathbf w_{n-1})$ is dropped). As with GD, we assume that the learning rates $\alpha_n,\beta_n$ are either constant or predefined $n$-dependent. Throughout the paper we assume, as is common, that $0\le \beta_n<1.$

\paragraph{Conjugate Gradients (CG)} is defined by the same formula as HB, but (as with SD) with adaptively chosen learning rates minimizing the loss at each step:
\begin{align*}
    \alpha_n,\beta_n ={}& \argmin _{\alpha,\beta} L(\mathbf w_n-\alpha\nabla L(\mathbf w_n)+\beta(\mathbf w_{n}-\mathbf w_{n-1})).
\end{align*}
For a quadratic loss, the optimal $\alpha_n,\beta_n$ are given  by
\begin{align}
    \label{eq:CG_quadratic_choice_1}
    \alpha_n
    ={}&\frac{\|\mathbf r_n\|^2\langle A\mathbf p_n,\mathbf p_n\rangle-\langle \mathbf r_n,\mathbf p_n\rangle\langle A\mathbf r_n,\mathbf p_n\rangle}{\langle A\mathbf r_n,\mathbf r_n\rangle\langle A\mathbf p_n,\mathbf p_n\rangle-\langle A\mathbf r_n,\mathbf p_n\rangle},\\
    \label{eq:CG_quadratic_choice_2}
    \beta_n ={}& \frac{\|\mathbf r_n\|^2\langle A\mathbf r_n,\mathbf p_n\rangle-\langle \mathbf r_n,\mathbf p_n\rangle\langle A\mathbf r_n,\mathbf r_n\rangle}{\langle A\mathbf r_n,\mathbf r_n\rangle\langle A\mathbf p_n,\mathbf p_n\rangle-\langle A\mathbf r_n,\mathbf p_n\rangle},\\
    \mathbf r_n ={}& \nabla L(\mathbf w_n)=A\mathbf w_n-\mathbf b,\quad \mathbf p_n = \mathbf w_n-\mathbf w_{n-1}
\end{align}
(see \cite{Polyak87}, Section 3.2.2). The fundamental importance of CG lies in the fact that, for quadratic problems, this algorithm is optimal among all first order methods generating new iterates $\mathbf w_{n+1}$ by shifting the initial point $\mathbf w_0$ along linear subspaces spanned by the previously computed gradients  $\nabla L(\mathbf w_0),\ldots,\nabla L(\mathbf w_{n})$.

\section{Overview of  results}\label{sec:overview}

\begin{table*}[t]
\caption{Summary of convergence rates of $L(\mathbf w_n)$ for different algorithms and  learning rate schedules under spectral assumptions \eqref{eq:rhozeta} (i.e., $\rho((0,\lambda])=O(\lambda^{\zeta})$) and \eqref{eq:lkcknu} (i.e., $\lambda_k=O(k^{-\nu})$). Assumption \eqref{eq:lkcknu} only matters in the case of CG: in all other cases the spectrum does not even need to be discrete. For CG, if only assumption \eqref{eq:rhozeta} holds, then $L(\mathbf w_n)=O(n^{-2\zeta})$; if additionally \eqref{eq:lkcknu} holds, then $L(\mathbf w_n)=O(n^{-(2+\nu)\zeta})$. Each of the bounds in the table is tight in the sense that the respective exponents $\zeta, 2\zeta, (2+\nu)\zeta$ cannot be improved. The subscripts indicate the sections where the respective results are presented: roman for upper bounds, \textit{italic} for lower bounds, and \textbf{bold} for both.}
\label{tab:bounds}
\begin{center}
\begin{tabular}{lcccc}\toprule
&
\multicolumn{4}{c}{Learning rates} \\
\cmidrule(lr){2-5} 
 & \multirow{2}{12mm}{constant} & & predefined & \multirow{2}{12mm}{adaptive}
\\
& & & $n$-dependent & 
\\
\midrule
\multirow{2}{3cm}{Single-step} 
& \multicolumn{3}{c}{\small \bf Gradient Descent (GD)} & {\small \bf Steepest Descent (SD)}\\
& $O(n^{-\zeta})$\textsubscript{\bf \tiny \ref{sec:constant_lr_main}} & & $O(n^{-2\zeta})$\textsubscript{\tiny {\it \ref{sec:exact_power-law_measure}}, \ref{sec:ub}, {\it\ref{sec:lb}}} & $O(n^{-\zeta})$\textsubscript{\bf\tiny \ref{sec:sd}} \\
\midrule 
\multirow{2}{3cm}{Multi-step} 
& \multicolumn{3}{c}{\small \bf Heavy Ball (HB)} & {\small \bf Conjugate Gradients (CG)}\\
& $O(n^{-\zeta})$\textsubscript{\bf \tiny \ref{sec:constant_lr_main}} & & $O(n^{-2\zeta})$\textsubscript{\tiny {\it \ref{sec:exact_power-law_measure}}, \ref{sec:ub}, {\it\ref{sec:lb}}} & $O(n^{-2\zeta}) \;|\; O(n^{-(2+\nu)\zeta})$\textsubscript{\tiny  {\it \ref{sec:exact_power-law_measure}}, \ref{sec:ub}, {\it\ref{sec:lb}}} \\
\bottomrule

\end{tabular}
\end{center}
\end{table*}

\paragraph{The complete picture of upper and lower bounds.} Our main theoretical result is an essentially complete picture of optimal convergence rates, summarized in Table \ref{tab:bounds}: for each of the algorithms GD, SD, HB, CG, for each type of learning rate schedule (constant, predefined step-dependent, adaptive), for each $\zeta>0$ we establish an upper bound of the form $L(\mathbf w_n)=O(n^{-\xi})$ and a respective lower bound showing that the exponent $\xi$ cannot be improved.

In all cases except CG, only our primary condition \eqref{eq:rhozeta1} matters for the convergence rate: adding the eigenvalue decay condition does not affect the rate. This is confirmed by the lower bounds, which are constructed to satisfy both conditions. CG is an exceptional case where adding the eigenvalue decay condition allows to improve the upper bound from $O(n^{-2\zeta})$ to $O(n^{-(2+\nu)\zeta}).$ We prove that both these bounds are tight. 

Note that adaptivity of learning rates does not improve convergence rate for single-step methods (GD vs. SD), but does improve it for multi-step methods (HB vs. CG). The exponents $2\zeta$ of faster algorithms are twice as large as the exponents $\zeta$ of the basic ones (cf. \eqref{eq:gd_bound_polyak}, \eqref{eq:cgd_bound_polyak}). In a $d$-dimensional setting with finite $d$ CG finds the exact solution after $d$ iterations; the analog of this in our setting is the increased exponent $(2+\nu)\zeta$. 

Though theoretically CG has the highest convergence rate $O(n^{-(2+\nu)\zeta}),$ its practical implementation is not so efficient because of a fast accumulation of numerical errors. The indicated rate requires the polynomials associated with CG (see Section \ref{sec:poly}) to have roots very close to the eigenvalues of $A$, which imposes strong requirements on the precision of computations. Also, the $O(n^{-2\zeta})$ convergence of GD with predefined schedule is very sensitive to non-quadratic perturbations of the problem. See experiments in Section \ref{sec:exp}. 

In Table \ref{tab:bounds} we have four instances which enjoy convergence rates accelerated from $O(n^{-\zeta})$ to $O(n^{-2\zeta})$ or $O(n^{-(2+\nu)\zeta})$. In all these cases the stated rates are achieved using constructions based on Jacobi polynomials $P^{(a,b)}_n$  (see Section \ref{sec:exact_power-law_measure}).

The classical bounds $O(n^{-1})$ and $O(n^{-2})$ for GD and CG, respectively (cf. Eqs. \eqref{eq:gd_bound_polyak}, \eqref{eq:cgd_bound_polyak}), are, up to the coefficients, special cases of the bounds $O(n^{-\zeta})$ and $O(n^{-2\zeta})$ when $\|\mathbf w_*\|<\infty$, since by Lemma \ref{lem:conditions_comp} (specifically, by Eq. \eqref{eq:sourcecond_to_ourcond}) our main spectral condition \eqref{eq:rhozeta1} holds with $\zeta=1$ in this case. 

In the multi-step predefined step-dependent scenario, our lower bound (Theorem \ref{ther:discrete_spectra_predefined_schedule_LB}) applies to any method linearly expressing current step in terms of past gradients. Accordingly, this bound covers not only Heavy Ball, but also its modifications such as Nesterov Accelerated Gradient (NAG, \citet{nesterov1983method}). We discuss NAG in Appendix \ref{sec:relwork}. 

As already mentioned, most bounds of Table \ref{tab:bounds} (or some closely related bounds) already appeared in some form in earlier research \citep{nemirovskiy1984iterative1, nemirovsky1984iterative2, brakhage1987ill, hanke1991accelerated, hanke1996asymptotics, gilyazov2013regularization}, albeit under the stronger classical source assumption \eqref{eq:sourcecond}. Below we discuss various new elements of our work which were not present in earlier research. 

\paragraph{Optimization with unattainable solutions.} In almost all previous research, only the case of attainable solutions $\|\mathbf w_*\|<\infty$ (i.e., $\zeta>1$ in Eq. \eqref{eq:rhozeta}) is considered. However, as already pointed out in Section \ref{sec:intro}, even simple realistic problems such as MNIST have unattainable solutions. In fact, one can argue that this non-attainability is typical for a wide range of problems. In particular, it is shown in \cite{velikanov2021explicit} that in the $d$-dimensional kernel regression with kernels having homogeneous singularities of degree $\alpha$, the task of fitting indicator functions corresponds to the exponents  $\nu=1+\frac{\alpha}{d}, \kappa=\tfrac{1}{d}$ and $\zeta=\tfrac{\kappa}{\mu}=\tfrac{1}{d+\alpha}$ in Eqs. \eqref{eq:rhozeta}, \eqref{eq:lkcknu}. ReLU neural networks in the NTK regime are effectively such kernels models \citep{jacot2018neural} with $\alpha=1,$ so in these scenarios we always have $\zeta=\tfrac{1}{d+1}<1$. Our bounds in Table \ref{tab:bounds} are valid for all $\zeta>0$ and show that the non-attainability of the solution is not an obstacle for successful optimization. 

Even more importantly, our lower bounds show that, regardless of the optimization algorithm, the exponent $\zeta$ in the loss power law $L(\mathbf w_n)= O(n^{-\zeta})$ will, in general, be close to 0 if $\zeta$ is close to 0, i.e. the optimization will inevitably be quite slow. This agrees with experiment and dispels the excessively optimistic theoretical expectations such as $L(\mathbf w_n)= O(n^{-1})$ and $L(\mathbf w_n)= O(n^{-2})$ that one might get from Eqs. \eqref{eq:gd_bound_polyak}, \eqref{eq:cgd_bound_polyak}. 

\paragraph{New bounds.} Our significant new technical contributions are the tight lower bounds $\Omega(n^{-\zeta})$ and $\Omega(n^{-(2+\nu)\zeta})$ for SD and CG (see Theorems \ref{th:sd}, \ref{th:cglbdiscr}). While the respective upper bounds were known from \citet{nemirovskiy1984iterative1, hanke1991accelerated, hanke1996asymptotics, gilyazov2013regularization}, the tight lower bounds were not available under any kind of power-law assumption. 

We consider our lower bound $\Omega(n^{-(2+\nu)\zeta})$ for CG with discrete spectrum to be especially important, because CG can be viewed as an ``ultimate'' iterative first order algorithm for quadratic objectives: it essentially reconstructs the objective on the nested sequence of whole Krylov subspaces exhausting the space $\mathcal H$, and so in a sense optimally exploits all the iteratively available zero- and first-order information about the objective. Our lower bound then shows that even this optimal exploitation will not generally give fast convergence if $\zeta$ and $\nu$ are small (unless the problem or the algorithm are improved using some additional information about the problem  -- e.g., by pre-conditioning).  

To the best of our knowledge, the only previously available lower bounds for CG in the power-law setting were given in \cite{hanke1996asymptotics} and only covered two special cases $\nu=1,2$ for which explicit orthogonal polynomials are known. Our approach is completely different: for each $\zeta, \nu>0$ we give a simple explicit example of the operator $J$ and target $\mathbf f_*$ for which spectral conditions \eqref{eq:rhozeta}, \eqref{eq:lkcknu} hold and $L(\mathbf w_n)=\Omega(n^{-(2+\nu)\zeta})$ (assuming $\zeta\notin\mathbb Z$; see Theorem \ref{th:cglbdiscr}). 

Our tight lower bound $\Omega(n^{-\zeta})$ for SD also seems to be new. We give a simple proof based on the limiting periodic behavior of SD (Theorem \ref{th:sd}). 

Finally, our simple construction of the step-dependent schedule ensuring the improved convergence $L(\mathbf w_n)=O(n^{-2\zeta})$ for GD (see Theorem \ref{th:ubgdnonconst}) does not seem to have been described in earlier literature.

\paragraph{Tighter bounds: a weaker spectral assumption.}
As already mentioned in Section \ref{sec:setting}, our ``target expansion'' condition \eqref{eq:rhozeta} is a weaker version of a more standard ``source condition'' \eqref{eq:sourcecond}. In Section \ref{sec:tighterub} we show that this difference between conditions can play a significant role. Specifically, we show for the MNIST quadratic optimization problem that the upper bounds based on classical condition \eqref{eq:sourcecond} poorly describe the actual loss trajectory, while the upper bound based on our condition \eqref{eq:rhozeta} matches the true trajectory much better. We confirm this empirical observation theoretically for problems with a power-law loss trajectory $L(\mathbf{w}_n)\propto n^{-\zeta}$. We prove that in such problems, the upper bound based on the classical condition acquires an additional logarithmic factor, $O(n^{-\zeta}\log n)$, while the bound based on our condition retains the correct rate $O(n^{-\zeta})$.

\paragraph{Tighter bounds: specifying the constant.} A simplest example of an optimization problem exhibiting a $O(n^{-\xi})$ convergence rate is the exact power-law measure $\rho_\zeta([0,\lambda])=\lambda^\zeta$, a boundary case of our condition \eqref{eq:rhozeta1}. We interpret the loss of this problem, $L_n^{(\zeta)}$, as a reference point for convergence rates. It is easy to derive its full loss asymptotic $L_n^{(\zeta)}\stackrel{n\to\infty}{=}C n^{-\zeta}(1+o(1))$ with a specific constant $C$. Then, we are able to provide upper bounds that are quite close to this typical performance. For example, for GD and HB with constant learning rates, the bound asymptotically matches the typical performance: $L(\mathbf{w}_n)\leq L_n^{(\zeta)}(1+o(1))$. As for accelerated HB with the rate $O(n^{-2\zeta})$, the bound is just a few times larger than the typical performance, e.g. $L(\mathbf{w}_n)\leq 4 L_n^{(\zeta)}$ for $\zeta=1$.

\paragraph{Unified picture of convergence bounds and acceleration.} We develop a new, general and transparent approach to simultaneously obtain an upper and a matching lower loss bounds for most of the considered algorithms (except SD) under spectral conditions like \eqref{eq:rhozeta1} (see Sections \ref{sec:worst_case_measure} -- \ref{sec:ub}). This is done by relating the convergence for general problems satisfying condition $\rho([0,\lambda])\leq G(\lambda)$ to the convergence of a ``solvable'' problem with ``smooth'' spectral measure $\rho(d\lambda)=G'(\lambda)d\lambda$. For the solvable problem the optimal learning rate schedule can be found analytically through the 3-term recurrence relation of the related system of polynomials orthogonal with weight $\lambda G'(\lambda)d\lambda$. Then we show that this schedule remains efficient for all problems subject to $\rho([0,\lambda])\leq G(\lambda)$.

\section{Upper and lower bounds: detailed results} \label{sec:detail}

\begin{figure}
    \centering

    \begin{tikzpicture}[scale=1]
    \def\r{0.07}
    
    \node[rectangle,draw,minimum size=6mm] (WorstCase) at (-7,0) {\ref{sec:worst_case_measure}};
    \node[rectangle,draw,minimum size=6mm] (Constant) at (-5.5,0) {\ref{sec:constant_lr_main}};
    \node[rectangle,draw,minimum size=6mm] (Acceleration) at (-4,0) {\ref{sec:exact_power-law_measure}};
    \node[rectangle,draw,minimum size=6mm] (GenUpper) at (-2.5,0) {\ref{sec:ub}};
    \node[rectangle,draw,minimum size=6mm] (Lower) at (0,0) {\ref{sec:lb}};
    \node[rectangle,draw,minimum size=6mm] (SD) at (2.5,0) {\ref{sec:sd}};
    \node[rectangle,draw,minimum size=6mm] (Compare) at (-7,-1.5) {\ref{sec:tighterub}};

    \node[text width=6cm, align=center] () at (-4.75, 1.1)  {Upper and some lower bounds \\ for GD, HB, CG};
    \node[text width=3cm, align=center] () at (0, 0.9)  {Lower bounds \\ for GD, HB, CG};
    \node[text width=3cm, align=center] () at (2.5, 0.7)  {SD};
    \node[text width=4cm, align=center] () at (-5, -1.5)  {Comparison of\\ spectral conditions};

    \draw[->] (WorstCase) to (Constant); 
    \draw[->] (WorstCase) to [out=25,in=155] (Acceleration);
    \draw[->] (WorstCase) to [out=25,in=155] (GenUpper);
    \draw[->] (WorstCase) to (Compare);
    \draw[->] (Constant) to (Compare);
    \draw[->] (Acceleration) to (GenUpper);
    
    \end{tikzpicture}
    
    \caption{Logical dependencies between sections with main results.}
    \label{fig:logical}
\end{figure}

The structure of our exposition is shown in Figure \ref{fig:logical}. We start with a block of four sections establishing our methodology of working with spectral condition \eqref{eq:rhozeta1}. In Section \ref{sec:worst_case_measure}, we connect the loss convergence in a general problem described by our condition \eqref{eq:rhozeta1} with the convergence for an exact power-law measure. Then, in Section \ref{sec:constant_lr_main}, we use this connection to establish upper and lower bounds for constant learning rate algorithms. In Section \ref{sec:exact_power-law_measure}, we obtain accelerating strategies for the exact power-law spectral measure. Finally, in Section \ref{sec:ub}, we derive a number of upper bounds based on previously obtained accelerating strategy. In Section \ref{sec:lb} we derive the lower bounds for general algorithms with predefined schedules and CG applied to a problem with discrete spectrum \eqref{eq:lkcknu}. Lastly, in Section \ref{sec:sd}, we consider the SD algorithm, which requires tools and reasoning different from the other algorithms.

Our proofs rely heavily on the spectral representation of optimization by residual polynomials, which is recalled in Section \ref{sec:poly}. Each of the algorithms of Section \ref{sec:algo} is represented by a sequence of residual polynomials $p_n(\lambda), \; p_n(0)=1$  so that the $n$-step solution $\mathbf{f}_n=J\mathbf{w}_n$ satisfies
\begin{equation}
    \mathbf{f}_*-\mathbf{f}_n = p_n(\widetilde{A}) \mathbf{f}_*.
\end{equation}
The loss at step $n$ can be expressed through $p_n$ as 
\begin{equation}\label{eq:loss_through_residual_poly}
    L(\mathbf{w}_n) = \frac{1}{2}\int p_n(\lambda)^2 \rho(d\lambda).
\end{equation}

\subsection{Worst-case measures under main spectral condition}\label{sec:worst_case_measure}
In non-adaptive GD algorithms, the polynomials $p_n(\lambda)$ are fixed and independent of the problem's measure $\rho$. In this case, the worst case loss under a condition of the type \eqref{eq:rhozeta1} has a special structure revealed in the following theorem.

\begin{theor}[see proof in Section \ref{sec:proof_flattened}]\label{ther:worst_case_loss} Let $G(x)$ be a nondecreasing absolutely continuous function on $[0,1]$ such that $G(0)=0$, and let $q(x)$ be any nonnegative polynomial on $[0,1]$. Consider the integral $\int q(x) \rho(dx)$ as a functional on measures $\rho$ supported on $[0,1]$ and  satisfying $\rho([0,x])\leq G(x)$ for all $x\in[0,1]$. Then  the maximum of this functional is given by  
    \begin{align}\label{eq:worst_case_loss_ourcond}
        \sup\limits_{\substack{\rho:\; \operatorname{supp}(\rho)\subset[0,1],\\ \rho([0,x])\leq G(x)\forall x}}\int q(x) \rho(dx) = \int \overline{q}(x) G'(x)dx, \qquad \overline{q}(x) = \underset{y \geq x}{\sup} \; q(y). 
    \end{align}
\end{theor}
We will refer to $\overline{q}(x)$ as a ``flattened polynomial''. Considering the case $G(\lambda)=\lambda^\zeta$, we see that this theorem allows to reduce the analysis of upper bounds under condition \eqref{eq:rhozeta1} to estimating the averages of flattened polynomials $\overline{p^2_n}(x)$ over the exact power-law measure $\rho(d\lambda)=d(\lambda^\zeta)$.  

\begin{figure}
    \centering
    {\includegraphics[scale=0.6, trim=4mm 4mm 3mm 4mm, clip]{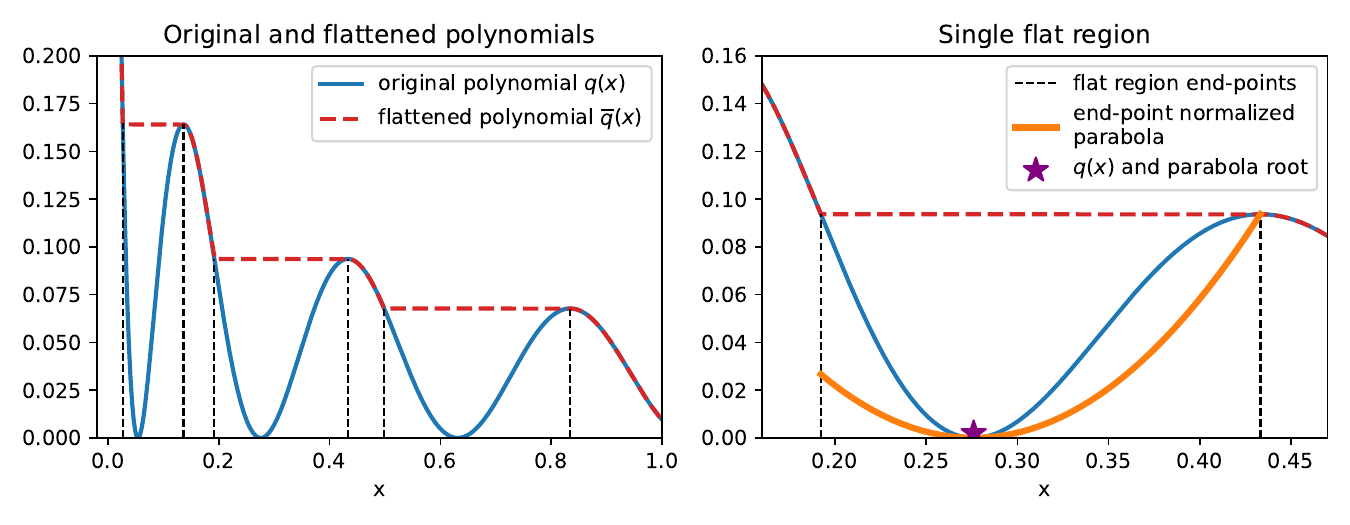}}
    
    \vspace{-1mm}
    \caption{The ``flattened polynomial'' $\overline{q}(x)$ associated with the worst-case spectral measure (see Theorem \ref{ther:worst_case_loss}). \textbf{Left:} The original polynomial $q(x)=p_n^2(x)$ and respective flattened polynomial $\overline{q}(x)$ for $p_n(x)$ associated with the Jacobi scheduled HB at step $n=7$ (see Section \ref{sec:exact_power-law_measure}). \textbf{Right:} Same as the left but zoomed in to a neighborhood of a single flat region of $\overline{q}(x)$. The orange parabola, placed at the respective root of $p_n(x)$ and normalized to match $q(x)$ at the right end of the flat region, is used to estimate the contribution of the flat region to the loss upper bound.}
    \label{fig:Flat_polynomial}
\end{figure}

The flattened polynomial $\overline{q}(x)$ can be simply characterized by considering the sequence of largest local maxima $0\leq x_1 < x_2 < \ldots x_m \leq 1$   of $q(x)$ on $[0,1]$ such that $\{q(x_i)\}_{i=1}^m$ is decreasing. Indeed, take any interval $[x_i, x_{i+1}]$ and denote $y_i\in (x_i,x_{i+1})$ the left most point such that $q(y_i)=q(x_{i+1})$. Then, it is straightforward to see that $\overline{q}(x)=q(x)$ on $[x_i,y_i]$ and $\overline{q}(x)=q(x_{i+1})$ on $[y_i,x_{i+1}]$, hence the name ``flattened''. See Figure \ref{fig:Flat_polynomial} for an illustration. 

Let us now outline the structure of convergence rate analysis that is suggested by  Theorem \ref{ther:worst_case_loss} and will be behind most of our bounds for the algorithms GD and HB. For our main spectral condition \eqref{eq:rhozeta1} we have $G(\lambda)=\lambda^\zeta$, and equation \eqref{eq:worst_case_loss_ourcond} leads to the exact power-law spectral measure $\rho_\zeta(d\lambda)=d(\lambda^\zeta)$ with cumulative distribution function
\begin{equation}\label{eq:exactpowerlaw2}
    \rho_\zeta((0,\lambda]) = \lambda^{\zeta}, \quad \lambda\in[0,1].
\end{equation}
For this measure, we define the \emph{pair} of the worst-case loss given by Eq. \eqref{eq:worst_case_loss} and the exact loss,
\begin{align}
    \label{eq:worst_case_loss}
    \overline{L_n^{(\zeta)}} &= \frac{1}{2}\int_0^1 \overline{p_n^2}(\lambda) d(\lambda^\zeta),\\
    \label{eq:exact_power_measure_loss}
    L_n^{(\zeta)} &= \frac{1}{2}\int_0^1 p_n^2(\lambda) d(\lambda^\zeta).
\end{align}
In our results, we will observe the following traits of this pair. First, the worst-case loss $\overline{L_n^{(\zeta)}}$ is not significantly worse than the exact power-law loss $L_n^{(\zeta)}$ and can be tightly bound to it. Then, $L_n^{(\zeta)}$ can be precisely described relying on a simple form of exact power-law measure $d(\lambda^\zeta)$ and properties of the chosen polynomials family $p_n(\lambda)$. Once the pair is characterized, we have the (tightest) upper bound $L(\mathbf{w}_n)\leq \overline{L_n^{(\zeta)}}$, and a general (e.g., without discreteness restriction \eqref{eq:lkcknu}) lower bound $L_n^{(\zeta)}$. 

\subsection[Constant learning rates]{Constant learning rates (Section \ref{sec:constant_lr})}\label{sec:constant_lr_main}
Suppose that the learning rate $\alpha_n\equiv\alpha>0$ and, if present, the momentum parameter $\beta_n\equiv\beta$. The respective residual polynomials for GD and HB  are given by (see Section \ref{sec:const_lr_HB_UB})
\begin{align}
    \label{eq:const_lr_GD_poly}
    &p_n(\lambda)=(1-\alpha \lambda)^n, \\
    \label{eq:const_lr_HB_poly}
    &p_n(\lambda)=\big(\sqrt{\beta}\big)^n\Big(U_n(z) -\sqrt{\beta}U_{n-1}(z)\Big), \quad z(\lambda)=\frac{1+\beta-\alpha\lambda}{2\sqrt{\beta}},
\end{align}
where $U_n$ are the Chebyshev polynomials of the second kind.

Following our strategy described in Section \ref{sec:worst_case_measure}, we analyze the pair of losses $\overline{L^{(\zeta)}_n}, \; L^{(\zeta)}_n$ for the constant learning rate GD and HB characterized by Eqs. \eqref{eq:const_lr_GD_poly} and \eqref{eq:const_lr_HB_poly}.
\begin{theor}\label{ther:constant_lr_bounds}
Define residual polynomials $p_n(\lambda)$ with \eqref{eq:const_lr_GD_poly} for $\beta=0$ and with \eqref{eq:const_lr_HB_poly} for $0<\beta<1$,  and assume $\alpha<2(1+\beta)$. Consider the pair $\overline{L^{(\zeta)}_n}, L^{(\zeta)}_n$ of the worst-case loss \eqref{eq:worst_case_loss} and the loss  \eqref{eq:exact_power_measure_loss} for the exact power-law spectral measure. Then, as $n\to\infty$,
\begin{equation}\label{eq:exact_power_measure_loss_const_lr_GD}
    L^{(\zeta)}_n=\frac{\Gamma(\zeta+1)}{2}\Big(\frac{2\alpha n}{1-\beta}\Big)^{-\zeta}(1+o(1))
\end{equation}
and
\begin{equation}\label{eq:worst_case_loss_constant_lr_GD}
    \overline{L^{(\zeta)}_n} = L^{(\zeta)}_n + 
    \begin{cases}
    \mathbbm{1}_{\alpha>1} \;O\big(u^n\big), \quad &\beta=0\\
    \mathbbm{1}_{z_1<1} O\big(n^2u^n\big), \quad &\beta>0
    \end{cases}
\end{equation}
where $z_1=z(\lambda=1)$ and $u$ is some value independent of $n$ and such that $0<u<1$.
\end{theor}
Let's make a few remarks about Theorem \ref{ther:constant_lr_bounds}. First, observe from \eqref{eq:worst_case_loss_constant_lr_GD} that the worst-case loss $\overline{L^{(\zeta)}_n}$ is just equal to the power-law loss $L^{(\zeta)}_n$ if $\alpha\leq1$ or $z_1\geq1$. The reason behind this is that the flattened polynomial from Theorem \ref{ther:worst_case_loss} is unchanged: $\overline{p_n^2}(\lambda)=p_n^2(\lambda)$, which holds for $p_n(\lambda)$ monotone decreasing on $[0,1]$. For vanilla GD ($\beta=0$) the monotonicity can be seen directly from \eqref{eq:const_lr_GD_poly}, while for HB ($\beta\ne 0$) it requires more care but intuitively is connected to the localization of the roots of Chebyshev polynomials $U_n(z)$ on $[-1,1]$.  

Next, note that the difference between $\overline{L_n^{(\zeta)}}$ and $L_n^{(\zeta)}$ becomes exponentially small at large steps $n$. The speed of this exponential decay is given by parameter $u$, for which we obtain an explicit expression in the proof of the theorem. In particular, in GD without momentum $u=(1-\alpha)^2$, and we can clearly observe that the convergence condition $\alpha<2$ is equivalent to the condition $u<1$ of exponential decay of the correction term. 

Finally, we observe from \eqref{eq:exact_power_measure_loss_const_lr_GD} that the constant $C$ in the asymptotic $L^{(\zeta)}_n=Cn^{-\zeta}(1+o(1))$ can be made arbitrarily small by taking $\beta \nearrow 1$. In other words, higher ``inertia'' leads to faster convergence. We will see a reflection of this behavior in Section \ref{sec:ub}, where accelerated covergence rate $L^{(\zeta)}_n=O(n^{-2\zeta})$ is reached with the schedule of momentum behaving as $\beta_n\nearrow 1, \; n\to \infty$.

Now, we complete the picture for constant learning rate algorithms by establishing the lower bound in the class of discrete problems characterized by \eqref{eq:lkcknu}.
\begin{theor}\label{ther:const_lr_discrete_lower_bound}
Consider the discrete spectral measure $\rho_{\zeta,\nu} = \sum_{k=1}^{\infty} \big(k^{-\zeta \nu}-(k+1)^{-\zeta\nu}\big)\delta_{k^{-\nu}}$. Then 1) $\rho_{\zeta,\nu}$ satisfies both conditions \eqref{eq:rhozeta} and \eqref{eq:lkcknu}; 2) the loss of constant learning rate GD and HB with $\alpha<2(1+\beta)$ applied to the problem characterized by $\rho_{\zeta,\nu}$ is given by the right-hand side of \eqref{eq:exact_power_measure_loss_const_lr_GD}.
\end{theor}
Basically, this result indicates that constant learning rate algorithms can not take advantage of discrete power-law spectrum $\lambda_k\leq k^{-\nu}$. This fully settles $L(\mathbf w_n)=O(n^{-\zeta})$ as the tight bound for GD and HB in the case of constant learning rates. 

\subsection[A guide to acceleration: exact power-law spectral measure  ]{A guide to acceleration: exact power-law spectral measure (Section \ref{sec:lbcgjacobi})}\label{sec:exact_power-law_measure}
If we want to accelerate GD/HB, in the sense of decreasing the worst-case loss values $\overline{L_n^{(\zeta)}}$, Theorem \ref{ther:worst_case_loss} naturally guides us how to do that.
Specifically, assume that in search for an accelerated algorithm we  end up with a good enough family of polynomials $p_n(\lambda)$ such that $\overline{L_n^{(\zeta)}}$ and the exact power-law loss $L_n^{(\zeta)}$ are not far from each other, e.g. as in Theorem \ref{ther:constant_lr_bounds}. Then, instead of minimizing $\overline{L_n^{(\zeta)}}$ we can focus on minimizing $L_n^{(\zeta)}$. The latter problem is well-defined and is given by
\begin{equation}\label{eq:CG_polynomials}
    p_n = \argmin_{q_n:\deg q_n = n,q_n(0)=1}\frac{1}{2}\int_0^1 q_n^2(\lambda)\rho_\zeta(d\lambda).
\end{equation}
Recall (see Section \ref{sec:poly} and specifically Eq. \eqref{eq:CG_polynomials_opt}) that the solution to \eqref{eq:CG_polynomials} is exactly the CG algorithm applied to the exact power-law measure $\rho_\zeta$. The corresponding optimal residual polynomial can be found by expressing its variation as $\delta p_n(\lambda) = \lambda r_{n-1}(\lambda)$ with arbitrary degree-$(n-1)$ polynomial $r_{n-1}$ and then equating the variation of loss \eqref{eq:exact_power_measure_loss} to zero:
\begin{equation}\label{eq:exact_power_law_loss_variation}
    \delta L_n^{(\zeta)} = \int_0^1 p_n(\lambda) r_{n-1}(\lambda) \lambda\rho_\zeta( d\lambda)=0,
\end{equation}
implying that $p_n$ is an orthogonal polynomial on $[0,1]$ w.r.t. the weight $\lambda \rho_\zeta(d\lambda)=\zeta \lambda^\zeta d\lambda$. Then $p_n$ is a shifted and normalized Jacobi polynomial $P_n^{(\zeta,0)}(x)$:
\begin{equation}\label{eq:power-law_measure_solution}
    p_n(\lambda)=\frac{P_n^{(\zeta,0)}(1-2\lambda)}{P_n^{(\zeta,0)}(1)}.
\end{equation}

This leads to the precise convergence rate of Conjugate Gradients under the exact power-law measure $\rho_\zeta$.
\begin{theor}\label{th:cgexactpower}
The losses of CG method applied to a problem with measure \eqref{eq:exactpowerlaw2} are

\begin{equation}\label{eq:cgexactpower}
\begin{split}
    L(\mathbf w_n) ={}& \frac{\Gamma^2(\zeta+1)n!^2}{2\Gamma^2(\zeta+n+1)}=\frac{ \Gamma^2(\zeta+1)}{2}n^{-2\zeta}(1+o(1))\quad (n\to\infty).
\end{split}
\end{equation}
\end{theor}
This result implies, in particular, that under the main spectral assumption \eqref{eq:rhozeta1} the CG loss $L(\mathbf w_n)$  will not, in general, decrease faster than $O(n^{-2\zeta})$. The same is also true for GD and HB, since their losses at any iteration are not less than the respective loss of CG.  

The CG solution \eqref{eq:power-law_measure_solution} suggests that other residual polynomials based on Jacobi polynomials might be good candidates for an accelerated GD method under power-law spectral conditions. Moreover, the prospects of applying \eqref{eq:power-law_measure_solution} to practical problems require the robustness of the results with respect to errors in estimating the exponent $\zeta$. To address these questions, we consider a $3$-parameter \emph{ansatz} of residual polynomials
\begin{equation}\label{eq:Jacobi_ansatz}
    q^{(a,b,r)}_n(\lambda) = \frac{P^{(a,b)}_n(1-r\lambda)}{P^{(a,b)}_n(1)},
\end{equation}
which contains \eqref{eq:power-law_measure_solution} with parameters $(a,b,r)$ set to $a=\zeta, \; b=0, \; r=2$. Then, we have
\begin{prop}\label{prop:Jacobi_anzatz_typical_convergence}
Consider residual polynomials $p_n(\lambda)$ given by \eqref{eq:Jacobi_ansatz} with $a,b>-\tfrac{1}{2}$ and $r<2$. The respective exact power-law measure loss \eqref{eq:exact_power_measure_loss} is given by 
\begin{equation}\label{eq:Jacobi_ansatz_powerlaw_integral}
    L^{(\zeta)}_n= \begin{cases}
    \frac{\zeta\Gamma^2(a+1)B(\zeta,2a-2\zeta+1)}{2^\zeta r^{-\zeta}\Gamma^2(a-\zeta+1)}n^{-2\zeta}\big(1+o(1)\big),  & a>\zeta-\frac{1}{2}\\
    \frac{2^\zeta\zeta\Gamma^2(a+1)B(\tfrac{r}{2};\zeta-a-\tfrac{1}{2}, b+\tfrac{1}{2})}{2\pi r^\zeta} n^{-2a-1} \big(1+o(1)\big), & a<\zeta-\frac{1}{2}
    \end{cases}
\end{equation}
\end{prop}
Observe that Eqs. \eqref{eq:cgexactpower} and \eqref{eq:Jacobi_ansatz_powerlaw_integral} are consistent with each other. But most importantly, the condition $a>\zeta-\tfrac{1}{2}$ is critical to ensure the optimal convergence rate $O(n^{-2\zeta})$. Once this condition is ensured, the dependence on parameters $(a,b,r)$ becomes \emph{soft}: their variation only smoothly changes the constant without changing the rate $O(n^{-2\zeta})$. 

Let us make explicit the connection between ansatz \eqref{eq:Jacobi_ansatz} and the associated HB method with $n$-dependent learning rates $\alpha_n,\;\beta_n$. The connection is enabled by $q_n^{(a,b,r)}$ being obtained from rescaled and $n$-independently shifted family of orthogonal polynomials. This implies that the sequence $q_n^{(a,b,r)}$ obeys a 3-term recurrence relation, which, due to the residual normalization, has exactly the form of momentum update: $p_{n+1}=p_n-\alpha_n\lambda p_n + \beta_n(p_n-p_{n-1})$. The resulting learning rates for the ansatz \eqref{eq:Jacobi_ansatz} are given by
\begin{equation}\label{eq:Jacobi_ansatz_parameters}
    \begin{cases}
    \alpha_n = r\frac{(2n+a+b+1)(2n+a+b+2)}{2(n+a+1)(n+a+b+1)} = 2r + O(n^{-1}),\\
    \beta_n = \frac{n(n+b)(2n+a+b+2)}{(n+a+1)(n+a+b+1)(2n+a+b)} = 1- \frac{2a+1}{n}+O(n^{-2}). 
    \end{cases}
\end{equation}
Special cases of \eqref{eq:Jacobi_ansatz} and \eqref{eq:Jacobi_ansatz_parameters} were previously considered in \cite{brakhage1987ill} with parameters $a=\zeta-\tfrac{1}{2}, \; b=-\tfrac{1}{2},\; r=2$, and in \cite{hanke1991accelerated} with parameters $a=\zeta, \; b=0,\; r=2$. Our general formula \eqref{eq:Jacobi_ansatz_parameters} allows to give an example of parameters $(a,b)$ different from the cases considered by these authors and having a much simpler expression for learning rates. Specifically, with $a=b$ the Jacobi polynomials in \eqref{eq:Jacobi_ansatz} reduce to the ultraspherical polynomials $q^{(a,a,r)}_n=C_n^{a-\tfrac{1}{2}}(1-r\lambda)\big/C_n^{a-\tfrac{1}{2}}(1)$, and the respective learning rates are
\begin{equation}
    \begin{cases}
    \alpha_n = 2r-\frac{2a+1}{n+2a+1},\\
    \beta_n = 1-\frac{2a+1}{n+2a+1}.
    \end{cases}
\end{equation}
Our experiments (see Section \ref{sec:exp}) suggest that the $O(n^{-2\zeta})$ performance is retained even if we simplify the learning rate expressions even further, to the leading terms $\alpha_n=2r, \; \beta_n=1-\frac{2a+1}{n}$ in Eq. \eqref{eq:Jacobi_ansatz_parameters}, but we do not have a proof of optimality in this case.

Importantly, \cite{brakhage1987ill} and \cite{hanke1991accelerated} used relatively indirect reasoning to arrive at their accelerated methods based on Jacobi polynomials. In contrast, our approach is straightforward -- given a spectral condition $\rho([0,\lambda])\leq G(\lambda)$, one simply needs to take the system of polynomials orthogonal w.r.t. the weight function $\lambda G'(\lambda)$. In particular, we expect that our approach can be  generalized to spectral conditions specified by functions $G(\lambda)$ other than power-laws.

\subsection{General upper bounds}\label{sec:ub}
\paragraph{Jacobi ansatz (Section \ref{sec:ubcggeneral}).}
The key intuition employed in the previous sections was that the GD method efficiently minimizing $L_n^{(\zeta)}$ would also work for all problems specified by \eqref{eq:rhozeta1}. We quantify this intuition in the following way:  
\begin{theor}\label{ther:Jacobi_worst_case_UB}
    Consider residual polynomials $p_n(\lambda)$ given by \eqref{eq:Jacobi_ansatz} with $a,b>-\tfrac{1}{2}$ and $r\leq\tfrac{2a+1}{a+b+1}$. Then, the worst-case loss \eqref{eq:worst_case_loss} is bounded in terms of exact power-law loss \eqref{eq:exact_power_measure_loss}:
    \begin{equation}\label{eq:Jacobi_ansatz_UB}
        \overline{L_n^{(\zeta)}} \leq C_\zeta L_n^{(\zeta)},
    \end{equation}
    where the constant $C_\zeta= C[\rho_\zeta]$ is given by the functional $C[\rho]$ of measure $\rho$ defined as
    \begin{equation}\label{eq:parabola_av_infinum}
    \frac{1}{C[\rho]} = \underset{c,x_l,x_r \in \supp \rho}{\operatorname{inf}} \left[\int_{x_l}^{x_r} \frac{(\lambda-c)^2}{\operatorname{max}\big((x_l-c)^2,(x_r-c)^2\big)} \rho(d\lambda)\bigg/\int_{x_l}^{x_r} \rho(d\lambda)\right]
\end{equation}
\end{theor}
The functional $C[\rho]$ has a simple geometric interpretation: the expression minimized in \eqref{eq:parabola_av_infinum} is a $\rho$-weighted average of a parabola with center at $c$ and normalized by its value at one of the edges $x_l, x_r$. The origin of this parabola is illustrated in Figure \ref{fig:Flat_polynomial} (right): if a flat region of polynomial $\overline{p_n^2}(\lambda)$ contains only a single root, the true polynomials $p_n^2(\lambda)$ can be lower-bounded by a such normalized parabola. Looking at the contribution to the losses \eqref{eq:worst_case_loss},\eqref{eq:exact_power_measure_loss} from this flat region reveals that their ratio is not worse than the ratio of the $\rho$-weighted averages of the constant and the normalized parabola over the flat region. Interestingly, the geometric picture depicted on Figure \ref{fig:Flat_polynomial} (right) requires only basic properties of polynomials $p_n(\lambda)$: non-degeneracy of the roots and monotonicity of local extrema.  We explicitly calculate the functional $C[\rho]$ for the exact power-law measure.

\begin{prop}\label{prop:parabola_powerlaw_av_minimum}
Let $\rho_\zeta$ be defined as in \eqref{eq:exactpowerlaw2}. Then
\begin{equation}\label{eq:powerlaw_parabola_infinum}
    C[\rho_\zeta]=\begin{cases}
    (\zeta+1)^{2}, \quad &\zeta\ge 1\\
    2+2/\zeta, \quad & \zeta \leq 1
    \end{cases}
\end{equation}
\end{prop}
Now, we denote the coefficient in the $a>\zeta-\tfrac{1}{2}$ case of Eq. \eqref{eq:Jacobi_ansatz_powerlaw_integral} by $R(a,r,\zeta)$, and summarize Theorem \ref{ther:Jacobi_worst_case_UB} and Propositions \ref{prop:Jacobi_anzatz_typical_convergence}, \ref{prop:parabola_powerlaw_av_minimum} as
\begin{corol}\label{corol:HB_upper_bound}
Let $a>\zeta-\tfrac{1}{2}, \;b>-\tfrac{1}{2}, \; r \leq \tfrac{2a+1}{a+b+1}$. Then the loss of HB method with the schedule \eqref{eq:Jacobi_ansatz_parameters} applied to a problem described by condition \eqref{eq:rhozeta1} is
\begin{equation}\label{eq:HB_upper_bound}
    L(\mathbf{w}_n) \leq  C_\zeta R(a,r,\zeta) n^{-2\zeta}(1+o(1))
\end{equation}
\end{corol}
The same bound obviously remains valid for CG, since its loss is dominated by the HB loss. 

\paragraph{GD with predefined schedule (Section \ref{sec:ubgdpredefined}).} The above result ensures an $O(n^{-2\zeta})$ convergence of HB with a suitable problem-independent learning rate schedule. We show that, theoretically, such a rate can also be achieved for GD (i.e., without using momentum): 
\begin{theor}\label{th:ubgdnonconst} 
Given $\zeta>0$, there exists a sequence $\alpha_n$ such that for any problem subject to spectral condition \eqref{eq:rhozeta1}, GD with this schedule $\alpha_n$ satisfies 
\begin{equation}
    L(\mathbf w_n)\le   C_\zeta R(a,r,\zeta) 4^{2\zeta}n^{-2\zeta} \big(1+o(1)\big),
\end{equation}
where the parameter $a$ and the constants $C_\zeta, R(a,r,\zeta)$ are as in Corollary \ref{corol:HB_upper_bound}.
\end{theor}
The idea of the proof is to consider a subsequence of polynomials \eqref{eq:Jacobi_ansatz} with growing degrees $2^l$, and choose the learning rates $\alpha_n$ as inverse roots of these polynomials. 

We remark, however, that this construction requires very large learning rates $\alpha_n$, which makes this algorithm, in contrast to HB with schedule \eqref{eq:Jacobi_ansatz_parameters}, fairly unstable for non-linear models (see experiments in Section \ref{sec:exp}).

\paragraph{Conjugate Gradients: discrete spectrum (Section \ref{sec:ubcgdiscrete}).} If the main spectral condition \eqref{eq:rhozeta1} is supplemented by eigenvalue decay condition \eqref{eq:lkcknu}, CG acquires quite different convergence rate $O(n^{-(2+\nu)\zeta})$:  
\begin{theor}\label{th:cgdiscr}
Assuming spectral conditions \eqref{eq:rhozeta1} and \eqref{eq:lkcknu}, the losses of CG  satisfy 
\begin{equation}\label{eq:cgdisclwn}
     L(\mathbf w_n) \leq  C_\zeta R(a,r,\zeta) \Lambda^\zeta (n/2)^{-(2+\nu)\zeta}(1+o(1)),
\end{equation}
where the parameter $a$ and the constants $C_\zeta, R(a,r,\zeta)$ are as in Corollary \ref{corol:HB_upper_bound}.
\end{theor}
The proof is based on assigning half of the roots of trial polynomials $q_n$ to the largest atoms $\lambda_k$ of the spectral measure $\rho$, and then adjusting the remaining roots on the segment $[0,\Lambda (n/2)^{-\nu}]$ by rescaling and invoking Corollary \ref{corol:HB_upper_bound}.

\subsection{Further lower bounds}\label{sec:lb}
\paragraph{Non-adaptive schedules (Section \ref{sec:lbpredefined}).}
If an optimization algorithm has a predefined (non-adaptive) learning rate schedule (as in our GD or HB), then it cannot in general  improve the exponent $2\zeta$ in the convergence rate $O(n^{-2\zeta})$, even if we additionally assume the discreteness of the spectrum with a particular power law decay: 

\begin{theor}\label{ther:discrete_spectra_predefined_schedule_LB}
Consider any optimization algorithm of the form 
\begin{equation}\label{eq:multistep}
    \mathbf w_{n+1} = \mathbf w_0+\sum_{j=0}^n \alpha_{nj}\nabla L(\mathbf w_j)
\end{equation}  
with fixed (problem-independent) $\alpha_{nj}.$ Then for any  $\zeta,\nu,\epsilon>0$ there exists a problem with a compact $A$ and $\mathbf b$ subject to  
\begin{align}
    \lambda_n ={}& n^{-\nu}(1+o(1)),\quad n\to\infty,\label{eq:lndisc1}\\
    \rho((0,\lambda]) ={}& \lambda^\zeta(1+o(1)),\quad \lambda\to 0+,\label{eq:lndisc2}
\end{align}
such that there is an infinite sequence $n_1<n_2<\ldots$ for which
\begin{equation}
    L(\mathbf w_{n_s})>n_s^{-2\zeta-\epsilon}.
\end{equation}
\end{theor}

\paragraph{CG with discrete spectrum (Section \ref{sec:lbcgdiscrete}).}
We give an explicit example showing that the bound $L(\mathbf w_n)=O(n^{-(2+\nu)\zeta})$ established in Theorem \ref{th:cgdiscr} for CG under two spectral conditions \eqref{eq:rhozeta1}, \eqref{eq:lkcknu} cannot generally be improved.  
For any constants $\nu>0$ and $\zeta>0$, consider the operator $J$ defined on the space $l^2$ of square-summable sequences $\mathbf w=(w_1,w_2,\ldots)$  by 
\begin{equation}\label{eq:j}
    (J\mathbf w)_n=\begin{cases}
    w_1,& n=1,\\
    n^{-\frac{\nu}{2}}w_n-(\tfrac{n}{n-1})^{\frac{1-(2+\nu)\zeta}{2}}(n-1)^{-\frac{\nu}{2}}w_{n-1},& n=2,3,\ldots
    \end{cases}
\end{equation}
Next, let $\mathbf f_*=\mathbf e_1=(1,0,\ldots)$. We will show that the quadratic problem \eqref{eq:lwj} defined by these $J$ and $\mathbf f_*$ is a desired example. 

Let us clarify the idea behind this choice of the operator $J$. Its two-diagonal form implies that the respective Krylov subspaces are just the standard coordinate subspaces, which allows to easily compute the exact loss trajectory $L(\mathbf w_n)$ (statement 1 of the following theorem). On the other hand, the coefficients in Eq. \eqref{eq:j} are adjusted to ensure the desired asymptotics of the eigenvalues $\lambda_k$ and the spectral measure $\rho$ (statements 2 and 3).

\begin{theor}\label{th:cglbdiscr}\hfill
\begin{enumerate}
\item The loss values of CG for the problem defined by the above $J$ and $\mathbf f_*$ are 
    \begin{align}
        L(\mathbf w_n)=\Big(2\sum_{m=1}^{n+1}m^{(2+\nu)\zeta-1}\Big)^{-1}=(1+o(1))\frac{(2+\nu)\zeta}{2}n^{-(2+\nu)\zeta},\quad n\to\infty.\nonumber
    \end{align} 
    \item For any $\nu>0$ and $\zeta>0$, $\widetilde A=JJ^\dagger$ is a compact  operator with eigenvalues $\lambda_k=O(k^{-\nu})$.
    \item For any non-integer $\zeta>0$, the spectral measure $\rho$ associated with $\widetilde A$ and $\mathbf f_*$ satisfies $\rho((0,\lambda])=O(\lambda^\zeta)$ as $\lambda\to 0+$.
\end{enumerate}
\end{theor}
The restriction to non-integer $\zeta$ in Statement 3 is due to our proof technique; it can probably be lifted using a more careful analysis. If $\zeta$ is non-integer, then Theorem \ref{th:cglbdiscr} gives precisely an example of a problem satisfying spectral conditions \eqref{eq:rhozeta1}, \eqref{eq:lkcknu} and a lower bound $L(\mathbf w_n)=\Omega( n^{-(2+\nu)\zeta})$. If $\zeta$ is an integer, then we can still use the theorem for a slightly weaker conclusion: considering operator \eqref{eq:j}  with $\zeta$ replaced by $\zeta+\epsilon$ with an arbitrary $0<\epsilon<1$, we get an example satisfying spectral conditions \eqref{eq:rhozeta1}, \eqref{eq:lkcknu} and a lower bound $L(\mathbf w_n)=\Omega(n^{-(2+\nu)(\zeta+\epsilon)})$.  

\smallskip
\begin{proof}[of Theorem \ref{th:cglbdiscr}]
As a preliminary observation, note that $J^\dagger$ is given by
\begin{equation}\label{eq:jstar}
    (J^\dagger\mathbf x)_n = n^{-\nu/2}(x_n-(\tfrac{n+1}{n})^{(1-(2+\nu)\zeta)/2}x_{n+1}),\quad n=1,2,\ldots
\end{equation}
\emph{Statement 1.} In the case of CG, $L(\mathbf w_n)$ is obtained by optimizing $L(\mathbf w)$ over the Krylov subspace spanned by $\{(J^\dagger J)^mJ^\dagger\mathbf e_1\}_{m=0}^{n-1}.$ Note that $(J^\dagger J)^mJ^\dagger=J^\dagger(JJ^\dagger)^m=J^\dagger\widetilde A^m$ and that $\widetilde A$ is three-diagonal, so that the vectors $\{\widetilde A^m\mathbf e_1\}_{m=0}^{n-1}$ span the coordinate subspace $\mathcal H_n$ spanned by $\mathbf e_1,\ldots,\mathbf e_n$. Therefore, 
\begin{equation}
    L(\mathbf w_n)=\min_{\mathbf x\in\mathcal H_n}\tfrac{1}{2}\|JJ^\dagger\mathbf x-\mathbf e_1\|^2=\min_{\mathbf x\in\mathcal H_n}\tfrac{1}{2}\|\widetilde A\mathbf x-\mathbf e_1\|^2.
\end{equation}
Consider the vector $\mathbf v=(v_1,v_2,\ldots)$ defined by
\begin{equation}
    v_m=\begin{cases}
    m^{((2+\nu)\zeta-1)/2},& m=1,\ldots,n+1,\\
    0,& m > n+1.
    \end{cases}
\end{equation}
Then, using Eq. \eqref{eq:jstar}, $(J^\dagger\mathbf v)_m=0$ for $m=1,\ldots,n.$ Accordingly, $\langle \widetilde A\mathbf x,\mathbf v\rangle=\langle J^\dagger\mathbf x,J^\dagger\mathbf v\rangle=0$ for any $\mathbf x\in\mathcal H_n.$ On the other hand, it is easy to see that if a vector $\mathbf u$ in the coordinate subspace $\mathcal H_{n+1}$ is orthogonal to this $\mathbf v$, then $\mathbf u=\widetilde A\mathbf x$ for some $\mathbf x\in \mathcal H_n.$ It follows that
\begin{align}
    L(\mathbf w_n)=\min_{\mathbf x\in\mathcal H_{n+1}\ominus\mathbf v}\tfrac{1}{2}\|\mathbf x-\mathbf e_1\|^2=\frac{\langle \mathbf e_1,\mathbf v\rangle^2}{2\|\mathbf v\|^2}=\Big(2\sum_{m=1}^{n+1}m^{(2+\nu)\zeta-1}\Big)^{-1},
\end{align}
as desired.

\medskip
\noindent 
\emph{Statement 2} is implied by the following (more detailed) characterization of the spectrum of~ $\widetilde A$. 
\begin{lem}
The operator $\widetilde A$ is compact, and the sorted positive eigenvalues $\lambda_1\ge \lambda_2\ge \ldots > 0$ satisfy 
\begin{equation}
    (2k)^{-\nu}\le \lambda_k \le 5k^{-\nu}.
\end{equation}
\end{lem}
\begin{proof}
The compactness follows since $J$ is approximated in norm by the finite-dimensional operators obtained by truncating the assignment \eqref{eq:j}. As a result of compactness, the spectrum of $\widetilde A$ is discrete and consists of nonnegative eigenvalues; the positive eigenvalues can be sorted in decreasing order. To lower bound the eigenvalues, use the minimax principle:
\begin{equation}
    \lambda_k = \max_{\substack{\mathcal H_k\subset l^2:\\ \dim \mathcal H_k=k}}\min_{\mathbf x\in\mathcal H_k}\frac{\langle \widetilde A\mathbf x,\mathbf x\rangle}{\|\mathbf x\|^2}.
\end{equation}
Choosing the subspace $\mathcal H_k$ spanned by $\mathbf e_2,\mathbf e_4,\ldots,\mathbf e_{2k}$, we get
\begin{align}
    \lambda_k \ge{}&\min_{\mathbf x\in\mathcal H_k}\frac{\langle \widetilde A\mathbf x,\mathbf x\rangle}{\|\mathbf x\|^2}=\min_{\mathbf x\in\mathcal H_k}\frac{\|J^\dagger\mathbf x\|^2}{\|\mathbf x\|^2}\\
    ={}&\min_{\mathbf x\in\mathcal H_k}\frac{1}{\|\mathbf x\|^2}\sum_{m=1}^{k}[(\tfrac{2m}{2m-1})^{1-(2+\nu)\zeta} (2m-1)^{-\nu}+(2m)^{-\nu}]x_{2m}^2\\
    \ge{}& (2k)^{-\nu}.
\end{align}
To upper bound $\lambda_k$ use the minimax principle in a different form:
\begin{equation}
    \lambda_k = \min_{\substack{\mathcal G_k\in l^2:\\ \dim (l^2\ominus\mathcal G_k)=k-1}}\max_{\mathbf x\in\mathcal G_k}\frac{\langle \widetilde A\mathbf x,\mathbf x\rangle}{\|\mathbf x\|^2}.
\end{equation}
Choosing $\mathcal G_k$ spanned by $\mathbf e_k,\mathbf e_{k+1},\ldots$, we get
\begin{align}
    \lambda_k \le{}&\max_{\mathbf x\in\mathcal G_k}\frac{\langle \widetilde A\mathbf x,\mathbf x\rangle}{\|\mathbf x\|^2}=\max_{\mathbf w\in\mathcal G_k}\frac{\|J^\dagger\mathbf x\|^2}{\|\mathbf x\|^2}\\
    ={}&\max_{\mathbf x\in\mathcal G_k}\frac{1}{\|\mathbf x\|^2}\Big[k^{-\nu}x_k^2+\sum_{m=k+1}^{\infty}(m^{-\nu/2}x_m-(\tfrac{m}{m-1})^{(1-(2+\nu)\zeta)/2}(m-1)^{-\nu/2}x_{m-1})^2 \Big]\nonumber\\
    \le{}&\max_{\mathbf x\in\mathcal G_k}\frac{1}{\|\mathbf x\|^2}\Big[k^{-\nu}x_k^2+\sum_{m=k+1}^{\infty}\big(2m^{-\nu}x_m^2+4(m-1)^{-\nu}x_{m-1}^2\big) \Big]\\
    \le{}& \max_{\mathbf x\in\mathcal G_k}\frac{1}{\|\mathbf x\|^2}\Big[5k^{-\nu}\sum_{m=k}^\infty x_m^2\Big]\\
    ={}& 5k^{-\nu}.
\end{align}
\end{proof}

\noindent
\emph{Statement 3} relies on the following resolvent bounds.

\begin{prop}\label{lm:resolv2} \hfill
\begin{enumerate}
\item Assuming $2m<\zeta<2m+1$ for some integer $m\ge 0$, the vectors $\widetilde A^{-m}\mathbf e_1$ and $\widetilde A^{-m}(\widetilde A+\epsilon)^{-1}\mathbf e_1$ exist as elements of $l^2$ and 
\begin{equation}\label{eq:resolv21} 
    \langle \widetilde A^{-m}\mathbf e_1, \widetilde A^{-m}(\widetilde A+\epsilon)^{-1}\mathbf e_1\rangle = O(\epsilon^{\zeta-2m-1}),\quad \epsilon\to 0+.
\end{equation}
\item Assuming $2m+1<\zeta<2m+2$ for some integer $m\ge 0$, the vectors $J^{-1}\widetilde A^{-m}\mathbf e_1$ and $J^{-1}\widetilde A^{-m}(\widetilde A+\epsilon)^{-1}\mathbf e_1$ exist as elements of $l^2$ and 
\begin{equation} \label{eq:resolv22} 
    \langle J^{-1}\widetilde A^{-m}\mathbf e_1, J^{-1}\widetilde A^{-m}(\widetilde A+\epsilon)^{-1}\mathbf e_1\rangle = O(\epsilon^{\zeta-2m-2}),\quad \epsilon\to 0+.
\end{equation}
\end{enumerate}
\end{prop}

\noindent
The proof of this proposition is quite lengthy, and we defer it to Sections \ref{sec:lmresolv} and \ref{sec:lmresolv2}. Let us show how it implies the desired spectral bound.

Assume first that $2m<\zeta<2m+1$ for some integer $m\ge 0$. By definition of the spectral measure,
\begin{align}
    \langle \widetilde A^{-m}\mathbf e_1, \widetilde A^{-m}(\widetilde A+\epsilon)^{-1}\mathbf e_1\rangle ={}& \int_0^\infty \lambda^{-m}\cdot\lambda^{-m}(\lambda+\epsilon)^{-1} \rho(d\lambda)\\
    \ge{}&\int_{0}^\epsilon \epsilon^{-2m}(2\epsilon)^{-1}\rho(d\lambda)\\
    ={}&\tfrac{1}{2}\epsilon^{-1-2m}\rho((0,\epsilon]).
\end{align}
It follows then by Statement 1 of Proposition \ref{lm:resolv2} that
\begin{equation}
    \rho((0,\lambda]) \le 2\epsilon^{1+2m} O(\epsilon^{\zeta-2m-1})=O(\epsilon^\zeta),
\end{equation}
as desired.

The case $2m+1<\zeta<2m+2$ is analyzed similarly, using part 2 of the proposition and the observation 
\begin{align}
    \langle J^{-1}\widetilde A^{-m}\mathbf e_1, J^{-1}\widetilde A^{-m}(\widetilde A+\epsilon)^{-1}\mathbf e_1\rangle ={}& \int_0^\infty \lambda^{-m-1/2}\cdot\lambda^{-m-1/2}(\lambda+\epsilon)^{-1} \rho(d\lambda).
\end{align}

\end{proof}

\begin{figure}
    \centering
    {\includegraphics[scale=0.75, trim=6mm 4mm 2mm 9mm, clip]{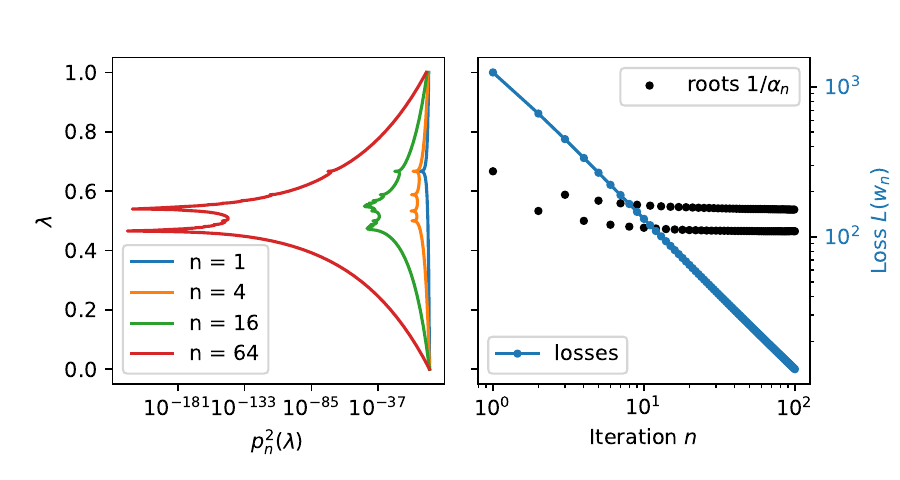}}
    
    \vspace{-1mm}
    \caption{SD applied to the uniform spectral distribution ($\rho((0,\lambda])=Q\lambda$) on $[0,1]$ converges to a period-2 oscillatory regime.
    }
    \label{fig:SD}
\end{figure}

\subsection{Steepest descent}\label{sec:sd}
Our analysis of SD is based on the remarkable asymptotic periodicity of this algorithm: as $n$ increases, the adaptive learning rates $\alpha_n$ start to perform approximate period-2 oscillations, and the subsequences $\alpha_{2n}$ and $ \alpha_{2n+1}$ converge (see Figure \ref{fig:SD}). This effect was first established, for finite-dimensional problems, in \cite{akaike1959successive}. We will use a generalization to infinite-dimensional spaces proved in \cite{pronzato2001renormalised}.   

Denote by $\lambda_{\min}$ and $\lambda_{\max}$ the left and right ends of the support of spectral measure~$\rho$:
\begin{equation}
    \lambda_{\min} = \sup\{\lambda: \rho((-\infty,\lambda))=0\},\qquad \lambda_{\max} = \inf\{\lambda: \rho((\lambda,\infty))=0\}.
\end{equation}
We will assume that $\lambda_{\min}\ne \lambda_{\max}$ (excluding the trivial case of a Dirac delta), so 
\begin{equation}
    0\le \lambda_{\min}<\lambda_{\max}<\infty.
\end{equation}
It is convenient to introduce the inverses $b_n$ of the learning rates $\alpha_n$:
\begin{equation}
    b_n=1/\alpha_n.
\end{equation}
The values $b_n$ are the roots of the residual polynomials $p_n$ associated with the iterates of SD (see Section \ref{sec:poly}):
\begin{equation}\label{eq:sdqn}
    p_n(\lambda) = \prod_{s=0}^{n-1}(1-\lambda/b_s).
\end{equation}
By definition of SD, $\alpha_n$ is obtained by optimizing 
\begin{equation}
    \int_{\lambda_{\min}}^{\lambda_{\max}} (1-\alpha\lambda)^2p^2_{n}(\lambda)\rho(d\lambda)\to \min_{\alpha}.
\end{equation}
This gives
\begin{equation}\label{eq:anlambdaq1}
    \alpha_n=\frac{\int_{\lambda_{\min}}^{\lambda_{\max}}  \lambda p_n^2(\lambda)\rho(d\lambda)}{\int_{\lambda_{\min}}^{\lambda_{\max}}  \lambda^2 p_n^2(\lambda)\rho(d\lambda)},\qquad 
    b_n=\frac{\int_{\lambda_{\min}}^{\lambda_{\max}}  \lambda^2 p_n^2(\lambda)\rho(d\lambda)}{\int_{\lambda_{\min}}^{\lambda_{\max}}  \lambda p_n^2(\lambda)\rho(d\lambda)}.
\end{equation}
Let us introduce the probability measure $\sigma_n$ by
\begin{equation}
    \sigma_n(d\lambda)=Z_n^{-1}\lambda p_n^2(\lambda)\rho(d\lambda),
\end{equation}
where $Z_n=\int_{\lambda_{\min}}^{\lambda_{\max}} \lambda p_n^2(\lambda)\rho(d\lambda)$ is the normalizing factor.   Eq. \eqref{eq:anlambdaq1} shows that $b_n$ is the mean of $\sigma_n$: 
\begin{equation}
    b_{n}=\int_{\lambda_{\min}}^{\lambda_{\max}}\lambda\sigma_n(d\lambda).
\end{equation}
Moreover, using Eq. \eqref{eq:sdqn}, the evolution of the measures $\sigma_n$ with SD iterations is given simply by
\begin{align}\label{eq:sigmait}
    \sigma_{n+1}(d\lambda)=D_n^{-1}(\lambda-b_n)^2\sigma_n(d\lambda),
\end{align}
where $D_n=\int_{\lambda_{\min}}^{\lambda_{\max}} (\lambda-b_n)^2 \sigma_n(d\lambda)$ is the variance of $\sigma_n$. 

By our assumptions, $0$ is not an eigenvalue of $\widetilde A$ and so is not an isolated atom of the measure $\rho$. It follows that the measure $\sigma_0(d\lambda)=Z_0^{-1}\lambda\rho(d\lambda)$ has the same end points $\lambda_{\min}, \lambda_{\max}$ of its support as the measure $\rho$. 

Evolution \eqref{eq:sigmait} admits a simple family of special period-2 solutions parameterized by $q\in(0,1)$: 
\begin{equation}
    \sigma_{2n}=q\delta_{\lambda_{\min}}+(1-q)\delta_{\lambda_{\max}},\qquad \sigma_{2n+1}=(1-q)\delta_{\lambda_{\min}}+q\delta_{\lambda_{\max}}.
\end{equation}
The following result shows that any sequence of iterates $\sigma_n$ is attracted to one of these special solutions.

\begin{theor}[Theorem 2 in \cite{pronzato2001renormalised}]\label{th:period2}
Consider iterations \eqref{eq:sigmait} starting from some compactly supported Borel probability measure $\sigma_0$ with end points $\lambda_{\min}< \lambda_{\max}$ of its support.\footnote{The statement of this theorem in \cite{pronzato2001renormalised} also includes the condition $\lambda_{\min}>0$, but it is clear that this condition can be dropped since evolution \eqref{eq:sigmait} is translation invariant.} Then there exists  $q\in (0,1)$ such that for any $\lambda\in(\lambda_{\min},\lambda_{\max})$ 
\begin{equation}
    \sigma_{2n}([\lambda_{\min},\lambda])\stackrel{n\to\infty}{\longrightarrow}q,\qquad \sigma_{2n+1}([\lambda_{\min},\lambda])\stackrel{n\to\infty}{\longrightarrow}1-q.
\end{equation}
\end{theor}
This result implies, in particular, that
\begin{equation}\label{eq:anlim}
    b_{2n}\stackrel{n\to\infty}{\longrightarrow} q\lambda_{\min}+(1-q)\lambda_{\max},\qquad b_{2n+1}\stackrel{n\to\infty}{\longrightarrow}(1-q)\lambda_{\min}+q\lambda_{\max}.
\end{equation}
Using Theorem \ref{th:period2} and asymptotics \eqref{eq:anlim}, it is easy to connect the convergence rates of the SD evolution to those of GD with constant rates. The case $\lambda_{\min}>0$ is discussed in Section 5 of \cite{pronzato2001renormalised}; it is shown there that in this case the convergence of SD is (like that of GD) exponentially fast:
\begin{equation}
    L(\mathbf w_n)=O\Big(\Big(\frac{\lambda_{\max}-\lambda_{\min}}{\lambda_{\max}+\lambda_{\min}}+\epsilon\Big)^{2n}\Big),\quad n\to\infty,
\end{equation}
for any $\epsilon>0$. Consider now the case $\lambda_{\min}=0$. The loss $L(\mathbf w_n)$ can be written in terms of $\sigma_n$ as 
\begin{equation}
    L(\mathbf w_n)=\frac{Z_n}{2}\int_{0}^{\lambda_{\max}}\lambda^{-1}\sigma_n(d\lambda).
\end{equation}
Applying Theorem \ref{th:period2}, the leading contribution to this integral comes from small neighborhoods of $\lambda=0$: for any $\widetilde \lambda\in(0,\lambda_{\max})$ 
\begin{equation}
    \int_{0}^{\lambda_{\max}}\lambda^{-1}\sigma_n(d\lambda)=(1+o(1))\int_{0}^{\widetilde\lambda}\lambda^{-1}\sigma_n(d\lambda),\quad n\to\infty,
\end{equation}
and accordingly
\begin{equation}
    L(\mathbf w_n)=(1+o(1))\frac{1}{2}\int_{0}^{\widetilde \lambda}p_n^2(\lambda)\rho(d\lambda),\quad n\to\infty.
\end{equation}
Now choose $\widetilde \lambda=\tfrac{1}{2}\inf_n b_n$. Using convergence \eqref{eq:anlim} of the values $b_n$, we have $\widetilde \lambda>0$. Recalling that the values $b_n$ are the roots of the residual polynomials $p_n$, we can find constants $c_1,c_2>0$ such that for any  $\lambda\in[0,\widetilde \lambda]$ and $n$
\begin{equation}
    e^{-nc_1\lambda}\le p_n^2(\lambda)\le e^{-nc_2\lambda}
\end{equation}
and so
\begin{equation}\label{eq:sdlnsand}
     (1+o(1))\frac{1}{2}\int_0^{\widetilde \lambda}e^{-nc_1\lambda}\rho(d\lambda)\le L(\mathbf w_n)\le (1+o(1))\frac{1}{2}\int_0^{\widetilde \lambda}e^{-nc_2\lambda}\rho(d\lambda),\quad n\to\infty.
\end{equation}
Integrating by parts and making the change of variable $nc \lambda=t$,
\begin{equation}
  \int_0^{\widetilde \lambda}e^{-nc\lambda}\rho(d\lambda)=e^{-nc\widetilde \lambda}\rho((0,\widetilde \lambda])+\int_0^{nc\widetilde \lambda}e^{-t}\rho((0,\tfrac{t}{cn}])dt.  
\end{equation}
The first term falls off exponentially, while in the case of the power law measure $\rho((0,\lambda])=\min(\lambda^{\zeta},\lambda_{\max}^{\zeta})$ the second term equals $\Gamma(\zeta+1)(cn)^{-\zeta}(1+o(1))$. Combined with Eq. \eqref{eq:sdlnsand}, this immediately implies the desired upper and lower loss bounds:
\begin{theor}\label{th:sd}
    Assuming the main spectral condition \eqref{eq:rhozeta1}, the SD loss obeys $L(\mathbf w_n)=O(n^{-\zeta})$. On the other hand, if we assume a lower bound $\rho((0,\lambda])=\Omega(\lambda^{\zeta})$, then $L(\mathbf w_n)=\Omega(n^{-\zeta})$.
\end{theor}
Recall the discrete measure $\rho_{\zeta,\nu} = \sum_{k=1}^{\infty} (k^{-\zeta \nu}-(k+1)^{-\zeta\nu})\delta_{k^{-\nu}}$ that appeared in Theorem \ref{ther:const_lr_discrete_lower_bound} and satisfies both main spectral condition  \eqref{eq:rhozeta1} and eigenvalue decay condition \eqref{eq:lkcknu}. It is easy to see that $\rho((0,\lambda])\ge 2^{-\zeta\nu}\lambda^{-\zeta}$ for $0<\lambda\le 1$, so both statements of Theorem \ref{th:sd} are applicable to $\rho_{\zeta,\nu}$. It follows that the loss convergence bound $O(n^{-\zeta})$ is tight even if the main spectral condition \eqref{eq:rhozeta1} is supplemented by the eigenvalue decay condition \eqref{eq:lkcknu}.

We remark that a $O(n^{-\zeta})$ upper bound for the loss was obtained previously by a different method, based on moment inequalities, in \cite{gilyazov2013regularization} (see their Theorem 2.2.5). However, that method seems to require the stronger source condition \eqref{eq:sourcecond} and does not produce tight lower bounds.

\section{Comparison of spectral conditions}\label{sec:tighterub}
As discussed in Section \ref{sec:probdef}, our target expansion condition \eqref{eq:rhozeta} is a variant of the more standard source condition \eqref{eq:sourcecond1}. In this section we compare the two versions and argue that our condition \eqref{eq:rhozeta} can be more convenient and natural in applications. We have already shown in Lemma \ref{lem:conditions_comp} and Section \ref{sec:detail} that our condition \eqref{eq:rhozeta} with a particular exponent $\zeta$ is slightly weaker than the respective source condition \eqref{eq:sourcecond1}, but leads to similar power-law loss bounds $O(n^{-\zeta}), O(n^{-2\zeta}), O(n^{-(2+\nu)\zeta})$. We will argue now that, moreover, our condition generally better fits practical power-law spectra and produces tighter bounds when optimized over spectral parameters.   

\paragraph{Upper bounds for classical source condition.}
We briefly recap the classical technique used for obtaining loss upper bounds under the classical source condition \eqref{eq:sourcecond1} (see, e.g. \cite{Polyak87, nemirovskiy1984iterative1, brakhage1987ill}). Recall that the loss is given by $L_n=\tfrac{1}{2}\int_0^1p_n^2(\lambda)\rho(d\lambda)$ with a residual polynomial $p_n$ associated with a particular optimization algorithm. Consider $p_n$ as fixed and the loss $L_n=L_n(\rho)$ as a function of measure $\rho$. Under the classical source condition with parameters $\zeta',Q'$, the largest value of $L_n$ is
\begin{equation}\label{eq:source_cond_worst_case}
    \sup_{\rho \in \mathrm{P}'(\zeta',Q')}L_n(\rho)=\sup_{\rho \in \mathrm{P}'(\zeta',Q')}\frac{1}{2}\int_0^1 [\lambda^{\zeta'} p_n^2(\lambda)]\lambda^{-\zeta'}\rho(d\lambda)=\frac{Q'}{2} \underset{0\leq \lambda \leq 1}{\operatorname{sup}} \big[\lambda^{\zeta'} p_n^2(\lambda)\big].
\end{equation}
The value $\omega(\zeta',p_n)\equiv\sup_{0\le\lambda\le 1} \big[\lambda^{\zeta'} p_n^2(\lambda)\big]$ is the main object studied in \cite{Polyak87, nemirovskiy1984iterative1, brakhage1987ill} and other related works to characterize convergence rates. Note that the loss in \eqref{eq:source_cond_worst_case} is maximized at the rescaled Dirac delta $\rho^*=Q'(\lambda^*)^{\zeta'}\delta_{\lambda^*}$, where $\lambda^* = \underset{0\leq \lambda \leq 1}{\argmax} [\lambda^{\zeta'}p_n^2(\lambda)]$. This shows that the tightest upper bound under the source condition is
\begin{align}\label{eq:lnubcl}
    L^{UB}_n(\zeta',Q')=\sup_{\rho \in \mathrm{P}'(\zeta',Q')}L_n(\rho)=\frac{Q'}{2} \underset{0\leq \lambda \leq 1}{\operatorname{sup}} \big[\lambda^{\zeta'} p_n^2(\lambda)\big],
\end{align}
and the bound is especially accurate for measures close to the delta measure $\rho^*$. The value $\lambda^*$ is $n$-dependent, so for any fixed measure $\rho\in \mathrm P(\zeta',Q')$ the bound \eqref{eq:lnubcl} is necessarily suboptimal for all steps $n$ except for a finite number of them. 

This result is in stark contrast to its counterpart for our condition \eqref{eq:rhozeta1} described by Theorem \ref{ther:worst_case_loss}. Specifically, if $p_n^2(\lambda)$ is monotone decreasing, the loss $L_n$ is maximized by the exact power-law measure $\rho(d\lambda)=Qd\lambda^\zeta$. In the more general case of non-monotone $p_n^2(\lambda)$, the mass of the worst-case measure becomes partially redistributed towards the local maxima of $p_n^2(\lambda)$ while still being rather well-distributed overall (see proof of Theorem \ref{ther:worst_case_loss} for details). For problems with approximately power-law spectral measures, such well-distributed character of the worst-case measure results in accurate upper bounds for all steps $n$.

As an example of application of Eq. \eqref{eq:lnubcl}, consider vanilla GD with learning rate $\alpha<2$. The respective polynomial is $p_n(\lambda)=(1-\alpha\lambda)^n$. The position of the Dirac delta can be found exactly by differentiating $\lambda^{\zeta'}(1-\alpha\lambda)^{2n}$ and is given by $\lambda^*=\alpha^{-1}\tfrac{\zeta'}{2n+\zeta'}$. Substituting this into \eqref{eq:lnubcl} gives
\begin{equation}\label{eq:sourcecond_GD_UB}
    L^{UB}_n(\zeta',Q') = \frac{Q'}{2} \left(1-\frac{\zeta'}{2n+\zeta'}\right)^{2n}\left(\frac{\zeta'}{\alpha(2n+\zeta')}\right)^{\zeta'} \overset{n\to\infty}{=\joinrel=} \frac{Q'}{2}\left(\frac{\zeta'}{2\alpha e}\right)^{\zeta'} n^{-\zeta'}(1+o(1)).
\end{equation}
This $O(n^{-\zeta'})$ bound seems reasonable, but we will see later that it is suboptimal: it can only hold when the true loss does not have a power-law behavior with the same exponent~$\zeta'$. 

\begin{figure}
     \centering
     {\includegraphics[scale=0.45, trim=33mm 82mm 35mm 7mm, clip]{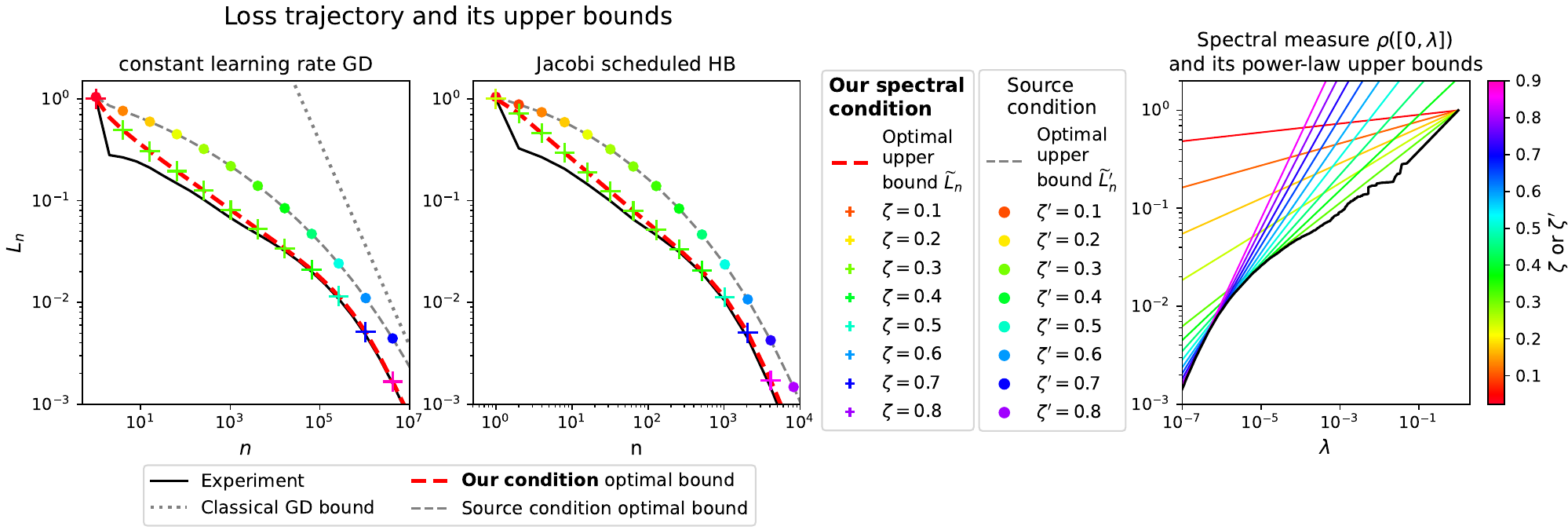}}
    
     \caption{Comparison of the experimental loss and different upper bounds for a kernel regression on the MNIST dataset. \textbf{Left:} Loss trajectories and respective bounds for constant learning rate GD (left subfigure) and Jacobi scheduled HB (right subfigure). For both GD and HB, the two upper bound curves are given by the functions $\widetilde L_n, \widetilde L'_n$ defined in Eqs. \eqref{eq:optimal_UB_our_condition},\eqref{eq:optimal_UB_source_condition}; in the GD case we additionally show the crude classical bound \eqref{eq:gd_bound_polyak} corresponding to $\zeta=1$ and requiring a very large constant $C$. The colors of the dots reflect the optimized values of $\zeta,\zeta'$ in Eqs. \eqref{eq:optimal_UB_our_condition}, \eqref{eq:optimal_UB_source_condition}. \textbf{Right:} The actual spectral distribution $\rho$ and different spectral bounds $\rho((0,\lambda])\le Q\lambda^{\zeta}$ with varying $\zeta$ and respective optimal $Q.$
     }
     \label{fig:MNIST_UBs}
\end{figure}

\paragraph{A practical example.}
The above arguments suggest that our spectral condition and respective bounds should be more efficient than the classical source condition and respective bounds for problems with approximate power-law spectra. In Figure  \ref{fig:MNIST_UBs} we verify this conclusion experimentally on a kernel regression problem for the MNIST dataset, optimized either with constant learning rate GD or HB with Jacobi-based schedule \eqref{eq:Jacobi_ansatz_parameters} (see Section \ref{sec:exp_details} for further details).

For each step $n$ and a given distribution $\rho$,  we compute the respective optimal bounds $\widetilde{L}_n(\rho), \widetilde{L}'_n(\rho)$ obtained with our and classical source condition by
\begin{align}
    \label{eq:optimal_UB_our_condition}
    \widetilde{L}_n(\rho) ={}& \inf_{\zeta,Q:\rho\in \mathrm P(\zeta,Q)}\sup_{\widetilde \rho \in \mathrm{P}(\zeta,Q)}L_n(\widetilde \rho),\\ \label{eq:optimal_UB_source_condition}
    \widetilde{L}'_n(\rho) ={}&\inf_{\zeta',Q':\rho\in \mathrm P'(\zeta',Q') }\sup_{\widetilde \rho \in \mathrm{P'}(\zeta',Q')}L_n(\widetilde \rho).
\end{align}
In either case, in the inner supremum we choose the tightest upper bound available for given parameters $\zeta,Q$ or $\zeta',Q'$, and then in the outer infimum optimize it over all admissible parameters.

We observe in Figure \ref{fig:MNIST_UBs} that the curves $\widetilde L_n$ corresponding to our spectral condition lie much closer to the actual loss trajectory than the curves $\widetilde L_n'$ corresponding to the classical source condition, in agreement with our prediction.  Accordingly, when using our spectral condition, the optimal $\zeta$ stays the same until the late stages of training ($n\sim 10^5$ for GD and $n\sim 10^3$ for HB), meaning that a single spectral condition with fixed $\zeta, Q$ can efficiently describe the loss evolution. In contrast, for the classical source condition \eqref{eq:sourcecond1}, the optimal parameters $\zeta', Q'$ are constantly changing along the whole optimization trajectory.

\paragraph{Theoretical suboptimality of the classical source condition.} We state now the theoretical suboptimality result announced earlier and corroborating theoretical expectations and the experimental observations. 
\begin{theor}\label{ther:source_cond_suboptimality}
    Assume that, for a certain spectral measure $\rho$, the sequence of the loss values under GD with constant learning rate $\alpha<1$ is $L_n = C n^{-\xi}(1+o(1))$. Then, the respective optimal upper bounds $\widetilde{L}_n',\widetilde{L}_n$ defined in Eqs. \eqref{eq:optimal_UB_our_condition}, \eqref{eq:optimal_UB_source_condition} are given by
    \begin{align}\label{eq:optima}
        \widetilde{L}_n &= \left[\frac{Q}{2}\Gamma(\xi+1)(2\alpha)^{-\xi}\right] n^{-\xi}(1+o(1)), \quad Q=\underset{\lambda\in(0,1]}{\sup} \rho([0,\lambda])/\lambda^\xi < \infty,\\
        \label{eq:optima_source}
        \widetilde{L}_n' &= \left[C\frac{\xi^{\xi+1}}{\Gamma(\xi+1)e^{\xi-1}}\right] \log(n) n^{-\xi} (1+o(1)).
    \end{align}
\end{theor}
This result shows that if the actual loss decreases as a power law, then the optimal upper bound \eqref{eq:optima} based on our spectral condition will agree with the actual loss up to a constant factor, while the optimal bound \eqref{eq:optima_source} based on the classical source condition will be off by at least a factor of $\log n$, even when we optimize the bound over the parameters $Q',\zeta'$. 

In the remainder of this section, let us outline the proof of Theorem \ref{ther:source_cond_suboptimality} (see Section \ref{sec:proof_th_our_vs_classical} for details). First, we show by tauberian-type arguments that the loss asymptotic $L_n = C n^{-\xi}(1+o(1))$ implies a respective power-law asymptotic of the spectral measure: $\rho([0,\lambda])=Q_\rho\lambda^\xi(1+o(1)),$ with $Q_\rho=2C \tfrac{(2\alpha)^\xi}{\Gamma(\xi+1)}$. One can think of this as a partial converse (for $\beta=0$) of theorem \ref{ther:constant_lr_bounds}, hence the value of the constant $Q_\rho$. 

Next, consider the exact power-law measure $\rho_\xi([0,\lambda])=Q_\rho\lambda^\xi$. While the full proof needs to carefully take into account the correction $\rho-\rho_\xi$ at finite $\lambda$ (in particular, leading to $Q>Q_\rho$ in \eqref{eq:optima}), the exact power-law measure captures the essence of the optimal bounds \eqref{eq:optima}, \eqref{eq:optima_source}. The optimal bound \eqref{eq:optima} for our condition is basically given by $L^{(\xi)}_n$ from theorem \ref{ther:constant_lr_bounds}, since for the exact power-law measure $\rho_\xi$ we have $\widetilde{L}_n(\rho_\xi)=Q_\rho \overline{L^{(\xi)}_n}$. 

Turning to the second result \eqref{eq:optima_source}, we note that the inner supremum in \eqref{eq:optimal_UB_source_condition} is already derived in \eqref{eq:sourcecond_GD_UB}. As for the outer infimum in \eqref{eq:optimal_UB_source_condition}, the smallest possible $Q'$ at a given $\zeta'$ can be inferred from lemma \ref{lem:conditions_comp}: $Q'(\zeta')=Q_\rho \frac{\xi}{\xi-\zeta'}$. From this point, we only need to estimate the optimal $\zeta'$ at a given iteration $n$: 
\begin{equation}\label{eq:optima_source_exact_powerlaw}
    \widetilde{L}'_n(\rho_\xi) = \inf_{0<\zeta'<\xi} \frac{Q_\rho}{2} \frac{\xi}{\xi-\zeta'} \left(\frac{\zeta'}{2\alpha e}\right)^{\zeta'} n^{-\zeta'}(1+o(1)) \stackrel{(*)}{=}\frac{Q_\rho \xi^{\xi+1}}{2(2\alpha e)^\xi} n^{-\xi}(1+o(1))\inf_{0<\varepsilon<\xi}\frac{n^\varepsilon}{\varepsilon}.
\end{equation}
Here in $(*)$, we took out all the factors that behave regularly at $\zeta'=\xi$, while the last infimum over $\varepsilon=\xi-\zeta'$ captures the essential tradeoff within the classical source condition: higher values of $\zeta'$ are more favorable on the level of the rate $O(n^{-\zeta'})$ but they come at a price of a large constant $\propto \frac{1}{\xi-\zeta'}$. The logarithm in \eqref{eq:optima_source} appears as a result of this tradeoff:
\begin{equation}
    \varepsilon^*_n\equiv\argmin_{\varepsilon>0} \frac{n^\varepsilon}{\varepsilon} = \frac{1}{\log n},\qquad \inf_{\varepsilon>0} \frac{n^\varepsilon}{\varepsilon}=\frac{n^{\varepsilon^*_n}}{\varepsilon^*_n}=e\log n.
\end{equation}

\section{Experiments}\label{sec:exp}

\begin{figure*}[t]
    \centering
    {\includegraphics[width=0.9\textwidth,clip,trim= 2mm 93mm 2mm 5mm ]{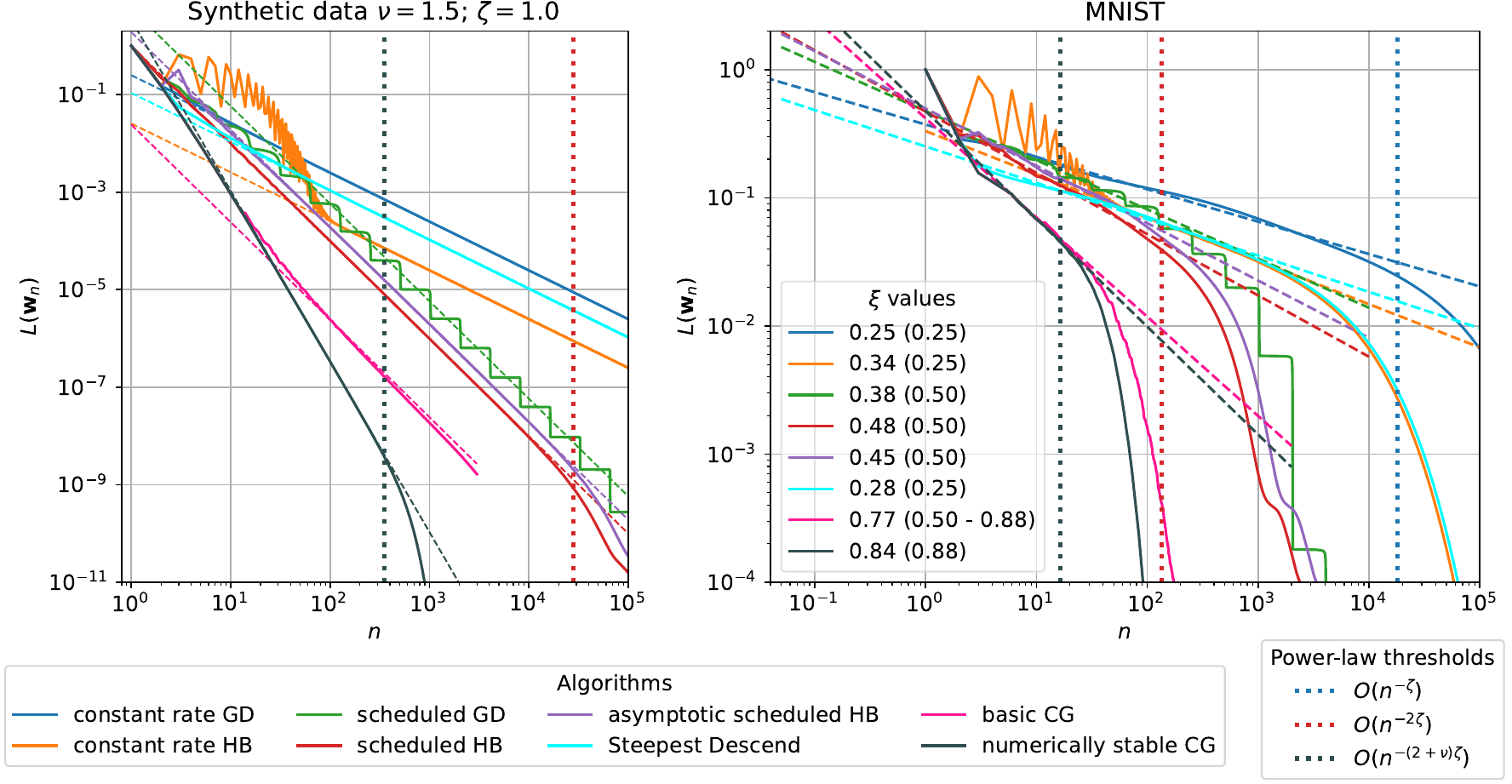}}
    
    \vspace{-1mm}
    \caption{Loss trajectories of different optimization algorithms for the artificial diagonal problem with $\nu=1.5$ and $\zeta=1$ \textbf{(Left)}, and for (the $\{0,\ldots,9\}$-valued version of) MNIST learned by the NTK kernel of shallow ReLU network \textbf{(Right)}. The dashed lines are the fitted power-laws; the fitted (and calculated theoretically) exponents are shown in the legend. Dotted vertical lines correspond to the estimated threshold of validity of loss power-law \eqref{eq:power_validity_thresholds}. 
    } 
  \label{fig:Losses}
\end{figure*}

\begin{figure}[t]
    \centering
    \setlength{\fboxsep}{0pt}
    {\includegraphics[width=0.9\textwidth,clip,trim= 1mm 70mm 0mm 1mm ]{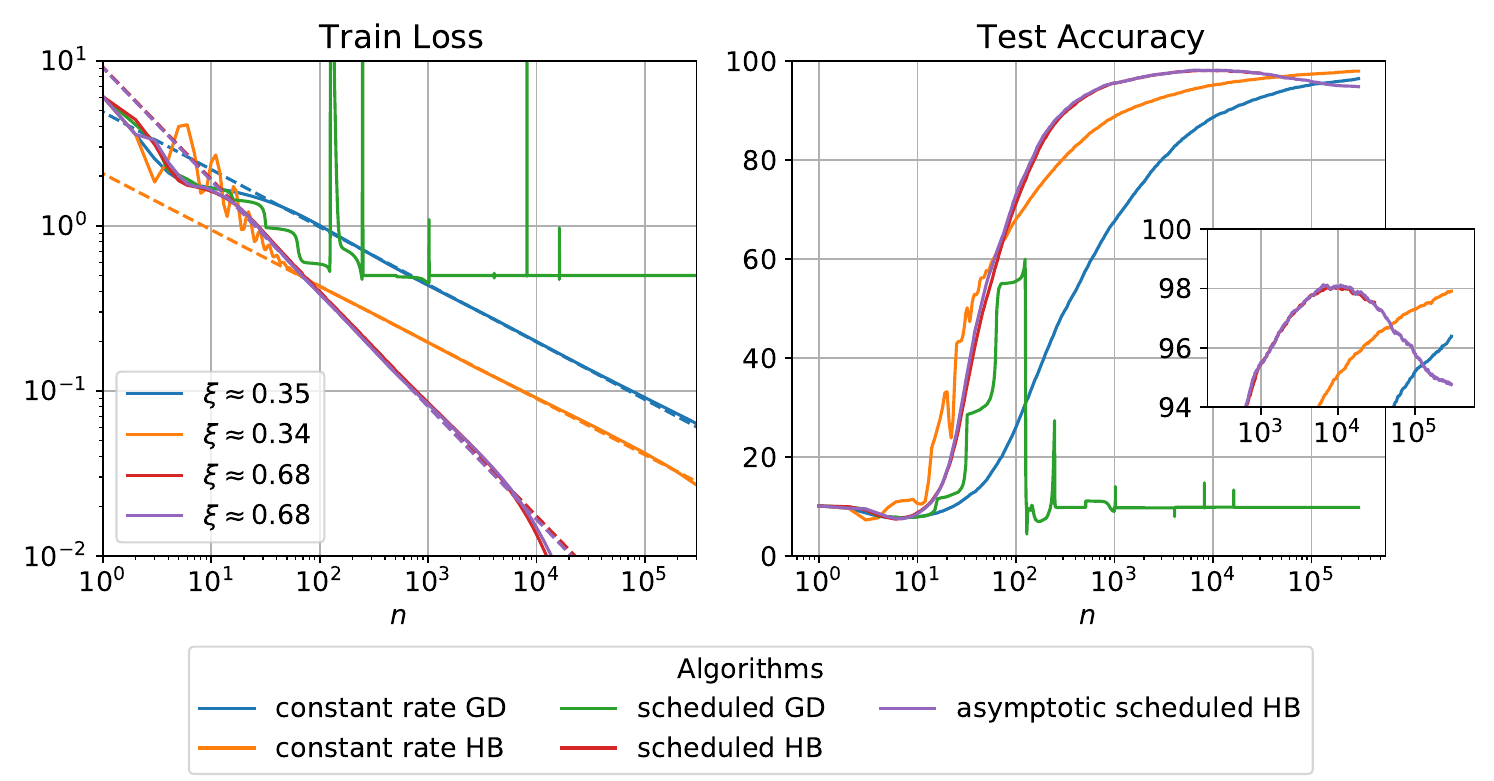}}
    
    \vspace{-1mm}
    \caption{Train loss and test accuracy of different algorithms on MNIST learned by a shallow width-1000 neural network.} 
  \label{fig:MNIST_net}
\end{figure}

\begin{figure}[thb!]
    \centering
    \setlength{\fboxsep}{0pt}
    {\includegraphics[width=0.9\textwidth,clip,trim= 5mm 5mm 3mm 0mm ]{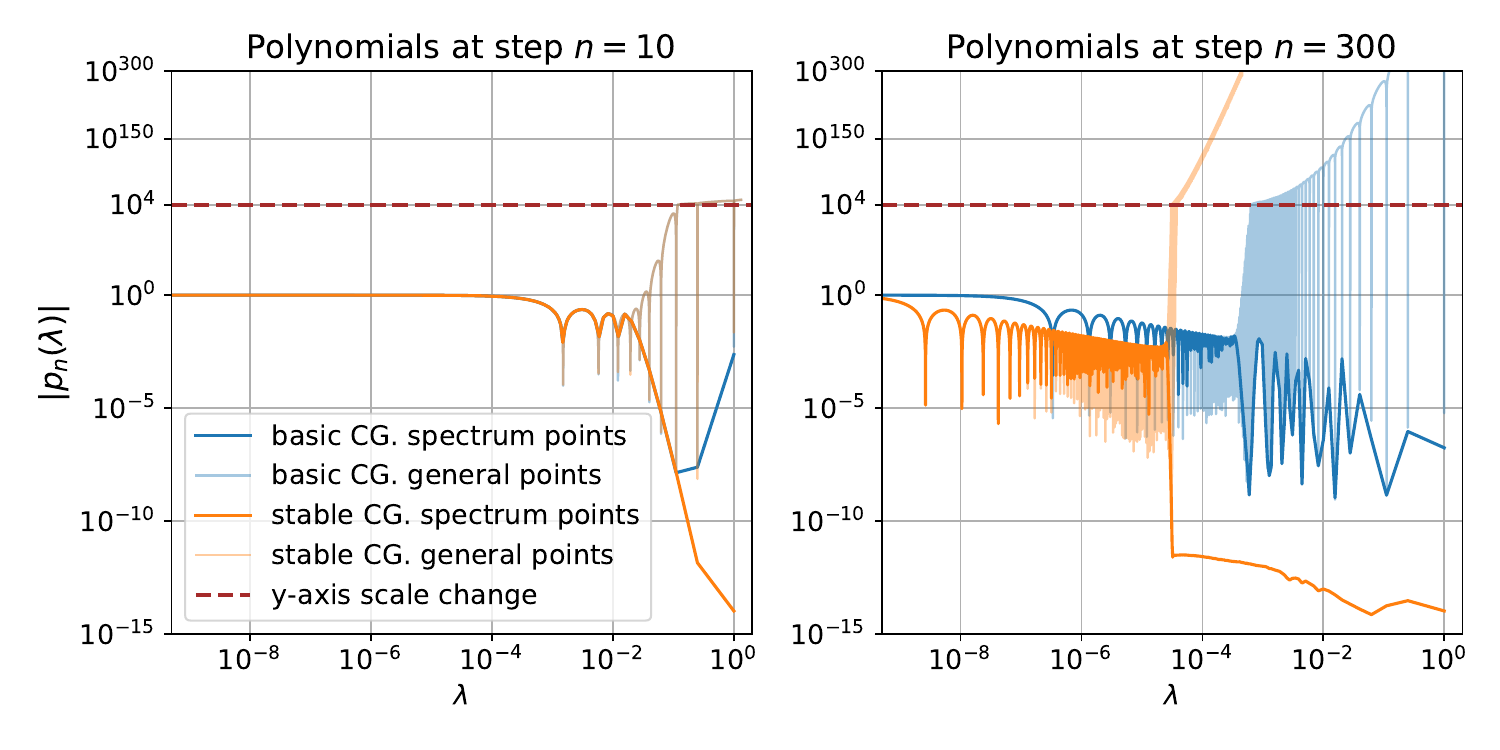}}
    
    \vspace{-3mm}
    \caption{Residual polynomials for stable and unstable CG.} 
  \label{fig:CG_polynomials}
\end{figure}

\paragraph{Diagonal matrices.} We start with an artificial quadratic problem in which we can directly control the exponents $\zeta$ and $\nu$: $\widetilde A$ is diagonal with eigenvalues $\lambda_k = k^{-\nu}$, and the respective coefficients of $\mathbf f_*$ are $c_k = k^{-\frac{\zeta\nu-1}{2}}$. The size of $\widetilde{A}\in\mathbb{R}^{M\times M}$ is $M=10^6$. The optimization results are shown in Figure \ref{fig:Losses} (Left). For all considered algorithms except CG, the losses have power-law rates with exponents $\xi$ in accordance with Table \ref{tab:bounds} (shown by dashed lines). In Figure \ref{fig:Losses} the asymptotic scheduled HB algorithms are defined using the simplified versions of learning rate and momenta, obtained by discarding the $O(\ldots)$ terms in Eq. \eqref{eq:Jacobi_ansatz_parameters}. While we do not have a theoretical convergence rate for this method, we see that it has the same rate $O(n^{-2\zeta})$ as the full scheduled HB. This suggests that the correct asymptotic of $\alpha_n,\beta_n$ at $n\to\infty$ is a deeper reason for acceleration.   

CG has the expected $\sim n^{-(2+\nu)\zeta}$ asymptotic only up to iteration $n_e \approx 20$, around which the asymptotic switches to $\sim n^{-2\zeta}$. This happens because of numerical errors (see further discussion in  paragraph ``CG polynomials'' below). A version of CG modified to ensure stability exhibits the $\sim n^{-(2+\nu)\zeta}$ convergence to the very end. 

Scheduled GD has a ``staircase'' shape because the schedule consists of size-$2^l$ chunks  (see Section \ref{sec:ubgdpredefined}). 

Note that for faster algorithms, such as CG or Jacobi scheduled HB, the power-law behavior of the loss breaks down at sufficiently large iteration $n$. This iteration can be estimated theoretically, as we explain below, and is depicted by vertical dotted lines in Figure \ref{fig:Losses}.

\paragraph{Intervals of validity of loss power laws.} When applied to real-life problems with approximately power-law spectra, the respective power-law behavior of optimization loss trajectories can be expected to hold only for moderately large iterations. In a real-life finite dimensional problem, the infinite-dimensional approximation $\rho([0,\lambda])\sim \lambda^\zeta$ breaks down for  $\lambda\lesssim\lambda_{\text{low}}$ with some characteristic value $\lambda_{\text{low}}$ (e.g., the minimal positive eigenvalue). Under optimization, the loss is given by $L(\mathbf w_n)=\tfrac{1}{2}\int_{0}^{\lambda_{\max}}p_n^2(\lambda)\rho(d\lambda)$ with a suitable residual polynomial $p_n$. At large $n$, under the assumption of a power-law measure $\rho$,  the leading contribution to this integral comes from the spectral interval $(0,\lambda_{\text{low}}).$ Accordingly, the loss power law $L(\mathbf w_n)\propto n^{-s\zeta}$, where $s=1,2$ or $2+\nu,$ breaks down for $n\gtrsim n_{\text{th}}$ with some characteristic iteration number $n_{\text{th}}$.
 In Section \ref{sec:powerlaw_end} we derive a (non-rigorous) estimate of $n_{\text{th}}:$
\begin{equation}\label{eq:power_validity_thresholds}
    n_\text{th} \propto 
    \begin{cases}
    \vspace{2mm}
    \frac{1-\beta}{\alpha \lambda_{\text{low}}}, \qquad& \text{for constant rate algorithms }(s=1)\\
    \vspace{2mm}
    \frac{1}{\sqrt{\lambda_{\text{low}}}}, \qquad& \text{for algorithms based on Jacobi polynomials }(s=2)\\
    \vspace{2mm}
    \lambda_{\text{low}}^{-\frac{1}{\nu+2}}, \qquad& \text{for (numerically stable) Conjugate Gradients }(s=2+\nu)
    \end{cases}
\end{equation}
In the experiments, we choose $\lambda_{\text{low}}$ either as the minimum eigenvalue (in the artificial power-law problems) or as a value at which we experimentally observe the breakdown of the spectral power-law (for MNIST).

\paragraph{Realistic quadratic problems.} As an example of a realistic quadratic problem we take a subset of MNIST (of size $M=30000$) and consider the scalar regression problem with targets given by the numerical values of corresponding digits $y\in\{0,1,\ldots,9\}$. The matrix $\widetilde{A}$ is  the NTK of an infinitely wide, single-hidden layer network. The results are depicted in Figure \ref{fig:Losses} (right). Again, we observe power-law dependencies up to the estimated thresholds. The numerical entries in the legend have the form $\xi_{\text{exp}} (\xi_{\text{theor}})$, where $\xi_{\text{exp}}$ is the ``experimental'' exponent  estimated directly from the loss trajectory, and $\xi_{\text{theor}}$ is the respective ``theoretical'' exponent  given by $\zeta, 2\zeta$, or $(2+\nu)\zeta$. Here the values $\nu\approx 1.37$ and $\zeta\approx 0.25$ are in turn estimated from the empirically found $\rho$ (Figure \ref{fig:MNIST_loss_and_spectral_distributions} (right)) and the eigenvalues $\lambda_k$ and partial sums of target expansion coefficients (Figure \ref{fig:MNIST_loss_and_spectral_distributions} (center)). We see a reasonable agreement between $\xi_{\text{exp}}$ and $\xi_{\text{theor}}$. Like with synthetic data, the asymptotic scheduled HB performs similarly to its full counterpart.

\paragraph{Neural networks.} We consider a shallow fully-connected ReLU network with 1000 hidden neurons and train it on the full MNIST with MSE loss calculated on one-hot encoded classes. Note that this is no longer a quadratic problem. We restrict ourselves to optimization algorithms with predefined schedules due to their computational efficiency compared to adaptive algorithms (in which the 1D nonlinear problem of step optimization has to be solved in each iteration). Also, we use full-batch gradient descent in accordance with the rest of the paper. The results are shown in Figure~\ref{fig:MNIST_net}. 

For all algorithms except scheduled GD we observe behavior similar to the quadratic case, and in particular asymptotic HB is very close to its full counterpart. The relation between the fitted exponents $\xi$ holds true: they are twice as large for scheduled methods as for constant learning rate methods. 

The unstable behavior of scheduled GD is explained by large step-sizes $\alpha_n$ present in the schedule (at steps $n\approx 2^l$). When the problem is quadratic, large $\alpha_n$ are compensated by smaller ones chosen at other steps $n$, but non-quadratic perturbations break this compensation mechanism.

\paragraph{CG polynomials.} In Figure \ref{fig:CG_polynomials} we plot CG polynomials $p_n(\lambda)$ for the basic and the numerically stable algorithms, calculated either at the spectral points $\lambda_k$, or also between them. At step $n=10$ the two polynomials mostly coincide except for big $\lambda$. At step $n=300$ the polynomials are different, and for either of them we observe two $\lambda$--regions with a sharp transition point $\widetilde\lambda$. For $\lambda>\widetilde\lambda$, the values of $p_n(\lambda)$ vanish at the spectral points but are extremely large in between, meaning that the roots of $p_n(\lambda)$ are located exactly at the spectral points $\lambda_k$. The rest of the roots are located at $\lambda<\widetilde\lambda$ and seem to optimize the overall envelope of $p_n(\lambda)$ instead of only root positions. This agrees with construction used in upper bound \eqref{eq:cgdisclwn}. As, due to numerical errors, the polynomial of basic CG places its roots in the region $\lambda>\widetilde\lambda$ with lower precision, the value of $\widetilde\lambda$ is higher in this case and hence convergence on $[0,\widetilde\lambda]$ is worse.        

\section{Conclusion}
We have considered a wide range of first-order optimization methods including Gradient Descent, Steepest Descent, Heavy Ball, and Conjugate Gradients, with constant, non-constant predefined, and adaptive learning rates. Under power-law spectral assumptions with target exponent $\zeta$ and eigenvalue exponent $\nu$ the convergence rates of these methods are given by $O(n^{-\xi})$, where $\xi=\zeta,2\zeta$ or $(2+\nu)\zeta$, depending on the method. The basic rate with $\xi=\zeta$ applies to Gradient Descent with constant learning rates and also to Steepest Descent. To reliably achieve the first accelerated rate $2\zeta$ with Heavy Ball, a specific Jacobi-based schedule of learning rate and momenta is required, with $\beta_n$ approaching 1 so that $1-\beta_n \propto n^{-1}$. Finally, the fastest rate $(2+\nu)\zeta$ is achieved by Conjugate Gradients -- the only method out of those we have considered that can take advantage of the discreteness of the problem spectrum by exactly fitting the target function in certain eigenspaces.    

We prove that all our upper bounds are tight. For each upper bound we provide an example problem whose convergence rate matches that of the upper bound, and in some cases also  has a very close coefficient. An important aspect of our approach is a power-law spectral assumption that is somewhat different from the classical source condition. We show, both experimentally and theoretically, that our condition much better describes problems whose actual loss trajectory is well approximated by a power-law. Specifically, for a problem with power-law loss asymptotic $L(\mathbf{w}_n)\sim n^{-\xi}$ our condition provides the matching bound $L(\mathbf{w}_n)\leq const \cdot n^{-\xi}$ while the best usage of the classical source condition can only provide a bound with additional logarithmic factor, $L(\mathbf{w}_n)\leq const \cdot n^{-\xi} \log n$.  

Our theoretical results are confirmed by experiments with both simulated and real problems, including classifying MNIST by a neural network (which is only an approximately quadratic problem). In all experiments we observe a clear power law dependence of the loss on the optimization step $n$ for steps that are neither too large nor too small, i.e. whenever both the infinite-dimensional approximation and asymptotic formulas are applicable. The respective exponents and their mutual relations agree well with theoretical predictions (unless the method is affected strongly by noise, as with CG, or by non-quadratic corrections, as with the optimally scheduled GD applied to a neural network).

Finally, let us outline a few natural topics for future research. First, as discussed in Section \ref{sec:exact_power-law_measure}, Heavy Ball with various Jacobi-based schedules with the asymptotic form $\alpha_n\sim const, \; 1-\beta_n \propto n^{-1}$ can ensure the same $O(-2\zeta)$ convergence rate. We hypothesize that under the general spectral condition $\rho([0,\lambda])\leq G(\lambda)$, the asymptotic of $G(\lambda)$ at small eigenvalues $\lambda\to0$ can be translated into a certain asymptotic of $1-\beta_n$ at large iterations $n$ for optimal HB. Second, it would be interesting to investigate whether weak non-quadratic perturbations of quadratic problems allow to retain the accelerated rate $O(n^{-2\zeta})$. Our experiments with a neural network on MNIST confirm this possibility. Third, it would be interesting to include stochasticity into consideration, as mini-batch stochastic gradient descent is a necessary requirement for any GD method to be used in modern deep learning applications.

\acks{We acknowledge support from the Russian Ministry of Science and Higher Education, grant No. 075-10-2021-068.}

\newpage

\appendix

\addcontentsline{toc}{section}{Appendix}

\section{Related work}\label{sec:relwork} 
\paragraph{Optimization by GD, SD, HB and CG under power law spectral assumptions.} The first study of GD, HB and CG under power-law spectral assumptions (in a form somewhat different from ours; see discussion at the end of Section \ref{sec:probdef}) was performed in \cite{nemirovskiy1984iterative1} (upper bounds) and \cite{nemirovsky1984iterative2} (lower bounds). These two works proved or conjectured some of the bounds appearing in our Table \ref{tab:bounds}. While these two papers only considered scheduled HB based on Chebyshev polynomials, \cite{brakhage1987ill} generalized it to a ``$\nu$-method'' based on general Jacobi polynomials, which allowed him to obtain the tight $O(n^{-2\nu})$ upper bound analogous to our Corollary \ref{corol:HB_upper_bound} for HB with predefined schedules. SD was analyzed in \cite{gilyazov2013regularization} who proved a $O(n^{-\zeta})$ upper bound (their Theorem 2.2.5). However, the proof of its tightness (supplemented in our Theorem \ref{th:sd}) does not seem to have been known prior to our work. Various aspects of optimization by HB and CG were discussed in \cite{hanke1991accelerated} and \cite{hanke1996asymptotics}. In particular, the latter paper gave a proof of the lower bound for CG in the special case of exponents $\nu=1,2$.
All of these works relied on the classical source condition and only considered problems with attainable solutions.  

The recent work \cite{Berthier2020_gossip}, although focusing on a specific application domain of gossip problem, uses a spectral condition (see their Proposition 5.5 or Definition I.2) which is different from the classical source condition and much closer to our condition, and also considers a Jacobi-based optimization algorithm. However, both \cite{Berthier2020_gossip} and earlier works \cite{brakhage1987ill,hanke1991accelerated} rely on classical asymptotic properties of Jacobi polynomials for the proofs of upper bounds, e.g. Theorem 7.32.2 of \cite{Szvego1939}. This approach quickly provides the desired $O(n^{-2\zeta})$ rate but does not specify the constant. In contrast, our flattened polynomial construction of Theorem \ref{ther:worst_case_loss} followed by accurate estimations in Theorem \ref{ther:Jacobi_worst_case_UB} and Proposition \ref{prop:parabola_powerlaw_av_minimum} lead to an explicit and tight constant in the convergence bound (e.g. overestimation by at most a factor of $C_\zeta=4$ for $\zeta=1$).

\paragraph{SGD.} Analogs of our power law spectral conditions \eqref{eq:rhozeta} and \eqref{eq:lkcknu} are well-known in literature on kernel methods, regularized regression and SGD \citep{caponnetto2007optimal, steinwart2009optimal, varre2021last}. Convergence of SGD under these or similar conditions has been studied in \cite{berthier2020tight, zou2021benign, varre2021last, velikanov2022view}. SGD subsumes GD as a special case of noiseless gradient evaluation, but is in a sense more complex than all the algorithms we discuss in this paper because even for linear models the loss evolution under SGD is not generally expressible in terms of only spectral data.  The most common version of SGD is mini-batch SGD in which the stochasticity is due to random sampling of the underlying data. In contrast to GD, SD and HB (cf. Table \ref{tab:bounds}), convergence rates of SGD do depend directly, in general, on the eigenvalue decay exponent $\nu$. In particular, for mini-batch SGD with constant learning rates the respective exponent equals $\min(\zeta, 2-1/\nu)$; moreover, optimization diverges if $\nu<1$.

\paragraph{Kernel methods and NTK.}
Power law eigenvalue decay bounds are known to generally hold for integral operators with kernels satisfying suitable regularity assumptions \citep{widom1963asymptotic, kuhn1987eigenvalues, ritter1995multivariate, ferreira2009eigenvalues, Birman_1970, williams2006gaussian}. 

In the NTK regime of training wide neural networks the network model essentially becomes a kernel model \citep{neal2012bayesian, jacot2018neural} with explicit kernels \citep{NIPS2009_kernel, NEURIPS2019_0d1a9651}. Several recent studies empirically verify and exploit power law assumptions for the NTK spectrum \citep{bahri2021explaining, canatar2021spectral, lee2020finite, nitanda2021optimal, jin2021learning}. Specific powers of eigenvalue decay and eigenfunction expansion coefficients for ReLU networks and some classes of target functions are derived in \cite{velikanov2021explicit}.

\paragraph{Steepest Descent.} See \cite{kantorovich1964functional} for a general introduction to Steepest Descent. In a general non-strongly convex case, convergence of the iterates to a solution (if it exists) was proved in \cite{fridman1952steepest}. In \cite{kammerer1971steepest} an explicit $\|\mathbf w_n-\mathbf w_*\|^2=O(n^{-1})$ bound was proved in the non-strongly convex case under assumption $\|\widetilde A^{-1}\mathbf f_*\|<\infty$. The $O(n^{-\zeta})$ convergence upper bound under a power-law spectral condition was proved in \cite{gilyazov2013regularization} using moment inequalities from \cite{krasnoselskii52approximate}. Our approach in Section is \ref{sec:sd} is rather different from these works and relies on the observation that SD converges to a period-2 oscillatory regime. This effect was established by  \cite{akaike1959successive} in the finite-dimensional setting and by \cite{ pronzato2001renormalised} in the infinite-dimensional setting. Compared to \cite{gilyazov2013regularization}, our approach is applicable under our slightly weaker spectral assumption \eqref{eq:rhozeta1} and additionally proves the tightness of the loss upper bound.

\paragraph{Heavy Ball.} Multi-step methods have long been used in numerical linear algebra. As a method of optimization for general (non-quadratic) problems, Heavy Ball was proposed in \cite{polyak1964some}. HB can be interpreted as a simplest method with the momentum term \citep{qian1999momentum}.  \cite{flammarion2015averaging} introduced a general family of methods that includes HB with $\beta_n= 1-2/n$ as well as averaged GD \citep{polyak1992acceleration}. Some variants of GD with momentum are optimal with respect to averaged case optimization scenarios \citep{pedregosa2020acceleration, lacotte2020optimal}.

\paragraph{Conjugate Gradients.} Method of Conjugate Gradients was proposed in \cite{Hestenes&Stiefel:1952} and extensively studied afterwards  \citep{daniel1971approximate, hestenes2012conjugate}. The extension of the method to non-quadratic problems was first proposed in \cite{fletcher1964quadratically}. Stability of CG is a complex issue that has also been analyzed extensively \citep{Hestenes&Stiefel:1952, bjorck1998stability, meurant2006lanczos, fischer2011polynomial}. A $\|\mathbf w_n-\mathbf w_*\|^2=O(n^{-1})$ convergence bound for CG in a gapless infinite-dimensional setting was proved in \cite{kammerer1972convergence}. A version of the bound $L(\mathbf w_n)=O(n^{-2\zeta})$ was proved in \cite{nemirovskiy1984iterative1}, and in the same paper it was observed that this rate can be improved if the spectrum is discrete. \cite{hanke1991accelerated,hanke1996asymptotics} gave a version of the $L(\mathbf w_n)=O(n^{-(2+\nu)\zeta})$ bound and proved its tightness in the cases $\nu=1,2$, for which a classical system of orthogonal polynomials is available.
Our general proof of the tightness of the $O(n^{-(2+\nu)\zeta})$ bound for CG under the power law eigenvalue decay assumption (Section \ref{sec:lb}) is inspired by Theorem 2.1.7 in \cite{nesterov2003introductory} which proves the tightness of the bound $L(\mathbf w_n)=O(n^{-2})$ in a setting of finite norm solution $\|\mathbf w_*\|<\infty$. However, the proof of our bound is significantly more difficult. 

\paragraph{Nesterov Accelerated Gradient (NAG).} NAG \citep{nesterov1983method} is a modification of Heavy Ball \eqref{eq:mgd1} in which the gradient is computed after applying the momentum term rather than before:
\begin{align}
    \mathbf w_{n+1}
    ={}&\mathbf w_n+\beta_n(\mathbf w_{n}-\mathbf w_{n-1})-\alpha_n\nabla L(\mathbf w_n+\beta_n(\mathbf w_{n}-\mathbf w_{n-1})).
\end{align}
For quadratic problems, the analog of Eq. \eqref{eq:mgd2} then reads
\begin{align}
    \mathbf w_{n+1}
    ={}&\mathbf w_n+\beta_n(\mathbf w_{n}-\mathbf w_{n-1})-\alpha_n[(A \mathbf w_n-\mathbf b)+\beta_nA(\mathbf w_n-\mathbf w_{n-1})].
\end{align}
NAG is a practically widely used method and it is known to provide improved upper bounds for general convex problems \citep{nesterov1983method}. However, it does not seem to improve on Heavy Ball in the purely quadratic case considered in the present paper, at least in terms of the optimal convergence exponent. Specifically, assuming that the coefficients $\alpha_n,\beta_n$ are non-adaptive (predefined), both NAG and Heavy Ball are subject to our Theorem \ref{ther:discrete_spectra_predefined_schedule_LB} showing that they cannot generally have a rate $L(\mathbf w_n)=O(n^{-\xi})$ with $\xi>2\zeta$, while the rate $L(\mathbf w_n)=O(n^{-2\zeta})$ is attained by Heavy Ball by Corollary \ref{corol:HB_upper_bound}.

\section{Background on polynomials for optimization} 
\label{sec:poly}
\paragraph{The polynomial representation of optimization updates.} 
The optimization algorithms of Section \ref{sec:algo} and their properties can be conveniently expressed in terms of polynomials of the operator $A$ (or $\widetilde A$). Suppose first for simplicity that our optimization problem $L(\mathbf w)=\tfrac{1}{2}\langle\mathbf w, A\mathbf w\rangle-\langle\mathbf w,\mathbf b\rangle+\tfrac{1}{2}\|\mathbf f_*\|^2\to\min_{\mathbf w}$ has a finite-norm optimizer $\mathbf w_*$ such that $A\mathbf w_*=\mathbf b$. Consider the deviations $\delta\mathbf w=\mathbf w-\mathbf w_*$ of the points $\mathbf w$ from the solution $\mathbf w_*$. For the basic GD or SD, we have
\begin{align}
    \delta \mathbf w_{n}
    ={}&\mathbf w_{n-1}-\alpha_{n-1}(A\mathbf w_{n-1}-\mathbf b)-\mathbf w_* \\
    ={}&(1-\alpha_{n-1}A)\delta\mathbf w_{n-1},
\end{align}
and so, by iterating,
\begin{align}
    \delta \mathbf w_{n} = p_n(A)\delta\mathbf w_0,
\end{align}
where $p_n$ is the degree-$n$ polynomial
\begin{equation}
    p_n(\lambda) = \prod_{s=1}^n(1-\alpha_{s-1}\lambda).
\end{equation}
The respective loss is
\begin{align}
    L(\mathbf w_n)={}&\frac{1}{2}\langle A\delta\mathbf w_n,\delta\mathbf w_n\rangle \\
    ={}&\frac{1}{2}\int \lambda p_n^2(\lambda)\rho_{A,\mathbf w_*}(d\lambda)\\
    ={}&\frac{1}{2}\int p_n^2(\lambda)\rho_{\widetilde A,\mathbf f_*}(d\lambda),\label{eq:lwnpnwtaf}
\end{align}
where $\rho_{A,\mathbf w_*}$ and $\rho_{\widetilde A,\mathbf f_*}$ are the spectral measures associated (as in Eq. \eqref{eq:rhodef}) with $A,\mathbf w_*$ and $\widetilde A,\mathbf f_*$, respectively.

Representation \eqref{eq:lwnpnwtaf} (with $\rho_{\widetilde A,\mathbf f_*}$) can alternatively be reached without assuming the existence of the solution $\mathbf w_*$, by considering the deviations $\delta \mathbf f =\mathbf f-\mathbf f_*$ in the target space
and similarly observing that 
\begin{align}
    \delta \mathbf f_{n} = p_n(A)\delta\mathbf f_0,
\end{align}
with the same polynomial $p_n$.

In the case of HB and CG, the iterations have the more general form
\begin{equation}
    \delta\mathbf w_{n+1} = (1-\alpha_n A)\delta\mathbf w_n+\beta_n(\delta\mathbf w_n-\delta\mathbf w_{n-1}).
\end{equation}
This again yields the polynomial representation $\delta \mathbf w_{n} = p_n(A)\delta\mathbf w_0$, but with a degree-$n$ polynomial $p_n$ depending on $\{\alpha_s,\beta_s\}_{s=0}^{n-1}$ in a more complicated way:
\begin{align}
    p_0={}&1,\\
    p_1={}&1-\alpha_0\lambda,\\
    \label{eq:HB_poly_recurrence}
    p_{n+1} ={}& (1-\alpha_n \lambda)p_n+\beta_n(p_n-p_{n-1}).
\end{align}
Note that $p_n$ is necessarily a \emph{residual} polynomial, in the sense that $p_n(0)=1$.

As mentioned in Section \ref{sec:algo}, CG has the important property of being optimal among all first order methods generating new iterates $\mathbf w_{n+1}$ by shifting the initial point $\mathbf w_0$ along linear subspaces spanned by the previously computed gradients  $\nabla L(\mathbf w_0),\ldots,\nabla L(\mathbf w_{n})$. In terms of the respective residual polynomials $p_n$, this means that they minimize the loss functional over all residual polynomials of given degree:
\begin{equation}\label{eq:CG_polynomials_opt}
    p_n = \argmin_{q_n:\deg q_n = n,q_n(0)=1}\frac{1}{2}\int q_n^2(\lambda)\rho_{\widetilde A,\mathbf f_*}(d\lambda).
\end{equation}

See the book \cite{fischer2011polynomial} for more details on the polynomial representation of optimization methods.

\paragraph{Jacobi polynomials.} 

As shown in Section \ref{sec:exact_power-law_measure}, Jacobi polynomials $ P_{n}^{(a,b)}$ arise as an optimal choice for power-law spectral measure. We heavily use these polynomials in many of our results. 

The appearance of Jacobi polynomials in our setting is related to their orthogonality w.r.t. power-law weight function:
\begin{equation}
    \int\limits_{-1}^1 (1-x)^a(1+x)^b P^{(a,b)}_n(x) P^{(a,b)}_m(x) dx= C_n \delta_{nm}.
\end{equation}
Here $\delta_{nm}$ is Kronecker delta function and $C_n$ are the constants depending on normalization of the polynomials. We adopt the standard normalization of Jacobi polynomials by their value at $x=1$:
\begin{equation}\label{eq:Jacobi_normalization}
    P^{(a,b)}_n(1) = \binom{n+a}{n}.
\end{equation}
Jacobi polynomials, like any system of orthogonal polynomials, enjoy three-term recurrence relations. Specifically,
\begin{equation}\label{eq:Jacobi_recurrence}
    \begin{split}
        2(n+1)(n+a+b+1)&(2n+a+b) P_{n+1}^{(a,b)}(x) = \\ 
        &(2n+a+b)(2n+a+b+1)(2n+a+b+2) xP_{n}^{(a,b)}(x) \\
        &+ (2n+a+b+1)(a^2-b^2) P_{n}^{(a,b)}(x) \\
        &- 2(n+a)(n+b)(2n+a+b+2)P_{n-1}^{(a,b)}(x).
    \end{split}
\end{equation}

\section{Main spectral condition}\label{sec:spec_cond_comparison}
In this section, we collect the proofs of the results concerning either general properties of our spectral condition \eqref{eq:rhozeta1} or its relation to the classical source condition \eqref{eq:sourcecond1}.

\subsection{Basic properties}\label{sec:attain}

\paragraph{Proof of Lemma \ref{lem:conditions_comp}.} \hfill

\noindent
\emph{Inclusion $\mathrm{P}(\zeta,Q) \subseteq \mathrm{P}'(\zeta',Q')$} (Eq. \eqref{eq:ourcond_to_sourcecond}). To test this inclusion for a certain pair of $\zeta,Q$ and $\zeta',Q'$, we need to check
\begin{equation}\label{eq:zetaQ_in_zetapQp_test}
    \underset{\rho \in \mathrm{P}(\zeta,Q)}{\operatorname{\sup}} \int_0^1 \lambda^{-\zeta'} \rho(d\lambda) \leq Q'.
\end{equation}
First, consider $\zeta'\geq\zeta$ and the exact power-law measure $\rho(d\lambda)=Qd(\lambda^\zeta) \in \mathrm{P}(\zeta,Q)$. Then, the integral in \eqref{eq:zetaQ_in_zetapQp_test} diverges as
\begin{equation}
    \underset{\varepsilon\to0}{\lim} \int_\varepsilon^1 \lambda^{-\zeta'}Q d(\lambda^\zeta) = \begin{cases}
         \underset{\varepsilon\to0}{\lim} \; \frac{Q\zeta}{\zeta'-\zeta}\left(\varepsilon^{\zeta-\zeta'}-1\right)=\infty, \quad &\zeta'>\zeta \\
         \underset{\varepsilon\to0}{\lim} \; Q\zeta\log(\varepsilon^{-1})=\infty, \quad &\zeta'=\zeta 
    \end{cases}
\end{equation}
which makes $\zeta'<\zeta$ a necessary condition for inclusion. Assuming this condition, the supremum in \eqref{eq:zetaQ_in_zetapQp_test} can be evaluated using integration by parts:
\begin{equation}\label{eq:ourcond_to_sourcecond_ibp}
    \int_0^1 \lambda^{-\zeta'} \rho(d\lambda) = \lambda^{-\zeta'}\rho([0,\lambda])\Big|_0^1 +\zeta'\int_0^1\lambda^{-\zeta'-1}\rho([0,\lambda])d\lambda.
\end{equation}
Note that both terms in \eqref{eq:ourcond_to_sourcecond_ibp} are well defined thanks to the constraint $\rho([0,\lambda])\leq Q\lambda^\zeta$. Importantly, the right-hand side of \eqref{eq:ourcond_to_sourcecond_ibp} is a pointwise positive linear functional of the cumulative distribution function $\rho([0,\lambda])$, which implies that the supremum in \eqref{eq:zetaQ_in_zetapQp_test} is reached at the exact power-law measure $\rho(d\lambda)=Qd(\lambda^\zeta)$, and its value is
\begin{equation}\label{eq:zetaQ_in_zetapQp_test2}
    \underset{\rho \in \mathrm{P}(\zeta,Q)}{\operatorname{\sup}} \int_0^1 \lambda^{-\zeta'} \rho(d\lambda) = Q+Q\frac{\zeta'}{\zeta-\zeta'}=Q\frac{\zeta}{\zeta-\zeta'}.
\end{equation}
This computation implies that $Q'\geq Q\frac{\zeta}{\zeta-\zeta'}$ is equivalent to the desired inclusion for $\zeta'<\zeta$, which completes the proof of \eqref{eq:ourcond_to_sourcecond}.

\medskip 
\noindent
\emph{Inclusion} $\mathrm{P}'(\zeta',Q') \subseteq \mathrm{P}(\zeta,Q)$ (Eq. \eqref{eq:sourcecond_to_ourcond}). First note that this inclusion cannot hold if $\zeta'<\zeta$. Indeed, in that case the equivalence \eqref{eq:ourcond_to_sourcecond} would imply $\mathrm{P}(\widetilde{\zeta},\widetilde{Q})\subseteq\mathrm{P}(\zeta,Q)$ for any $\widetilde{\zeta}\in(\zeta',\zeta)$ and some $\widetilde{Q}$, which contradicts $\rho(d\lambda)=\widetilde{Q}d(\lambda^{\widetilde\zeta})\notin \mathrm{P}(\zeta,Q)$.

For $\zeta'\geq\zeta$ the inclusion can be tested with
\begin{equation}\label{eq:zetapQp_in_zetaQ_test}
    \underset{\rho \in \mathrm{P}'(\zeta',Q')}{\operatorname{\sup}} \; \left[\underset{\lambda\in(0,1]}{\sup}\rho([0,\lambda])/\lambda^\zeta \right] \leq Q,
\end{equation}
where we used that $\rho(\{0\})=0$ in our setting (see section \ref{sec:setting}) to account for $\lambda=0$ case of \eqref{eq:rhozeta1}. Note that the expression $\rho([0,\lambda])/\lambda^\zeta$ is bounded for $\rho \in \mathrm{P}'(\zeta',Q'), \; \lambda\in(0,1]$ as
\begin{equation}
    \lambda^{-\zeta}\rho([0,\lambda]) \leq \lambda^{\zeta'-\zeta}\int_0^{\lambda}\lambda_1^{-\zeta'}\rho(d\lambda_1)\leq Q'
\end{equation}
Actually, this bound is tight, as can be shown by taking $\rho=Q'\delta_1\in \mathrm{P}'(\zeta',Q')$ and $\lambda=1$. This makes the value of the supremum in \eqref{eq:zetapQp_in_zetaQ_test} equal to $Q'$, thus establishing equivalence \eqref{eq:sourcecond_to_ourcond}. This completes the proof of Lemma \ref{lem:conditions_comp}.

\paragraph{Attainability.} Let $\rho$ be the spectral measure supported on $[0,1]$ and satisfying our main spectral condition $\rho((0,\lambda])\le Q\lambda^\zeta$ with some $Q,\zeta>0$. Recall that the attainability condition reads $\|\mathbf w_*\|^2=\|J^{-1}\mathbf f_*\|^2=\int_0^1 \lambda^{-1}\rho(d\lambda)<\infty$. If $\zeta\le 1$, then, in general, the solution is not attainable, as can be seen by considering the exact power law $\rho((0,\lambda])= \lambda^\zeta$. On the other hand, if $\zeta>1$ then, by Lemma \ref{lem:conditions_comp}, $\mathrm{P}(\zeta,Q) \subseteq \mathrm{P}'(1,Q\tfrac{\zeta}{\zeta-\zeta'})$, implying that the solution is attainable.

\paragraph{Scaling properties.} An important property of our quadratic optimization problem is its transformation under rescaling of the input data by  $J\mapsto c J$ or by  $\mathbf f_*\mapsto c\mathbf f_*$. Under these rescalings, all the optimization algorithms of Section \ref{sec:algo} and the spectral conditions \eqref{eq:rhozeta} and \eqref{eq:lkcknu} retain their structure, but the quantities appearing in their description get rescaled by $u\mapsto c^au$ with various scaling exponents $a$. In Table \ref{tab:scaling} we list these scaling exponents. 

As an application of this observation, if we have a result for a special case when two scalar parameters are fixed, we can derive the corresponding general result by rescaling $J\mapsto cJ$ and $\mathbf f_*\mapsto c'\mathbf f_*$ with suitable $c$ and $c'$. In particular, suppose that we have a bound for $L(\mathbf w_n)$ when $\lambda_{\max}=1$ and $Q=1$. Then the corresponding bound for general $\lambda_{\max}$ and $Q$ can be obtained by taking $c=\lambda_{\max}^{1/2}$ and $c'=Q^{1/2}\lambda_{\max}^{\zeta/2}$: we see that the loss will be rescaled by $L(\mathbf w_n) \mapsto (c')^2L(\mathbf w_n)=Q\lambda_{\max}^{\zeta}L(\mathbf w_n).$

\begin{table}
    \centering
    \caption{Scaling exponents $a$ in the transformations $u\mapsto c^au$ of various quantities $u$ appearing in the descriptions of optimization algorithms (Section \ref{sec:algo}) and spectral conditions \eqref{eq:rhozeta} and \eqref{eq:lkcknu} under the transformations $J\mapsto c J$ and $\mathbf f_*\mapsto c\mathbf f_*$.}
    \label{tab:scaling}
    
    \medskip
    \begin{tabular}{lccccccccccccc}\toprule
        & $J$ & $\mathbf f_*$ & $A$ & $\mathbf b$ & $\mathbf w_n$ &  $\alpha_n$ & $\beta_n$ & $L(\mathbf w_n)$ & $Q$ & $\Lambda$ & $\lambda_{\max}$ & $\zeta$ & $\nu$\\
    \midrule
       $J\mapsto c J$  & 1 & 0 & 2 & 1 & -1 & -2 & 0 & 0 & $-2\zeta$ & 2 & 2 & 0 & 0\\
       $\mathbf f_*\mapsto c\mathbf f_*$ & 0 & 1 & 0 & 1 & 1 & 0 & 0 & 2 & 2 & 0 & 0 & 0 & 0\\
    \bottomrule
    \end{tabular}
\end{table}

\subsection{Proof of Theorem \ref{ther:worst_case_loss}}\label{sec:proof_flattened}
First, lets us examine the structure of the function $\overline{q}(x)=\sup_{y\geq x} q(y)$ introduced in Section \ref{sec:worst_case_measure}. Since $q(x)$ is a polynomial, it has a finite number of local maxima on $[0,1]$, from which we choose a maximal length sequence $0\leq x_1 < x_2 < \ldots < x_m \leq 1$ such that the values at subsequent local maxima are decreasing: $q(x_i)>q(x_{i+1})$. Then, picking $m$ points $y_i$ such that $y_i, \; i=1\ldots m-1,$ is the leftmost point in $(x_i,x_{i+1})$ satisfying $q(y_i)=q(x_{i+1})$ and $y_m=1$, allows to characterize $\overline{q}(x)$ as
\begin{equation}\label{eq:flattened_poly_description}
    \overline{q}(x)=\begin{cases}
        q(x_1), \quad &0\leq x < x_1 \\
        q(x), \quad &x_i\leq x \leq y_i, \quad i=1\ldots m\\
        q(x_{i+1}), \quad &y_i < x < x_{i+1}, \quad i=1\ldots m-1
    \end{cases}
\end{equation}
The representation \eqref{eq:flattened_poly_description} can be verified by direct comparison with the definition $\overline{q}(x)=\sup_{y\geq x} q(y)$ in each of the three cases.

Now, assume that the original polynomial $q$ is upper bounded, $q(x)\leq g(x)$, by some absolutely continuous and non-increasing $g(x)$. Then, integrating by parts, the respective ``loss'' integral can be upper-bounded as
\begin{equation}\label{eq:loss_integral_UB}
\begin{split}
    \int_0^1 q(x) \rho(dx) &\leq \int_0^1 g(x) \rho(dx) =\int_0^1 g(x) d\rho([0,x])\\
    &={}g(1)\rho([0,1]) + \int_{0}^{1} (-g'(x)) \rho([0,x])dx\\
    &\stackrel{(1)}{\leq} g(1)G(1) + \int_{0}^{1} (-g'(x)) G(x) dx =\int_0^1 g(x) G'(x) dx
\end{split}
\end{equation}
where in $(1)$ we used that $-g'(x)\geq0$ due to $g(x)$ being non-decreasing, and that $g(1)\ge 0$ since the polynomial $q$ is by assumption nonnegative. Note that $\overline{q}(x)$ given by \eqref{eq:flattened_poly_description} is absolutely continuous and non-decreasing. Thus, the bound \eqref{eq:loss_integral_UB} applies with $g(x)=\overline{q}(x)$ which sets the r.h.s of \eqref{eq:worst_case_loss_ourcond} as an upper bound for the loss integral.

Next, we show that the obtained upper bond is reached with a specific spectral measure
\begin{equation}
    \rho^* = G(x_1)\delta_{x_1} + \sum_{i=1}^{m-1} \big(G(x_{i+1})-G(x_i)\big) \delta_{x_{i+1}} + \sum_{i=1}^m \rho^*_i, \quad \rho^*_i(dx)=\mathbbm{1}_{[x_i,y_i]}G'(x)dx
\end{equation}
which is a mix of Dirac delta measures $x_i$ and ``smooth'' measures with density $\rho^*_i$, supported on $[x_i,y_i]$. Note that $\rho^*$ satisfies the required condition $\rho^*([0,x])\leq G(x)$. Direct substitution of $\rho^*$ into the loss integral gives
\begin{equation}
    \begin{split}
        \int_0^1 q(x)\rho^*(dx) &= q(x_1) G(x_1) + \sum_{i=1}^{m-1}q(x_{i+1})\int_{y_i}^{x_{i+1}}G'(x)dx + \sum_{i=1}^{m}\int_{x_i}^{y_i} q(x) G'(x)dx \\
        &=\int_0^1 \overline{q}(x) G'(x)dx
    \end{split}
\end{equation}

\subsection{Proof of Theorem \ref{ther:source_cond_suboptimality}}\label{sec:proof_th_our_vs_classical}
Our proof consists of three steps. In Step 1 we will show that a power-law asymptotic of the loss implies a power-law asymptotic of the spectral measure. Then, in Step 2 we derive the asymptotic of the bound $\widetilde L'_n,$ and in Step 3 the asymptotic of the bound $\widetilde L_n.$

\subparagraph{Step 1.} We will use the following general lemma.
\begin{lem}\label{lemmma:loss_to_measure_asym}
Suppose that $\rho$ is a Borel measure on the segment $[0,1]$, and $a>0$ is a constant. Assume that $\int_0^1(1-a\lambda)^{2n}\rho(d\lambda)=n^{-\xi}(1+o(1))$ as $n\to\infty$, with some constant $\xi>0$. Then $\rho([0,\lambda])=(\Gamma(\xi+1))^{-1}(2a\lambda)^{\xi}(1+o(1))$ as $\lambda\searrow 0$.
\end{lem}
\begin{proof}
This lemma can be derived from the general theory of abelian--tauberian power-law relations (\cite{feller1991introduction}, Section XIII.5), but we find it simpler to just give a direct proof mimicking original Karamata's arguments \citep{karamata1930certains}.  

We argue that, under the hypotheses of the lemma, for all sufficiently regular functions $g:[0,1]\to\mathbb R$ holds
\begin{equation}\label{eq:taubi}
    \lim_{n\to\infty}n^{\xi}\int_0^1(1-a\lambda)^{2n} g((1-a\lambda)^{2n}) \rho(d\lambda) = I(g), 
\end{equation}
where
\begin{equation}
    I(g)=(\Gamma(\xi+1))^{-1}\int_0^\infty e^{-y} g(e^{-y})d y^\xi.
\end{equation}
Indeed, for monomials $g(x)=x^k$ both sides of Eq. \eqref{eq:taubi} equal $(k+1)^{-\xi}$. By linearity, Eq. \eqref{eq:taubi} then holds for all polynomials. 

Now observe that the integral on the l.h.s. of Eq. \eqref{eq:taubi} is monotone in $g$ -- in the sense that if $g_1(x)\le g_2(x)$ for all $x\in[0,1]$, then the same inequality holds for the respective integrals. 

Suppose next that a function $g$ is such that for any $\epsilon>0$ one can find polynomials $g_\pm$ for which $g_-(x)\le g(x)\le g_+(x)$ on $[0,1]$ and $I(g_+)-I(g_-)<\epsilon$. Then, using the above mentioned monotonicity, Eq. \eqref{eq:taubi} holds for the function $g$, too.

Clearly, this condition holds for the function 
\begin{equation}
    g(x)=\begin{cases}
        1/x,&x\in[e^{-1},1],\\
        0,&\text{otherwise}.
    \end{cases}
\end{equation}
Substituting in Eq. \eqref{eq:taubi}, we find
\begin{equation}
    \lim_{n\to\infty} n^{\xi} \rho([0,a^{-1}(1-e^{-1/(2n)})])=(\Gamma(\xi+1))^{-1},
\end{equation}
implying the claim of the lemma.
\end{proof}

Recalling that the loss of Gradient Decent with constant learning rate $\alpha$ is given by $L_n=\frac{1}{2}\int_0^1 (1-\alpha\lambda)^{2n}\rho(d\lambda)$, the asymptotic $L_n = C n^{-\xi}(1+o(1))$ and lemma \ref{lemmma:loss_to_measure_asym} imply \begin{equation}\label{eq:qrho}
    \rho([0,\lambda])=Q_\rho\lambda^\xi(1+o(1)), \quad Q_\rho=2C \tfrac{(2\alpha)^\xi}{\Gamma(\xi+1)}.
\end{equation} 

\subparagraph{Step 2.} We use spectral asymptotic \eqref{eq:qrho} derived above to calculate the optimal upper bound $\widetilde{L}_n'$ as defined in Eq. \eqref{eq:optimal_UB_source_condition}:
\begin{align}
    \label{eq:optimal_UB_source_condition_repeat}
    \widetilde{L}'_n(\rho) ={}&\inf_{\zeta',Q':\rho\in \mathrm P'(\zeta',Q') }\sup_{\widetilde \rho \in \mathrm{P'}(\zeta',Q')}L_n(\widetilde \rho).
\end{align}
Note that Eq. \eqref{eq:sourcecond_GD_UB} already gives the supremum $L^{UB}_n(\zeta',Q')=\underset{\widetilde{\rho} \in \mathrm{P}'(\zeta',Q')}{\sup} L_n(\widetilde{\rho})$, and we only need to optimize it over $Q'$ and $\zeta'<\xi$. At a given $\zeta'$, the minimal possible $Q'$ is simply $Q'(\zeta')=\int_0^1 \lambda^{-\zeta'}\rho(d\lambda)$, so optimization reduces to that over $\zeta'$ with this $Q'(\zeta')$. Expecting the need to take $\zeta'\nearrow\xi$ at large $n$, we denote $\varepsilon=\xi-\zeta'$ and calculate  
\begin{equation}\label{eq:Qprime_calc}
\begin{split}
     \lim_{\varepsilon \searrow 0} \varepsilon Q'(\zeta')={}&\lim_{\varepsilon \searrow 0} \varepsilon\left[\rho([0,1])+\zeta'\int_0^1 \lambda^{-\zeta'-1}\rho([0,\lambda])d\lambda\right] \\
     ={}&Q_\rho\lim_{\varepsilon \searrow 0} \varepsilon \zeta'\int_0^1 \lambda^{\varepsilon-1}(1+o(1))d\lambda\\
     ={}&Q_\rho\lim_{\varepsilon \searrow 0} \zeta' \int_0^1(1+o(1))d\lambda^{\varepsilon}\\
     ={}& Q_\rho \xi,     
 \end{split}
\end{equation}
where in the first line we integrated by parts and in the second used Eq. \eqref{eq:qrho}. It follows that $Q'(\zeta')$ asymptotically behaves as
\begin{equation}\label{eq:q'zq}
    Q'(\zeta')=Q_\rho\frac{\xi}{\xi-\zeta'}(1+o(1)),\quad \zeta'\nearrow \xi.
\end{equation}
Recalling the form of the upper bound \eqref{eq:sourcecond_GD_UB}, we calculate $\widetilde{L}_n'$ as
\begin{equation}\label{eq:sourcecond_opt_bound}
    \begin{split}
        \widetilde{L}_n' ={}& \inf_{0<\zeta'<\xi}L^{UB}_n(\zeta',Q'(\zeta'))\\
        ={}& \inf_{0<\zeta'<\xi}\frac{Q'(\zeta')}{2}\left(\frac{\zeta'}{2\alpha e}\right)^{\zeta'} n^{-\zeta'}(1+o_n(1))\\
        ={}&\underset{0<\zeta'<\xi}{\inf} \; \frac{Q_\rho\xi}{2}\left(\frac{\zeta'}{2\alpha e}\right)^{\zeta'}\frac{n^{-\zeta'}}{\xi-\zeta'}(1+o_n(1))(1+o_{\zeta'}(1)) \\
        \stackrel{\zeta'=\xi-\varepsilon}{=\joinrel=}{}& \frac{Q_\rho\xi}{2}\left(\frac{\xi}{2\alpha e}\right)^{\xi} n^{-\xi}(1+o_n(1))\underset{0<\varepsilon<\xi}{\inf} \frac{n^{\varepsilon}}{\varepsilon} (1+o_\varepsilon(1)).
    \end{split}
\end{equation}
Here we added subscripts to distinguish different $o(1)$ corrections, and used that the $o_n(1)$ correction from \eqref{eq:sourcecond_GD_UB} is in fact uniform for $\zeta'\in[0,c_1]$ with any finite $c_1$. 

Recall the optimal bound $\widetilde{L}_n'(\rho_\xi)$ for exact power-law measure given in \eqref{eq:optima_source_exact_powerlaw}. Substitution of the infimum $\inf_{0<\varepsilon<\xi} \frac{n^\varepsilon}{\varepsilon}=e \log n, \; n>e^{\frac{1}{\xi}}$ and the expression for $Q_\rho$ into \eqref{eq:optima_source_exact_powerlaw} gives the desired statement \eqref{eq:optima_source} of the theorem.

However, we still need to argue that this result is not affected by the factor $1+o_\varepsilon(1)$ appearing in $\inf_{0<\varepsilon<\xi}\tfrac{n^{\varepsilon}}{\varepsilon} (1+o_\varepsilon(1))$. To this end, it clearly suffices to show that the optimal $\varepsilon\to 0$ as $n\to\infty$. By tracing back our expression $1+o_\varepsilon(1)$ to formula  \eqref{eq:q'zq}, this expression is bounded away from 0 on the interval $[0,\xi]$. Then, on any interval $[c,\xi]$ with $c>0$ we get a power-law lower bound
\begin{equation}
    \inf_{c\le \varepsilon<\xi}\tfrac{n^{\varepsilon}}{\varepsilon} (1+o_\varepsilon(1))=\Omega(n^c),\quad n\to\infty.
\end{equation}
This shows by comparison with the logarithmic expression $\inf_{0<\varepsilon<\xi} \frac{n^\varepsilon}{\varepsilon} =e \log n$ that the values $\varepsilon$ bounded away from 0 are indeed asymptotically suboptimal. This completes the computation of $\widetilde{L}_n'$.

\subparagraph{Step 3.} Finally, we calculate the optimal bound $\widetilde{L}_n$ under our source condition \eqref{eq:rhozeta1}, as defined in Eq. \eqref{eq:optimal_UB_our_condition}:
\begin{align}
    \label{eq:optimal_UB_our_condition_repeat}
    \widetilde{L}_n(\rho) ={}& \inf_{\zeta,Q:\rho\in \mathrm P(\zeta,Q)}\sup_{\widetilde \rho \in \mathrm{P}(\zeta,Q)}L_n(\widetilde \rho). 
\end{align}
First, recall that the inner supremum here is given by theorem \ref{ther:worst_case_loss}, where for GD with $\alpha \leq 1$ the flattened polynomial $\overline{(1-\alpha\lambda)^{2n}}=(1-\alpha\lambda)^{2n}$. Thus, we have
\begin{equation}\label{eq:ourcond_GD_upperbound}
\begin{split}
        \underset{\widetilde{\rho}\in P(\zeta,Q)}{\sup} L_n(&\widetilde{\rho}) = \frac{Q}{2}\int_0^1 (1-\alpha\lambda)^{2n}d(\lambda^\zeta) = \frac{Q}{2}\int_0^{\alpha^{-1}} (1-\alpha\lambda)^{2n}d(\lambda^\zeta) + O((1-\alpha)^{2n})\\
        &=\frac{Q}{2}\zeta\alpha^{-\zeta} \frac{\Gamma(2n+1)\Gamma(\zeta)}{\Gamma(2n+\zeta+1)} +O((1-\alpha)^{2n})\overset{n\to\infty}{=\joinrel=} \frac{Q}{2}\Gamma(\zeta+1)(2\alpha n)^{-\zeta} (1+o(1)),
\end{split}
\end{equation}
where we recognized the integral $\int_0^1(1-z)^{2n}z^{\zeta-1}dz$ as a Beta function and substituted its expression in terms of Gamma functions.\footnote{Actually, the same computation is performed in the proof of the theorem \ref{ther:constant_lr_bounds}, see eq. \eqref{eq:const_lr_exact_powerlaw_measure_loss}. We repeat it here simply for convenience.} 

To optimize this expression over $Q$ and $\zeta$, note that we can take any $\zeta\leq\xi$, and at the given $\zeta$ the minimal constant $Q$ is 
\begin{equation}\label{eq:ourcond_optimal_costant}
    Q(\zeta)=\underset{\lambda\in(0,1]}{\sup} \rho([0,\lambda])/\lambda^\zeta = \underset{\lambda\in(0,1]}{\sup} Q_\rho\lambda^{\xi-\zeta}(1+o(1)).
\end{equation}
We note a couple of properties of $Q(\zeta)$:
\begin{enumerate}
    \item $Q(\zeta)\nearrow Q(\xi)$ as $\zeta \nearrow \xi$, because for any $\lambda\in(0,1]$ the function $\zeta\mapsto \rho([0,\lambda])/\lambda^\zeta$ is monotone non-decreasing and converging to $\rho([0,\lambda])/\lambda^\xi$ as $\zeta  \nearrow \xi$.
    \item $Q(\zeta)$ is bounded away from 0 on the interval $0\le \zeta\le \xi,$ since $Q(\zeta)\ge \rho([0,1])>0$.
\end{enumerate}
Property 2) and representation \eqref{eq:ourcond_optimal_costant} imply that the infimum of $\underset{\widetilde{\rho}\in P(\zeta,Q)}{\sup} L_n(\widetilde{\rho})$ over $\zeta$ and $Q(\zeta)$ is attained at a $\zeta$ deviating from $\xi$ by at most $O(1/\log n);$ in particular the optimal $\zeta$ converges to $\xi$ as $n\to\infty$. But then, using property 1) we get the desired asymptotic  \eqref{eq:optima}:
\begin{equation}
    \widetilde{L}_n(\rho) = \frac{Q(\xi)}{2}\Gamma(\xi+1)(2\alpha n)^{-\xi} (1+o(1)).
\end{equation}
This completes the proof of the theorem.

\section{Constant learning rates}\label{sec:constant_lr}
\subsection{Proof of Theorem \ref{ther:constant_lr_bounds}: the case of GD ($\beta=0$)}
First, we express worst-case loss \eqref{eq:worst_case_loss} through exact power-law loss \eqref{eq:exact_power_measure_loss}. 

If $\alpha\leq 1$, the polynomial $p_n^2(\lambda)=(1-\alpha\lambda)^{2n}$ is monotone decreasing and therefore $\overline{p_n^2}(\lambda)=p_n^2(\lambda)$. This implies that $\overline{L_n^{(\zeta)}}=L_n^{(\zeta)}$. If $1<\alpha<2$, the flattened polynomials $\overline{p_n^2}(\lambda)$ differ from $p_n^2(\lambda)$ on a single flat region and are given by
\begin{equation}
    \overline{p_n^2}(\lambda) = 
    \begin{cases}
    p_n^2(\lambda), \quad &\lambda<\frac{2-\alpha}{\alpha}\\
    (\alpha-1)^{2n}, \quad &\frac{2-\alpha}{\alpha}\leq\lambda\leq1
    \end{cases}
\end{equation}
The associated worst-case loss is 
\begin{equation}
\begin{split}
     \overline{L_n^{(\zeta)}}&=\frac{1}{2}\int_0^{\frac{2-\alpha}{\alpha}}p_n^2(\lambda)d(\lambda^\zeta)+\frac{1}{2}\int_{\frac{2-\alpha}{\alpha}}^1(\alpha-1)^{2n}d(\lambda^\zeta)\\
     &=\frac{1}{2}\int_0^{\frac{2-\alpha}{\alpha}}p_n^2(\lambda)d(\lambda^\zeta)+O\big((\alpha-1)^{2n}\big)\\
     &=\frac{1}{2}\int_0^{1}p_n^2(\lambda)d(\lambda^\zeta)+O\big((\alpha-1)^{2n}\big), 
\end{split}
\end{equation}
which is exactly the $\beta=0$ part of \eqref{eq:worst_case_loss_constant_lr_GD} with $u=(1-\alpha)^2$. Finally, we calculate the loss under exact power-law measure as
\begin{equation}\label{eq:const_lr_exact_powerlaw_measure_loss}
    \begin{split}
        L_n^{(\zeta)} &= \frac{1}{2} \int_0^1 (1-\alpha\lambda)^{2n}d(\lambda^\zeta) =  \frac{1}{2} \zeta\int_0^{\frac{1}{\alpha}}(1-\alpha\lambda)^{2n}\lambda^{\zeta-1}d\lambda + O\big((\alpha-1)^{2n}\big) \\
        &=\frac{1}{2}\zeta\alpha^{-\zeta}\frac{\Gamma(2n+1)\Gamma(\zeta)}{\Gamma(2n+1+\zeta)}+O\big((\alpha-1)^{2n}\big) \\
        &= \frac{1}{2} \Gamma(\zeta+1) (2n\alpha)^{-\zeta}(1+o(1)) + O\big((\alpha-1)^{2n}\big)
    \end{split}
\end{equation}
Here in the second line, we recognized the integral representation of the Beta function $B(a,b)=\int_0^1(1-z)^{a-1}z^{b-1}dz$ and expressed it through the Gamma functions. In the last line, we used $x\Gamma(x)=\Gamma(x+1)$ and and asymptotic of Gamma function $\Gamma(x+a)=\Gamma(x)x^a(1+o(1))$.

\subsection{Proof of Theorem \ref{ther:constant_lr_bounds}: the case of HB ($\beta\ne 0$)}\label{sec:const_lr_HB_UB}
\paragraph{Structure of HB residual polynomials.}
We start with deriving expression for residual polynomial corresponding to HB method with step-size $\alpha$ and momentum $\beta$. These residual polynomials satisfy recurrence relation with constant coefficients
\begin{equation}\label{eq:HB_recursion}
    p_{n+1}(\lambda) = p_n(\lambda) - \alpha \lambda p_n(\lambda) + \beta (p_n(\lambda)-p_{n-1}(\lambda)), \qquad p_0(\lambda)=p_{-1}(\lambda)=1.
\end{equation}
Linear transformations of the polynomials $p_n(\lambda) = c^n q_n(z), z=ax+b$ lead to new polynomials $q_n$ with different constants in their recurrence relations, which we choose to be that of Chebyshev polynomials.
\begin{align}
    \label{eq:HB_residual_poly_through_Chebyshev}
    p_n(\lambda) &= (\sqrt{\beta})^n q_n(z(\lambda)), \quad z(\lambda) = \frac{1-\alpha\lambda+\beta}{2\sqrt{\beta}} \\ 
    \label{eq:Chebyshev_recurrence}
    q_{n+1}(z) &= 2zq_n(z) - q_{n-1}(z), \quad q_0(z)=1, \quad q_{-1}(z)=\sqrt{\beta}
\end{align}
The initial conditions in \eqref{eq:Chebyshev_recurrence} are satisfied with $q_n(z)=U_n(z)-\sqrt{\beta}U_{n-1}(z)$, where $U_n(z)$ are Chebyshev polynomials of second kind 
\begin{equation}
    \label{eq:Chebyshev_polynomials_2}
    U_n(z) = \begin{cases}
    \frac{\sin((n+1) \varphi)}{\sin\varphi}, \; &|z|\leq 1\\
    \frac{\left((z+\sqrt{z^2-1})^{n+1}-(z-\sqrt{z^2-1})^{n+1}\right)}{2\sqrt{z^2-1}}, \; &|z| \geq 1
    \end{cases}
\end{equation}
Here $\cos\varphi = z$. Thus, we derived representation \eqref{eq:const_lr_HB_poly} for HB residual polynomials.

Let's list properties of $q_n(z)$ which will be useful in the subsequent parts of the proof.
\begin{enumerate}
    \item \textit{Monotonocity w.r.t. $z$}: \\ $q_n(z)^2$ is monotone decreasing for $z\in(-\infty, -1]$\ and monotone increasing for $z\in[1,\infty)$.
    \item \textit{Monotonocity w.r.t. $n$}: 
    \begin{equation}\label{eq:HB_polynomials_n_monotonicity}
        (\sqrt{\beta})^{n+1} q_{n+1}(z)\leq (\sqrt{\beta})^{n} q_{n}(z) \quad \text{for} \; z\in[1,\frac{1+\beta}{2\sqrt{\beta}}]
    \end{equation}
\end{enumerate}
The first property follows from the fact that all $n-1$ zeros of the derivative $\tfrac{d}{dz}q_n(z)$ are located between $n$ roots of $q_n(z)$, which in turn  are located on $(-1,1)$. To get the latter, note that the zero of $q_n(z)=U_n(z)-\sqrt{\beta}U_{n-1}(z)=\Big(\sin((n+1)\varphi) - \sqrt{\beta}\sin n\varphi\Big)/\sin \varphi$ is equivalent to
\begin{equation}\label{eq:zeros_of_HB_polynomial}
\begin{cases}
    \tan n\varphi=-\frac{\sin \varphi}{\cos \varphi-\sqrt{\beta}}, \quad & \cos \varphi \ne \sqrt{\beta} \\
    \cos n\varphi=0, \quad & \cos \varphi=\sqrt{\beta}
\end{cases}
\end{equation}
Here the first equation has at least $n-1$ solutions: a single solution on each interval $\tfrac{\pi}{2}+\pi k < n\varphi < \tfrac{\pi}{2}+\pi (k+1), \quad k=0,\ldots,n-2$. The remaining solution can be found in the interval containing $\cos \varphi =\sqrt{\beta}$, or exactly on the boundary if the second equation in \eqref{eq:zeros_of_HB_polynomial} is satisfied.

The obtain the second property, note that it is equivalent to $r_n(z)\leq 1$ where $r_{n+1}(z)=\tfrac{\sqrt{\beta}q_{n+1}(z)}{q_n(z)}$ and satisfies $r_{n+1}(z)=2\sqrt{\beta}z-\tfrac{\beta}{r_n(z)}$ due to \eqref{eq:Chebyshev_recurrence}. Observing that $r_0(z)=1$ we proceed by induction and assume that $r_n(z)\leq 1$ for $z\in[1,\frac{1+\beta}{2\sqrt{\beta}}]$. Then, using that all $q_n(z)$ are positive for $z\geq 1$ and therefore $r_n(z)>0$, we get $r_{n+1}\leq 2\sqrt{\beta}z-\beta\leq1$ for $z\in[1,\frac{1+\beta}{2\sqrt{\beta}}]$. 

\paragraph{Bounding the worst-case loss.} 
First, let's bound $q_n(z)$ inside the oscillatory region $z\in[-1,1]$. Since $|U_n(z)|\leq n+1$ for $z\in[-1,1]$, we get $|q_n(z)|=|U_n(z)-\sqrt{\beta}U_{n-1}(z)|\leq 2n+1$.  

Next, we bound $q_n(z)$ to the left of oscillatory region: $z<-1$. For convenience, we denote $z_{\pm}=z\pm \sqrt{z^2-1}$, and write
\begin{equation}\label{eq:HB_Chebyshev_explicit_form}
\begin{split}
    q_n(z) &= \frac{1}{2\sqrt{z^2-1}}\left[z_+^n\Big(z+\sqrt{z^2-1}-\sqrt{\beta}\Big)+z_-^n\Big(-z+\sqrt{z^2-1}+\sqrt{\beta}\Big)\right]\\
    &=\frac{z_+^n+z_-^n}{2}+(z-\sqrt{\beta})\sum_{k=0}^{n-1}z_+^kz_-^{n-1-k}
\end{split}
\end{equation}
using the representation above and the fact that $|z_-|\geq|z_+|$ for $z<-1$, we get
\begin{equation}
\begin{split}
    |q_n(z)|&\leq \frac{|z_+|^n+|z_-|^n}{2}+(\sqrt{\beta}-z)\sum_{k=0}^{n-1}|z_+|^k|z_-|^{n-1-k} \\
    &\leq |z_-|^n\left(1+\frac{n(\sqrt{\beta}-z)}{|z_-|}\right) \leq (2n+1)|z_-|^n
\end{split}
\end{equation}
Now, we are ready to bound the flattened HB polynomial 
\begin{equation}
    \overline{p_n^2}(\lambda)=(\sqrt{\beta})^n \overline{q_n^2}\big(z(\lambda)\big), \quad \overline{q_n^2}(z)=\underset{z_1\leq y \leq z}{\sup} q_n^2(y)
\end{equation}
where $z_1=z(\lambda=1)=\tfrac{1-\alpha+\beta}{2\sqrt{\beta}}$. Now, recall the monotonicity properties of $q_n^2(z)$ on $(-\infty,-1]$ and $[1,\infty)$. Then, for $z_1\geq1$ we immediately get $\overline{q_n^2}(z)=z_n^2(z)$, while for $z_1<1$ we first get a single bound on $[z_1,1]$ as 
\begin{equation}\label{eq:HB_polynomial_exp_bound}
    |q_n^2(z)| \leq \operatorname{max}\Big\{2n+1, \; \mathbbm{1}_{z_1<-1}(2n+1)|z_-(z_1)|^n \Big\} = (2n+1) \left|z_1-\sqrt{z_1^2-1}\right|^n
\end{equation}
where for $z_1\in (-1,1)$ the square root $\sqrt{z_1^2-1}$ is understood in the complex sense. 

Combining the obtained bounds, we can compactly characterize the flattened polynomial as
\begin{equation}
    \overline{q_n^2}(z)= q_n^2(z)+\mathbbm{1}_{z_1<1} O\Big(n^2\left|z_1-\sqrt{z_1^2-1}\right|^{2n}\Big),
\end{equation}
implying for the worst-case loss
\begin{equation}
        \overline{L_n^{(\zeta)}}=\frac{1}{2}\int_0^1 p_n^2(\lambda)d(\lambda^\zeta) + \int_0^1 \mathbbm{1}_{z_1<1} O\Big(n^2\beta^n\left|z_1-\sqrt{z_1^2-1}\right|^{2n}\Big) d(\lambda^\zeta)
\end{equation}
which is exactly the momentum case of \eqref{eq:worst_case_loss_constant_lr_GD} with $u=\beta \left|z_1-\sqrt{z_1^2-1}\right|^2$.

\paragraph{Calculating the loss under the exact power-law measure.}
While this can be done in a number of ways, we choose the approach based on the generating functions of $p_n^2(\lambda)$ and $L_n$. The approach is based on the connection between the asymptotic of the loss $L_n^{(\zeta)}$ and the singularity of its generating function 
\begin{equation}\label{eq:loss_gen_function}
    G_L(t)=\sum_{n=0}^\infty t^n L_n^{(\zeta)}
\end{equation}
at $t=1$. The two are connected by Tauberian theorem (\citet{feller1991introduction}, p. 445) which states that if generating function $G(t)=\sum_n t^n a_n$ of a sequence $a_n$ has asymptotic $G(1-\varepsilon)=C\varepsilon^{-\rho}(1+o(1)), \; \rho>0$, then \begin{equation}
    \sum_{k=1}^n a_k=\frac{C}{\Gamma(\rho+1)}n^\rho(1+o(1)),\quad n\to\infty.
\end{equation}  We will apply this theorem to the sequence $a_n=n^m L_n^{(\zeta)}$, where $m= \lfloor \zeta \rfloor$ is required to get a divergent behavior of the partial sums.

First, recall that thanks to \eqref{eq:HB_residual_poly_through_Chebyshev}, \eqref{eq:HB_Chebyshev_explicit_form} we can write HB residual polynomials in the form $p_n(\lambda) = f_+(z(\lambda))\big(\sqrt{\beta} z_+(\lambda)\big)^n + f_-(z(\lambda))\big(\sqrt{\beta} z_-(\lambda)\big)^n$. Then, generating function of $p_n^2(\lambda)$ can be immediately written as
\begin{equation}
\begin{split}
    G_p(t,\lambda)&\equiv\sum_{n=0}^\infty t^n p_n^2(\lambda) = \sum_{n=0}^\infty \Big[f_+^2(t\beta z_+^2)^n + f_-^2(t\beta z_-^2)^n + 2f_+f_-(t\beta z_+z_-)^n\Big] \\
    &=\frac{f_+^2}{1-t\beta z_+^2}+\frac{f_-^2}{1-t\beta z_-^2} + \frac{2f_+f_-}{1-t\beta z_+z_-}
\end{split}
\end{equation}
Substituting $z(\lambda)$ into $f_+(z),f_-(z),z_+(z),z_-(z)$ and straightforwardly simplifying the expression (e.g., using symbolic computer algebra software) reveals that $G_p(t,\lambda)$ is a rational function of its arguments equal to
\begin{align}
    \label{eq:p_n_generating_function}
    G_p(t,\lambda) &= \frac{(1-\beta t)(1-\beta^2t)+2\alpha\beta\lambda t}{(1-\beta t)\Big((1-t)(1-\beta^2 t)+\alpha\lambda t(2+2\beta-\alpha\lambda t)\Big)}\\
    \label{eq:p_n_generating_function_asym}
    &\stackrel{t=1-\varepsilon}{=}\frac{1}{\varepsilon+\tfrac{2\alpha}{1-\beta}\lambda}\Big(1+O(\varepsilon)+O(\lambda)\Big), \quad \text{as} \quad \varepsilon \searrow 0 \quad \text{and} \quad \lambda \searrow 0 
\end{align}
Here we observed from \eqref{eq:p_n_generating_function} that when stability condition $\alpha<2(1+\beta)$ is satisfied, $G_p(t,\lambda)$ on $[0,1]^2$ is regular everywhere except the singularity at $t=1$, $\lambda=0$.

Focusing on the contribution to the loss $L_{a,n}^{(\zeta)}=\frac{1}{2}\int_0^a p_n^2(\lambda)d(\lambda^\zeta)$ from $[0,a], \; a\leq1$ (to be specified later), we write $m$-th derivative of its generating function $G_{L,a}(t)=\sum_{n=0}^\infty t^n L_{a,n}^{(\zeta)}$ as
\begin{equation}
    \begin{split}
        \Big(t\frac{d}{dt}\Big)^m G_{L,a}(t) &=\frac{1}{2} \int_0^a \Big(t\frac{\partial}{\partial t}\Big)^m G_p(t,\lambda) d(\lambda^\zeta)\\
        &=\frac{m!}{2}\int_0^a \frac{\zeta \lambda^{\zeta-1}}{\left(\varepsilon + \frac{2\alpha}{1-\beta}\lambda\right)^{m+1}}\big(1+O(\varepsilon)+O(\lambda)\big)d\lambda\\
        &=\frac{m!\zeta}{2}\left(\frac{2\alpha}{1-\beta}\right)^{-\zeta}\varepsilon^{\zeta-m-1}\int_0^{\tfrac{2\alpha a}{\varepsilon(1-\beta)}}\frac{x^{\zeta-1}}{(1+x)^{m+1}}\big(1+(1+x)O(\varepsilon)\big)dx\\
        &=\frac{m!\zeta}{2}\left(\frac{2\alpha}{1-\beta}\right)^{-\zeta}\varepsilon^{\zeta-m-1}(1+O(\varepsilon))\int_0^\infty\frac{x^{\zeta-1}dx}{(1+x)^{m+1}}\\
        &=\frac{\Gamma(\zeta+1)\Gamma(m+1-\zeta)}{2}\left(\frac{2\alpha}{1-\beta}\right)^{-\zeta}\varepsilon^{\zeta-m-1}(1+O(\varepsilon))
    \end{split}
\end{equation}
where in the second-to-last line, we recognized the integral representation of Beta function $B(\zeta,m+1-\zeta)$ and subsequently expressed it in terms of Gamma functions. Observing that $\big(t\tfrac{d}{dt}\big)^m G_{L,a}(t)$ is the generating function of the sequence $n^m L_{a,n}^{(\zeta)}$, we apply Tauberian theorem to get asymptotic of the partial sums 
\begin{equation}\label{eq:partial_sum_asym}
    \sum_{k=0}^n k^m L_{a,k}^{(\zeta)} = \frac{\Gamma(\zeta+1)}{2(m+1-\zeta)}\left(\frac{2\alpha}{1-\beta}\right)^{-\zeta}n^{m+1-\zeta} (1+o(1))
\end{equation}
Now, we choose $a=\operatorname{min}\{1,(1-\sqrt{\beta})^2/\alpha\}$ where the second option corresponds to the border of the oscillating region $z(a)=1$ of polynomials $p_n(\lambda)$. Then, the monotonicity property \eqref{eq:HB_polynomials_n_monotonicity} imply monotonicity of $p_n^2(\lambda)$ on $[0,a]$, and therefore monotonicity of $L_{a,k}^{(\zeta)}$. This enables to use Lemma \ref{lem:partial_sums} below on partial sums \eqref{eq:partial_sum_asym} and get $L_{a,n}^{(\zeta)}=\tfrac{\Gamma(\zeta+1)}{2}(\tfrac{2\alpha n}{1-\beta})^{-\zeta}(1+o(1))$, which is the same as \eqref{eq:exact_power_measure_loss_const_lr_GD} thanks to exponentially suppressed (see eq. \eqref{eq:HB_polynomial_exp_bound}) contribution to the loss from $\lambda\in[a,1]$. 

\begin{lem}\label{lem:partial_sums}
Assume a sequence $a_n$ is monotonically decreasing, and there is $m\geq 0$ such that $\sum_{k=1}^n k^m a_k = n^\xi(1+o(1))$ with some $\xi>0$. Then, $a_n =\xi n^{\xi-m-1}(1+o(1))$.
\end{lem}
\begin{proof}
    Take a fixed $r>0$ and consider the partial sums $S_n = \sum_{k=n}^{n'-1}k^m a_k$ in the chunks $[n, n')$, $n'=\lfloor n(1+r)\rfloor$. In the limit $n\to\infty$ we have
    \begin{align}
        S_n &=  n^\xi \big[(1+r)^\xi-1\big](1+o(1))\\
        S_n &\leq a_n \sum_{k=n}^{n'-1}k^m = a_n n^{m+1}\frac{(1+r)^{m+1}-1}{m+1}(1+o(1))
    \end{align}
Combining these two estimates yields the bound
\begin{equation}\label{eq:partial_sums_lemma_bound}
    a_n \geq n^{\xi-m-1} \frac{(m+1)[(1+r)^\xi-1]}{(1+r)^{m+1}-1} (1+o(1))
\end{equation}
As $r$ was arbitrary, we take $r \searrow 0$ in \eqref{eq:partial_sums_lemma_bound} and get $a_n \geq \xi n^{\xi-m-1} (1+o(1))$. Next, we take a fixed $0<r<1$ and consider the partial sums in the chunks $[n', n)$, $n'=\lfloor n(1-r)\rfloor$. Then, similar reasoning gives $a_n \leq \xi n^{\xi-m-1} (1+o(1))$, thus completing the proof.
\end{proof}

\subsection{Proof of Theorem \ref{ther:const_lr_discrete_lower_bound}}
First, observe that the measure $\rho_{\zeta,\nu}$ trivially satisfies the condition \eqref{eq:lkcknu} since the eigenvalues corresponding to $\rho_{\zeta,\nu}$ are $\lambda_k=k^{-\nu}$. Next, we take  $\lambda\in[\lambda_k, \lambda_{k-1})$ and evaluate the respective cumulative distribution function $\rho_{\zeta,\nu}([0,\lambda])$ as 
\begin{equation}
    \rho_{\zeta,\nu}([0,\lambda]) = \sum_{l\geq k}^\infty \Big(l^{-\zeta\nu}-(l+1)^{-\zeta\nu}\Big) = k^{-\zeta\nu}=\lambda_k^{\zeta} \leq \lambda^\zeta,
\end{equation}
which confirms that $\rho_{\zeta,\nu}$ satisfies the main condition \eqref{eq:rhozeta}. 

To bound the loss under the measure $\rho_{\zeta.\nu}$, we first do so for $\rho_{\zeta,\nu}([0,\lambda])$. Take a $k_0\geq0$ such that $G_{\zeta,\nu}(\lambda)=\lambda^\zeta-\zeta\nu\lambda^{\zeta+1/\nu}$ is increasing on $[0,\lambda_{k_0}]$ and consider again $\lambda \in (\lambda_{k+1},\lambda_k], \; k\geq k_0$:
\begin{equation}
    G_{\zeta,\nu}(\lambda)\leq G_{\zeta,\nu}(\lambda_k) = k^{-\zeta\nu} - \zeta\nu k^{-\zeta\nu-1} \leq (k+1)^{-\zeta\nu}=\lambda_{k+1}^\zeta\leq\rho([0,\lambda]) 
\end{equation}
Thus, we established that $\rho_{\zeta,\nu}([0,\lambda])\geq G_{\zeta,\nu}(\lambda)$ for $\lambda \in [0,\lambda_{k_0}]$. Now, let $p_n(\lambda)$ be the residual polynomial of the considered GD algorithm and $\lambda_0$ be it's left-most zero. Since $p_n(\lambda)$ is monotone decreasing on $[0,\lambda_0]$ (see the proof of Theorem \ref{ther:constant_lr_bounds}), the contribution to the loss from $[0,\lambda^*], \; \lambda^*=\operatorname{min}(\lambda_{k_0},\lambda_0)$ is given by
\begin{equation}
\begin{split}
    \int_0^{\lambda^*}p_n^2&(\lambda)\rho_{\zeta,\nu}(d\lambda) =p_n^2(\lambda^*)\rho_{\zeta,\nu}([0,\lambda^*]) - \int_0^{\lambda^*} \Big(\frac{d}{d\lambda}p_n^2(\lambda)\Big)\rho_{\zeta,\nu}([0,\lambda])d\lambda\\
    &\geq p_n^2(\lambda^*)G_{\zeta,\nu}(\lambda^*) - \int_0^{\lambda^*} \Big(\frac{d}{d\lambda}p_n^2(\lambda)\Big)G_{\zeta,\nu}(\lambda)d\lambda=\int_0^{\lambda^*}p_n^2(\lambda)G_{\zeta,\nu}'(\lambda)d\lambda
\end{split}
\end{equation}
Referring to the proof of Theorem \eqref{ther:constant_lr_bounds} and eq. \eqref{eq:HB_Chebyshev_explicit_form} we observe that on any $[a,1], \; a>0$ the residual polynomials decay uniformly as $p_n^2(\lambda) = O(r^n), \; r<1$. Using this and \eqref{eq:exact_power_measure_loss_const_lr_GD} we bound the loss as
\begin{equation}
    \begin{split}
        L_n &= \frac{1}{2}\int_0^{1} p_n^2(\lambda)\rho_{\zeta,\nu}(d\lambda)= \frac{1}{2}\int_0^{\lambda^*} p_n^2(\lambda)\rho_{\zeta,\nu}(d\lambda) + O(r^n)\\
        &\geq \frac{1}{2}\int_0^{\lambda^*} p_n^2(\lambda)G_{\zeta,\nu}'(\lambda)d\lambda + O(r^n) = \frac{1}{2}\int_0^{1} p_n^2(\lambda)\Big[d(\lambda^\zeta)-\zeta\nu d(\lambda^{\zeta+1/\nu})\Big] + O(r^n)\\
        &=\frac{1}{2}\int_0^{1} p_n^2(\lambda)d(\lambda^\zeta)+O(n^{-\zeta-1/\nu}) + O(r^n) = \frac{\Gamma(\zeta+1)}{2}\left(\frac{2\alpha n}{1-\beta}\right)^{-\zeta}(1+o(1))
    \end{split}
\end{equation}

\section{Accelerated methods for exact power-law spectral measure}\label{sec:lbcgjacobi}
\paragraph{Proof of Theorem \ref{th:cgexactpower}.}
We substitute CG residual polynomial given by \eqref{eq:power-law_measure_solution} into the loss \eqref{eq:loss_through_residual_poly}
\begin{equation}\label{eq:pnviaPn2}
    p_n(\lambda)=\frac{P_n^{(\zeta,0)}(1-2\lambda)}{P_n^{(\zeta,0)}(1)}.
\end{equation}
Then, by a change of variables,
\begin{equation}
    L(\mathbf w_n)=\frac{1}{2}\int_0^1 p_n^2(\lambda)d\lambda^\zeta=\frac{\zeta }{2^{\zeta+1}\big(P_n^{(\zeta,0)}(1)\big)^2}\int_{-1}^1 (1-x)^{\zeta-1}(P_n^{(\zeta,0)}(x))^2dx.
\end{equation}
We will use Rodrigues' formula for $P_n^{(a,b)}:$
\begin{equation}
    P_n^{(a,b)}(x)=\frac{(-1)^n}{2^n n!}(1-x)^{-a}(1+x)^{-b}\frac{d^n}{dx^n}[(1-x)^{a+n}(1+x)^{b+n}].
\end{equation}
It gives (with $a=\zeta, b=0$)
\begin{equation}\label{eq:pnint}
    \int_{-1}^1 (1-x)^{\zeta-1}(P_n^{(\zeta,0)}(x))^2dx = \frac{(-1)^n}{2^n n!}\int_{-1}^1 (1-x)^{-1}\frac{d^n}{dx^n}[(1-x)^{\zeta+n}(1+x)^{n}] P_n^{(\zeta,0)}(x)dx.
\end{equation}
Observe that we can write 
\begin{equation}\label{eq:1xpn}
    (1-x)^{-1}P_n^{(\zeta,0)}(x)=P_n^{(\zeta,0)}(1)(1-x)^{-1}+q_{n-1}(x)
\end{equation}
with some polynomial $q_{n-1}$ of degree $n-1$. Suppose that we perform repeated integration by parts in the r.h.s. of \eqref{eq:pnint}, moving all the derivatives $\tfrac{d^n}{dx^n}$ from $(1-x)^{\zeta+n}(1+x)^{n}$ to $(1-x)^{-1}P_n^{(\zeta,0)}(x)$. Thanks to the condition $\zeta>0$, all the boundary terms will vanish. Moreover, since $d^n q_{n-1}/dx^n=0$, only the first term in the r.h.s. of Eq. \eqref{eq:1xpn} will give a nonvanishing contribution to the resulting integral, specifically
\begin{equation}
\begin{split}
    \int_{-1}^1 (1-x)^{-1}&\frac{d^n}{dx^n}[(1-x)^{\zeta+n}(1+x)^{n}] P_n^{(\zeta,0)}(x)dx \\
    &= (-1)^n n!P_n^{(\zeta,0)}(1)\int_{-1}^1(1-x)^{-n-1} [(1-x)^{\zeta+n}(1+x)^{n}] dx\\
    &=(-1)^n n!P_n^{(\zeta,0)}(1)2^{\zeta+n} \int_0^1 t^{\zeta-1}(1-t)^{n} dt\\
    &=(-1)^n n!P_n^{(\zeta,0)}(1)2^{\zeta+n} B(\zeta, n+1).
\end{split}
\end{equation}
Using the fact that
\begin{equation}
    P_n^{(\zeta,0)}(1) = \frac{\Gamma(\zeta+n+1)}{n!\Gamma(\zeta+1)}, 
\end{equation}
we finally obtain
\begin{align}
    L(\mathbf w_n)={}&\frac{\zeta }{2^{\zeta+1}\big(P_n^{(\zeta,0)}(1)\big)^2}\cdot \frac{(-1)^n}{2^n n!}\cdot (-1)^n n!P_n^{(\zeta,0)}(1)2^{\zeta+n}\cdot\frac{\Gamma(\zeta)\Gamma(n+1)}{\Gamma(\zeta+n+1)}\\
    ={}&\frac{\Gamma^2(\zeta+1)n!^2}{2\Gamma^2(\zeta+n+1)}\\
    ={}&\frac{ \Gamma^2(\zeta+1)}{2}n^{-2\zeta}(1+o(1))\quad (\lambda\to 0+).
\end{align}

\paragraph{Proof of Proposition \ref{prop:Jacobi_anzatz_typical_convergence}.}

The principal difference between $a>\zeta-\frac{1}{2}$ and $a<\zeta-\frac{1}{2}$ is that in the former case the dominating contribution to the integral comes from $\lambda \sim n^{-2}$ while in the latter case the dominant contribution comes from $\lambda\sim1$.

Let's start with $a>\zeta-\frac{1}{2}$. The classical asymptotic expansion of Jacobi polynomials $P^{(a,b)}_n(\cos\theta)$ at small $\theta$ (\cite{szego1959orthogonal}, Theorem 8.21.12.) states for a fixed $c,\varepsilon$ and $N=n+\frac{1}{2}(a+b+1)$
\begin{equation}\label{eq:Jacobi_classical_asymptotic}
\begin{split}
        \sin\left(\frac{\theta}{2}\right)^a \cos\left(\frac{\theta}{2}\right)^b P^{(a,b)}_n(\cos\theta) =  &N^{-a}\frac{\Gamma(n+a+1)}{n!}\left(\frac{\theta}{\sin\theta}\right)^{\frac{1}{2}} J_a(N\theta) \\
        &+ \begin{cases}
    \theta^{a+2}O(n^a), \quad & \theta<cn^{-1} \\
    \theta^\frac{1}{2}O(n^{-\frac{3}{2}}), \quad & cn^{-1}<\theta<\pi-\varepsilon
    \end{cases}
\end{split}
\end{equation}
Using that $z^{-a}J_a(z)$ is bounded and also $|J_a(z)| = O(z^{-\frac{1}{2}})$ uniformly, we adopt \eqref{eq:Jacobi_classical_asymptotic} to our needs and write an asymptotic form
\begin{equation}\label{eq:Jacobi_square_asym}
\begin{split}
    \Big(P^{(a,b)}_n&(\cos\theta)\Big/P^{(a,b)}_n(1) \Big)^2 = \\
    =&\frac{\Gamma^2(a+1)}{2N^{-2a-1}}\frac{N\theta J_a^2(N\theta)}{\big(\sin\tfrac{\theta}{2}\big)^{2a+1}\big(\cos\tfrac{\theta}{2}\big)^{2b+1}}+ \begin{cases}
    \theta^2 O(1),&\theta<cn^{-1} \\
    \theta^{-2a}O(n^{-2a-2}),&cn^{-1}<\theta<\pi-\varepsilon
    \end{cases}\\
    =& \Gamma^2(a+1)\left(\frac{N\theta}{2}\right)^{-2a}J_a^2(N\theta) + \begin{cases}
    \theta O(1),&\theta<cn^{-1} \\
    \theta^{-2a}O(n^{-2a-1}),&cn^{-1}<\theta<\pi-\varepsilon
    \end{cases}
\end{split}
\end{equation}
Next, we use coordinate transformation $\cos\theta = 1-r\lambda, \;\; d\lambda^\zeta=r^{-\zeta}2^{1-\zeta}\zeta\theta^{2\zeta-1}(1+O(\theta))d\theta$ and obtained asymptotic form to calculate the integral in the left-hand side of \eqref{eq:Jacobi_ansatz_powerlaw_integral}
\begin{equation}\label{eq:Jacobi_loss_big_a}
    \begin{split}
        \int_0^1 (q_n^{(a,b,r)}(\lambda))^2 &d\lambda^\zeta \stackrel{(1)}{=} \Gamma^2(a+1)r^{-\zeta}2^{2a+1-\zeta}\zeta\int_0^{\theta_r}\theta^{2\zeta-1-2a}N^{-2a}J_a^2(N\theta)d\theta \\
        &+O(1)\int_0^{cn^{-1}}\theta^{2\zeta}d\theta + O(n^{-2a-1})\int_{cn^{-1}}^{\theta_r}\theta^{2\zeta-1-2a}d\theta\\
        \stackrel{(2)}{=}&\Gamma^2(a+1)r^{-\zeta}2^{2a+1-\zeta}\zeta N^{-2\zeta}\int_0^{N\theta_r}z^{2\zeta-1-2a}J_a^2(z)dz + n^{-2\zeta}O(\delta_{\zeta,a}(n))\\
        =&\Gamma^2(a+1)r^{-\zeta}2^{2a+1-\zeta}\zeta N^{-2\zeta}\int_0^{\infty}z^{2\zeta-1-2a}J_a^2(z)dz + n^{-2\zeta}O(\delta_{\zeta,a}(n))\\
        \stackrel{(3)}{=}&\left(\frac{r}{2}\right)^{-\zeta} \frac{\zeta\Gamma^2(a+1)B(\zeta,2a-2\zeta+1)}{\Gamma^2(a-\zeta+1)}n^{-2\zeta} + n^{-2\zeta}O(\delta_{\zeta,a}(n))
    \end{split}
\end{equation}
where in (1) $\cos\theta_r=1-r$. In (2), error term $\delta_{\zeta,a}(n)$ comes from estimation of the last two integrals in (1) and is given by \begin{equation}\label{eq:Jacobi_ansatz_powerlaw_integral_errorterm}
    \delta_{\zeta,a}(n) = \begin{cases}
    n^{-1}, \hfill a>\zeta&\\
    n^{-1}\log n, \hfill  a=\zeta&\\
    n^{2\zeta-2a-1}, \quad \hfill \zeta-\frac{1}{2}<a<\zeta&
    \end{cases} \quad
\end{equation}
This error term gives more fine-grained characterization of the correction than $o(1)$ term in \eqref{eq:Jacobi_ansatz_powerlaw_integral}, where it was omitted for brevity. Finally, in (3) we used known integral for Bessel function, which can be found e.g. in \href{https://dlmf.nist.gov/}{DLMF} (\S10.22). 

Now we proceed with the second case $a<\zeta-\frac{1}{2}$. Using the first asymptotic in \eqref{eq:Jacobi_square_asym} and analyzing the error terms similarly to \eqref{eq:Jacobi_loss_big_a} we get 
\begin{equation}
    \begin{split}
        \int_0^1 (q_n^{(a,b,r)}(\lambda))^2 &d\lambda^\zeta= \frac{2^\zeta\zeta\Gamma^2(a+1)}{r^\zeta N^{2a+1}}\int_0^{\theta_r}\frac{N\theta J_a^2(N\theta)d(\sin\tfrac{\theta}{2})}{\big(\sin\tfrac{\theta}{2}\big)^{2a-2\zeta+2}\big(\cos\tfrac{\theta}{2}\big)^{2b+1}} +O(n^{-2a-2})\\
        \stackrel{(1)}{=}&\frac{2^\zeta\zeta\Gamma^2(a+1)}{\pi r^\zeta N^{2a+1}}\int_0^{\theta_r}\big(\sin\tfrac{\theta}{2}\big)^{2\zeta-2a-2}\big(\cos\tfrac{\theta}{2}\big)^{-2b-1}d\sin\tfrac{\theta}{2} +o(n^{-2a-1})\\
        \stackrel{(2)}{=}&\frac{2^\zeta\zeta\Gamma^2(a+1)}{2\pi r^\zeta n^{2a+1}}\int_0^{\tfrac{r}{2}}x^{\zeta-a-\tfrac{3}{2}}(1-x)^{-b-\tfrac{1}{2}}dx +o(n^{-2a-1}) \\
        =& \frac{2^\zeta\zeta\Gamma^2(a+1)B(\tfrac{r}{2};\zeta-a-\tfrac{1}{2}, b+\tfrac{1}{2})}{2\pi r^\zeta} n^{-2a-1} \big(1+o(1)\big)
    \end{split}
\end{equation}
Here in (1) we used the property that $\lim_{n\to\infty} \int_0^1 f(x) [nx J_a^2(nx)] dx = \pi^{-1}\int_0^1 f(x) dx$ for functions $f(x)$ integrable on $(0,1)$ and Lipschitz on any $(\varepsilon,1)$. This property follows from $zJ_a(z)^2$ being bounded, and asymptotic of Bessel function $zJ_a^2(z) = 2\pi^{-1}\cos^2(z-\alpha) + O(z^{-1})$. In (2) we changed integration coordinate to $x=\sin^2 \tfrac{\theta}{2}$. 

\paragraph{Learning rate schedule associated with Jacobi ansatz \eqref{eq:Jacobi_ansatz}.} In this section, we obtain the learning rate schedule \eqref{eq:Jacobi_ansatz_parameters}. Note that we can set $r=1$ in derivation but receiver it in the end since it always comes in combination $r\lambda$, therefore multiplicative modifying learning rate. 

Now, we start with standard recurrence relations \eqref{eq:Jacobi_recurrence} and first substitute $x=1-\lambda$:
\begin{equation}
    \begin{split}
        &2(n+1)(n+a+b+1)(2n+a+b) P_{n+1}^{(a,b)}(1-\lambda) = \\ &\qquad-(2n+a+b)(2n+a+b+1)(2n+a+b+2) \lambda P_{n}^{(a,b)}(1-\lambda) \\
        &\qquad + \Big[(2n+a+b+1)(a^2-b^2)+(2n+a+b)(2n+a+b+1)(2n+a+b+2)\Big] P_{n}^{(a,b)}(1-\lambda) \\
        &\qquad - 2(n+a)(n+b)(2n+a+b+2)P_{n-1}^{(a,b)}(1-\lambda).
    \end{split}
\end{equation}
Next step is to add normalization $P_n^{(a,b)}(1-\lambda)=P_n^{(a,b)}(1)p^{(a,b)}_n(\lambda)$, where according to \eqref{eq:Jacobi_normalization} $P_n^{(a,b)}(1) = \frac{\Gamma(n+a+1)}{\Gamma(n+1)\Gamma(a+1)}$. We get
\begin{equation}
    \begin{split}
        &p^{(a,b)}_{n+1}(\lambda) = \\ 
        &\qquad -\frac{(2n+a+b+1)(2n+a+b+2)}{2(n+a+1)(n+a+b+1)}\lambda p^{(a,b)}_{n}(\lambda) \\
        &\qquad+ \Big[\frac{(2n+a+b+1)(2n+a+b+2)}{2(n+a+1)(n+a+b+1)}+\frac{(2n+a+b+1)(a^2-b^2)}{2(n+a+1)(n+a+b+1)(2n+a+b)}\Big] p^{(a,b)}_{n}(\lambda) \\
        &\qquad - \frac{n(n+b)(2n+a+b+2)}{(n+a+1)(n+a+b+1)(2n+a+b)}p^{(a,b)}_{n-1}(\lambda)
    \end{split}
\end{equation}
Comparing with \eqref{eq:HB_poly_recurrence}, this gives exactly \eqref{eq:Jacobi_ansatz_parameters} with $r=1$. Then, $r$ is recovered by setting $\alpha_n \rightarrow r\alpha_n$. Finally, the asymptotic form in \eqref{eq:Jacobi_ansatz_parameters} is obtained by a simple Taylor expansion with respect to $\tfrac{1}{n}$.

\section{Non-constant learning rates: upper bounds}

\subsection{Accelerated Heavy Ball rates}\label{sec:ubcggeneral}
\paragraph{Proof of Theorem \ref{ther:Jacobi_worst_case_UB}}
From the properties of polynomials $q_n^{(a,b,r)}$, we will take only non-degeneracy of zeros and monotonicity of local maxima. The former follows directly from the same property of Jacobi polynomials. The monotonicity property is also inherited from Jacobi polynomials and the respective argument is implicitly given in Section 7.32 of \cite{szego1959orthogonal}. For completeness, we formulate and prove the monotonicity property here.
\begin{lem}\label{lemma:Jacobi_maxima_monotonicity}
Assume $a,b>-\tfrac{1}{2}$ and let $x_0 = \frac{b-a}{a+b+1}$. Next, denote $\{x_i\}_{i=1}^m$ the positions of local maxima of $|P_n^{(a,b)}(x)|$ on $(x_0,1)$ sorted in increasing order: $x_0<x_1<\ldots<x_m<1$. Then, the values at local maxima and at the endpoints form an increasing sequence
\begin{equation}
    |P_n^{(a,b)}(x_0)| < |P_n^{(a,b)}(x_1)| < \ldots < |P_n^{(a,b)}(x_m)| <|P_n^{(a,b)}(1)|
\end{equation}
\end{lem}
\begin{proof}
Recall that $y(x)=P_n^{(a,b)}(x)$ satisfy differential equation
\begin{equation}
    (1-x^2)y'' + \Big(b-a-(a+b+2)x\Big)y' + n(n+a+b+1)y=0
\end{equation}
Then, to characterize $y(x)$ at local extrema we introduce function $f(x)$ and calculate its derivative taking into account differential equation for $y(x)$. 
\begin{align}
    f(x) &=\Big(y(x)\Big)^2 + \frac{1-x^2}{n(n+a+b+1)}\Big(y'(x)\Big)^2\\
    f'(x) &= \frac{2(a+b+1)}{n(n+a+b+1)}(x-x_0)\Big(y'(x)\Big)^2
\end{align}
From the derivative expression we see that $f(x)$ is monotonously increasing on $[x_0,x]$. Now observe that $f(x)=y^2(x)$ at local minima $x_i$ and at endpoint $x=1$, which implies monotonicity of maxima $|P_n^{(a,b)}(x_1)| < \ldots < |P_n^{(a,b)}(x_m)| <|P_n^{(a,b)}(1)|$. For the left endpoint $x_0$ we notice that $y^2(x_0)\leq f(x_0) < f(x_1)=y^2(x_1)$ which completes the proof.   
\end{proof}
Note that according to \eqref{eq:Jacobi_ansatz}, restriction on $r$ means $\lambda\in[0,1]$ maps to $[1-r,1]\subseteq[x_0,1]$ in the argument of $P_n^{(a,b)}(x)$ with $x_0=\tfrac{b-a}{a+b+1}$. Then, according the lemma \ref{lemma:Jacobi_maxima_monotonicity}, for the local maxima $\{\lambda_i\}_{i=1}^m$ of $q_n^2(\lambda)$ on $(0,1)$ we have
\begin{equation}\label{eq:res_polynomial_maxima_monotonicity}
    |q_n(0)| > |q_n(\lambda_1)| >\ldots> |q_n(\lambda_m)|>|q_n(1)|
\end{equation}
From this point, we will not require any additional properties of $q_n^{(a,b,r)}$, and therefore denote $p_n(\lambda)\equiv q_n^{(a,b,r)}(\lambda)$. From the proof of Theorem \ref{ther:worst_case_loss}, we recall the structure of flattened polynomial $\overline{q_n^2}(\lambda)$ given by \eqref{eq:flattened_poly_description}. Monotonicity of local maxima of $q_n(\lambda)$ means that $x_i$ in \eqref{eq:flattened_poly_description} are simply local maxima of $q_n(\lambda)$, and, in particular, $x_1=0$.

Now, we focus on the contribution to the losses \eqref{eq:worst_case_loss} and \eqref{eq:exact_power_measure_loss} from a single flat region $[y_i,x_{i+1}]$. Let $c$ be the root of $p_n(\lambda)$ on $[y_i,x_{i+1}]$, and denote $\widetilde{p}_n(\lambda)=p_n(\lambda)/(\lambda-c)$. As $\widetilde{p}_n(\lambda)$ has all its roots outside of $[y_i,x_{i+1}]$, on this segment $\widetilde{p}_n^2(\lambda)$ is either 1) monotonically increasing and then decreasing 2) monotonically decreasing 3) monotonically increasing. Therefore, the minima of $\widetilde{p}_n^2(\lambda)$ on $[y_i,x_{i+1}]$ is attained at one of the ends of the segments. Taking into account that $\int_{y_i}^{x_{i+1}} \overline{p_n^2}(\lambda)d(\lambda^\zeta)=p_n^2(x_{i+1})(x_{i+1}^\zeta-y_i^\zeta)$, we have 
\begin{equation}
\begin{split}
    &\int_{y_i}^{x_{i+1}} p_n^2(\lambda) d(\lambda^\zeta)\bigg/\int_{y_i}^{x_{i+1}} \overline{p_n^2}(\lambda) d(\lambda^\zeta)\\
    =& \int_{y_i}^{x_{i+1}} \widetilde{p}_n^2(\lambda)(\lambda-c)^2 d(\lambda^\zeta)\bigg/\left(p_n^2(x_{i+1})\int_{y_i}^{x_{i+1}} d(\lambda^\zeta)\right)\\
    \geq&\int_{y_i}^{x_{i+1}} \operatorname{min}\big(\widetilde{p}_n^2(y_i),\widetilde{p}_n^2(x_{i+1})\big)(\lambda-c)^2 d(\lambda^\zeta)\bigg/\left(p_n^2(x_{i+1})\int_{y_i}^{x_{i+1}} d(\lambda^\zeta)\right)\\
    =& \int_{y_i}^{x_{i+1}} \frac{(\lambda-c)^2}{\operatorname{max}\big((y_i-c)^2,(x_{i+1}-c)^2\big)} d(\lambda^\zeta)\bigg/\int_{y_i}^{x_{i+1}} d(\lambda^\zeta) \geq \frac{1}{C[\rho_\zeta]}
\end{split}
\end{equation}
Here, we observed that the expression in the last line is a single realization of the expression minimized in \eqref{eq:parabola_av_infinum}, and therefore can be bounded with respective infimum $\frac{1}{C[\rho_\zeta]}$. Thus, we have bounded the ratio of integrals $\int \overline{p_n^2}(\lambda) d(\lambda^\zeta) \Big / \int p_n^2(\lambda) d(\lambda^\zeta)$ on $[y_i,x_{i+1}]$ with $C[\rho_\zeta]$. As the same bound trivially holds on $[x_i,y_i]$ (flattened and original polynomials are equal), and the respective segments cover the whole $[0,1]$, we get \eqref{eq:Jacobi_ansatz_UB}.

\paragraph{Proof of Proposition \ref{prop:parabola_powerlaw_av_minimum}.}
Let's denote the ratio of integrals under the infimum in \eqref{eq:parabola_av_infinum} as $C[\rho](c, x_l, x_r)$. Then, for the exact power law measure $\rho_\zeta([0,\lambda])=\lambda^\zeta$ the ratio becomes invariant under scaling transformations: $C[\rho_\zeta](\eta c, \eta x_l, \eta x_r)=C[\rho_\zeta](c, x_l, x_r), \forall \eta>0$. This scale invariance implies that it is sufficient only to consider the case $x_r=1$. Now, we can simply denote $x_l=x$. 

We reduce the space of $(x,c)$ potentially containing the infimum by noting that for $c \notin [x,1]$, it is always beneficial to move $c$ to the nearest endpoint of $[x,1]$. Next, we take advantage of monotonicity of the density $p_\zeta(\lambda)=\tfrac{d\rho_\zeta([0,\lambda])}{d\lambda}=\zeta \lambda^{\zeta-1}$ to further narrow down the search space: for any $c \in [x, 1]$ we compare it with its reflection $c'=x+1-c$ w.r.t. window center $c_0=(x+1)/2$
\begin{equation}
\begin{split}
        C[\rho]&(c', x, 1)-C[\rho](c, x, 1)\propto \int_x^1 \big[(\lambda-c')^2-(\lambda-c)^2\big]p(\lambda)d\lambda\\
        &\propto (c-c_0)\int_x^1 \big[\lambda-c_0\big]p(\lambda)d\lambda \propto (c-c_0)\int_0^{\frac{1-x}{2}} z \Big[p\big(c_0-z\big)-p\big(c_0+z\big)\Big]dz
\end{split}
\end{equation}
Here and in the remaining parts of the proof, the proportionality sign $\propto$ denotes equality up to a positive multiplicative factor. From the last line, we see that for increasing density $p(\lambda)$ it is always more beneficial to be in the right half of the window $c>c_0$, and vice versa for decreasing $p(\lambda)$. In the case of constant $p(\lambda)$, as for $\zeta=1$, both halves of the window $[x,1]$ are equivalent. 

The right (left) position of $c$ w.r.t. window center $c_0$ implies that parabola in \eqref{eq:parabola_av_infinum} is normalized by its left(right) endpoint. Slightly abusing the fact that after fixing the normalization endpoint, the positions of $c$ away from the intended half of the window are always suboptimal, we may write
\begin{align}
         &C_\zeta^{-1} = \underset{c,\; 0\leq x < 1}{\operatorname{inf}} C_\zeta(x, c)\\
         &C_\zeta(x, c) \equiv 
         \begin{cases}
            \langle(\lambda-c)^2\rangle_x \big/ \langle(1-c)^2\rangle_x, \quad 0<\zeta\leq1\\
            \langle(\lambda-c)^2\rangle_x \big/ \langle(x-c)^2\rangle_x,  \quad\zeta>1
        \end{cases}
\end{align}
where angle brackets denote the integral $\langle f(\lambda) \rangle_x \equiv \int_x^1 f(\lambda)p_\zeta(\lambda)d\lambda$. Now we proceed with finding the optimal point $(x^*,c^*)=\operatorname{argmin} C_\zeta(x,c)$ and respective value $C_\zeta$ separately for the cases $\zeta\leq1$ and $\zeta>0$. In both cases, it turns out that at the optimum $x^*=0$, which makes it easy to find respective $c^*$. However, showing that $x^*=0$ is technically challenging, and we had to use symbolic computation, e.g. Wolfram Mathematica \cite{Mathematica}. 

\paragraph{Decreasing density ($\zeta\leq1$).} First, let's find optimal $c=c^*(x)$ at a given $x$. Since $C_\zeta(x,c)$ is a rational function in $c$, the optimum is given by a zero of the derivative 
\begin{equation}
    \frac{\partial C_\zeta(x,c)}{\partial c} = 2\frac{\lav{(\lambda-c)^2}}{(1-c)^3 \lav{1}} - 2\frac{\lav{\lambda-c}}{(1-c)^2\lav{1}} \propto \sgn{1-c}\Big(c\lav{1-\lambda}-\lav{\lambda(1-\lambda)}\Big)
\end{equation}
From this expression, we see that the minimum is indeed unique and achieved at
\begin{equation}\label{eq:c_opt_decreasing_density}
c^*(x)=\frac{\lav{\lambda(1-\lambda)}}{\lav{1-\lambda}}   
\end{equation}
Next, as the global minimum of $C_\zeta(x,c)$ is located on the curve $c=c^*(x)$, we may analyze the derivative along the curve $\frac{d}{dx} C_\zeta(x, c^*(x))=\frac{\partial }{\partial x} C_\zeta(x, c^*(x))$
\begin{equation}\label{eq:dC_dl}
\begin{split}
    \frac{\partial }{\partial x} C_\zeta(x, c^*)&\propto \frac{\partial }{\partial x} \frac{\lav{(\lambda-c^*)^2}}{\lav{1}}=-\frac{(x-c^*)^2 p_\zeta(x)}{\lav{1}}+\frac{\lav{(\lambda-c^*)^2}}{(\lav{1})^2}p_\zeta(x)\\
    &\propto \lav{(\lambda-c^*)^2-(x-c^*)^2}=\lav{\lambda^2-x^2}-2c^*\lav{\lambda-x}\\
    &\propto \lav{\lambda^2-x^2}\lav{1-\lambda}-2\lav{\lambda-x}\lav{\lambda(1-\lambda)} \equiv g_1(x)
\end{split}
\end{equation}
Now we will show that $g_1(x)$, and therefore the derivative $\frac{d}{dx} C_\zeta(x, c^*(x))$, is non-negative for $x\in(0,1)$ implying that the global minimum is achieved at $x=0$. First, observe that $g_1(x)$ can be written as an explicit  function of $x$ by substituting moments $\lav{\lambda^n}=\tfrac{\zeta}{\zeta+n}(1-x^{\zeta+n})$. Next, we perform a \emph{top-down} step: use symbolic computations to evaluate several derivatives of $g_1(x)$ in the form of the following statements 
\begin{enumerate}
    \item $g_1(x)=\tfrac{d}{dx}g_1(x)=\big(\tfrac{d}{dx}\big)^2g_1(x)=\big(\tfrac{d}{dx}\big)^3g_1(x)=0$ at $x=1$.
    \item Denote $g_2(x)=x^{3-\zeta}\big(\tfrac{d}{dx}\big)^3g_1(x)$. Then $g_2(x)=\tfrac{d}{dx}g_2(x)=\big(\tfrac{d}{dx}\big)^2g_2(x)=0$ at $x=1$.
    \item $g_2(x=0)=-\zeta(4+\zeta(\zeta^2-4\zeta+7))/(\zeta+1)<0$ for $0<\zeta\leq1$.
    \item $\big(\tfrac{d}{dx}\big)^3g_2(x)=10(1-\zeta)\zeta^2$ at $x=1$, and $\big(\tfrac{d}{dx}\big)^4g_2(x)=8(1-\zeta)\zeta^2(1+2\zeta)x^{\zeta-2}$.
\end{enumerate}
Now we proceed with a \emph{bottom-up} step: use simple expressions of lowest derivatives to reconstruct the positivity of $g_1(x)$. It will be convenient to call \emph{sign signature} of a function the sequence of its signs on a given interval, e.g. $f(x)=(2x-1)^2-0.5$ has sign signature $(+-+)$ on interval $(0,1)$. Then 
\begin{enumerate}
    \item $\big(\tfrac{d}{dx}\big)^4g_2(x)>0$ and $\big(\tfrac{d}{dx}\big)^3g_2(x=1)>0$ implies that $\big(\tfrac{d}{dx}\big)^3g_2(x)$ has sign signature either $(-+)$ or $(+)$ on $(0,1)$.
    \item Sign signature of $\big(\tfrac{d}{dx}\big)^3g_2(x)$ and $\big(\tfrac{d}{dx}\big)^2g_2(x=1)=0$ implies that $\big(\tfrac{d}{dx}\big)^2g_2(x)$ has sign signature either $(+-)$ or $(-)$ on $(0,1)$.
    \item Sign signature of $\big(\tfrac{d}{dx}\big)^2g_2(x)$ and $\tfrac{d}{dx}g_2(x=1)=0$ implies that $\tfrac{d}{dx}g_2(x)$ has sign signature either $(-+)$ or $(+)$ on $(0,1)$.
    \item Sign signature of $\tfrac{d}{dx}g_2(x)$ implies that maximum of $g_2(x)$ on $[0,1]$ is reached either at $x=0$ or $x=1$. Since $g_2(1)=0$ and $g_2(0)<0$, we have $g_2(x) \leq 0$ and therefore $\big(\tfrac{d}{dx}\big)^3g_1(x)\leq0$ on $(0,1)$.
    \item $g_1(x)=\tfrac{d}{dx}g_1(x)=\big(\tfrac{d}{dx}\big)^2g_1(x)=\big(\tfrac{d}{dx}\big)^3g_1(x)=0$ at $x=1$ and $\big(\tfrac{d}{dx}\big)^3g_1(x)\leq0$ on $(0,1)$ implies that $g_1(x)\geq0$ on $(0,1)$, which completes the argument. 
\end{enumerate}

Finally, we can proceed with calculating the value at the global minimum $C_\zeta(x=0,c^*(x=0))$. When $x=0$, the moments are $\langle \lambda^n\rangle_0=\zeta/(\zeta+n)$, which after substitution into \eqref{eq:c_opt_decreasing_density} gives $c^*=\zeta/(\zeta+2)$. Then we again substitute the moments into $C_\zeta(0,c^*)$ and get
\begin{equation}
    C_\zeta^{-1}=C_\zeta(0, c^*(0))= \left(\frac{\zeta+2}{2}\right)^2 \left(\frac{\zeta}{\zeta+2}-2\frac{\zeta}{\zeta+1}\frac{\zeta}{\zeta+2}+\big(\frac{\zeta}{\zeta+2}\big)^2\right)=\frac{\zeta}{2(\zeta+1)}
\end{equation}

\paragraph{Increasing density ($\zeta>1$).}
Similarly to $\zeta\leq1$ case, we start with obtaining optimal $c$ at fixed $x$ by calculating the derivative
\begin{equation}
    \frac{\partial C_\zeta(x,c)}{\partial c} = 2\frac{\lav{(\lambda-c)^2}}{(x-c)^3 \lav{1}} - 2\frac{\lav{\lambda-c}}{(x-c)^2\lav{1}} \propto \sgn{x-c}\Big(\lav{\lambda(\lambda-x)}-c\lav{\lambda-x}\Big)
\end{equation}
which gives the optimal position of the parabola root
\begin{equation}
    c^*(x)=\frac{\lav{\lambda(\lambda-x)}}{\lav{\lambda-x}}.
\end{equation}
Next, we again search for the global minimum of $C_\zeta(x,c)$ on the curve $c=c^*(x)$, by analyzing the derivative along the curve $\frac{d}{dx} C_\zeta(x, c^*(x))=\frac{\partial }{\partial x} C_\zeta(x, c^*(x))$
\begin{equation}\label{eq:dC_dl_2}
\begin{split}
    \frac{\partial }{\partial x} C_\zeta(x, c^*)&= \frac{\partial }{\partial x} \frac{\lav{(\lambda-c^*)^2}}{\lav{(x-c^*)^2}}=p_\zeta(x)\frac{\lav{(\lambda-c^*)^2-(x-c^*)^2}}{(x-c^*)^2\lav{1}^2}-2\frac{\lav{(\lambda-c^*)^2}}{(x-c^*)^3\lav{1}}\\
    &\propto -p_\zeta(x)(x-c^*)\lav{\lambda^2-x^2-2(\lambda-x)c^*}+2\lav{(\lambda-c^*)^2}\lav{1}\\
    &=\frac{\lav{(\lambda-x)^2}}{\lav{\lambda-x}^2}\Big[2\lav{1}(\lav{\lambda^2}\lav{1}-\lav{\lambda}^2)-p_\zeta(x)\lav{\lambda-x}\lav{(\lambda-x)^2}\Big]\\
    &\propto 2\lav{1}(\lav{\lambda^2}\lav{1}-\lav{\lambda}^2)-p_\zeta(x)\lav{\lambda-x}\lav{(\lambda-x)^2} \equiv g_1(x)
\end{split}
\end{equation}
Continuing the same strategy as for the case $\zeta\leq1$, we will show $g_1(x)>0$ on $(0,1)$ by exploiting the explicit form of $g_1(x)$ and symbolic computations. \emph{top-down} step:
\begin{enumerate}
    \item $\frac{\tfrac{d}{dx}g_1(x)}{(1-x)x^{\zeta-2}}\equiv g_2(x)$ is a polynomial in variables $(x, x^\zeta)$.
    \item $g_1(x)=g_2(x)=\tfrac{d}{dx}g_2(x)=\big(\tfrac{d}{dx}\big)^2 g_2(x)=0$ at $x=1$.
    \item $\big(\tfrac{d}{dx}\big)^3 g_2(x)=2(\zeta-1)\zeta^3 (1-x)x^{\zeta-2}>0$ on $(0,1)$.
\end{enumerate}
Then, the \emph{bottom-up} argumentation is the following
\begin{enumerate}
    \item $\big(\tfrac{d}{dx}\big)^2g_2(x=1)=0$ and $\big(\tfrac{d}{dx}\big)^3g_2(x)>0$ on $(0,1)$ implies $\big(\tfrac{d}{dx}\big)^2g_2(x)<0$ on $(0,1)$.
    \item $\tfrac{d}{dx}g_2(x=1)=0$ and $\big(\tfrac{d}{dx}\big)^2g_2(x)<0$ on $(0,1)$ implies $\tfrac{d}{dx}g_2(x)>0$ on $(0,1)$.
    \item $g_2(x=1)=0$ and $\tfrac{d}{dx}g_2(x)>0$ on $(0,1)$ implies $g_2(x)<0$ on $(0,1)$, and, therefore, $\tfrac{d}{dx}g_1(x)<0$ on $(0,1)$.
    \item $g_1(x=1)=0$ and $\tfrac{d}{dx}g_1(x)<0$ on $(0,1)$ implies $g_1(x)>0$ on $(0,1)$.
\end{enumerate}
Having shown that at the minimum $x^*=0$, we find the optimal position of the parabola root to be $c^*=c^*(x=0)=\tfrac{\zeta+1}{\zeta+2}$. Plugging $c^*$ into $C_\zeta(0,c)$ gives
\begin{equation}
    C_\zeta^{-1}=C_\zeta(0, c^*)= \left(\frac{\zeta+2}{\zeta+1}\right)^2\left(\frac{\zeta}{\zeta+2}-2\frac{\zeta(\zeta+1)}{(\zeta+1)(\zeta+2)} + \left(\frac{\zeta+1}{\zeta+2}\right)^2\right)=\frac{1}{(\zeta+1)^2}
\end{equation}

\subsection{Gradient Descent with predefined schedule}\label{sec:ubgdpredefined}

\paragraph{Preliminaries: ``reduced'' polynomials.} We will use a construction based on ``reduced'' polynomials $p_{n,m}(x), 0\leq m\leq n$. Given a residual (equal to 1 at $x=0$) polynomial $p_n(x)$ of degree $n$ we define the corresponding reduced polynomials by 
\begin{equation}\label{eq:reduced_polynomial_bound}
    p_{n,m}(x) \equiv \prod_{i=1}^m \Big(1-\frac{x}{x_i}\Big),
\end{equation}
where $x_i$ are the roots of $p_n(x)$ sorted in the decreasing order $x_1\ge x_2\ge \ldots$. In particular, $p_{n,n}(x)=p_n(x)$. We will need the following technical lemma about residual polynomials 
\begin{lem}\label{ther:reduced_polynomials}
Let $p_n(x)$ be a residual polynomial of degree $n$ such that $|p_n(x)|\leq 1$ if $x\in [0,a]$. Then the same bound also holds for the corresponding reduced polynomials:
\begin{equation}
    |p_{n,m}(x)| \leq 1 \; \text{ for all } \; x\in[0,a] \; \text{ and } \; 0 \leq m \leq n. 
\end{equation}
\end{lem}
\begin{proof}
Let's fix $m$ and divide the segment $[0,a]$ into two parts: $[0,2x_m]$ and $[2x_m, a]$ We will prove bound \eqref{eq:reduced_polynomial_bound} separately for each part. (If $x_m < 0$ or $2x_m > a$, then there is only one nontrivial part that covers $[0,a]$, and we consider only the respective single case.) Recall that the initial polynomial $p_n(x)$ and reduced polynomial $p_{n,m}(x)$ can be written as
\begin{equation}
    p_n(x) = \prod_{i=1}^n\Big(1-\frac{x}{x_i}\Big), \quad p_{n,m}(x) = \prod_{i=1}^m\Big(1-\frac{x}{x_i}\Big).
\end{equation}
\begin{enumerate}
    \item \emph{Case $x\in [0,2x_m]$.} In this case we have $|1-\frac{x}{x_i}|\le 1, \; i\leq m,$ and thus $|p_{n,m}(x)|\leq1$.
    \item \emph{Case $x\in [2x_m,a].$} In this case we write
\begin{equation}\label{eq:reduced_polynomial_decomp}
    |p_{n,m}(x)| = \frac{|p_n(x)|}{\prod_{i=m+1}^n |1-\frac{x}{x_i}|}.
\end{equation}
Then for $x\in[2x_m, a]$ and $i>m$, if $x_i>0$, then
\begin{equation}
\Big|1-\frac{x}{x_i}\Big| = \frac{x}{x_i} - 1 \ge \frac{x}{x_m} - 1 \ge 1.
\end{equation}
The same inequality  $|1-\tfrac{x}{x_i}|\ge 1$ clearly also holds if $x_i<0$. Thus, in any case $|1-\tfrac{x}{x_i}|\ge 1$.
It follows then from \eqref{eq:reduced_polynomial_decomp}  that $p_{n,m}(x) \leq p_n(x) \leq 1$. 
\end{enumerate}
\end{proof}

\paragraph{Construction of learning rates $\alpha_n$.} Given $\zeta>0$, fix some $a>b>-\tfrac{1}{2}$ and $r\le 2$ and consider the residual polynomials $p_n$ obtained by shifting and normalizing the Jacobi polynomials as in Eq. \eqref{eq:Jacobi_ansatz}:     
\begin{equation}
    p_{n}(x) =  \frac{P^{(a,b)}_n(1-r\lambda)}{P^{(a,b)}_n(1)}.
\end{equation}

A well-known result from \cite{szego1959orthogonal} states that if $a>b\geq -\frac{1}{2},$ then the largest value of the Jacobi polynomial $P_n^{(a,b)}$ on the segment $[-1,1]$ is reached at $z=1$:
\begin{equation}
    \underset{|z|\leq 1}{\max} \left|P_n^{(a,b)}(z)\right| = P_n^{(a,b)}(1).
\end{equation}
It follows that our polynomials $p_n$  satisfy the condition $\max_{0\le x\le 1}|p_n(x)|=1$ of Lemma \ref{ther:reduced_polynomials}. 

Now we describe a construction of schedule $\{\alpha_i\}$ which gives the convergence rate $O(n^{-2\zeta})$ for GD.  Informally, we will build our GD polynomial $\widetilde{q}_k(x)$ by sequentially taking the roots of $p_1(x), p_2(x), p_4(x), p_8(x), \ldots$. More precisely, to determine $\alpha_i$ we first find the largest $l$ such that $i\ge 2^{l}$, and denote $l_i\equiv l$, $n_i\equiv2^l$, $m_i\equiv i - 2^l + 1$. Then we set 
\begin{equation}\label{eq:GD_schedule}
    \alpha_i = \frac{1}{x^{(n_i)}_{m_i}},
\end{equation} 
where $x^{(n)}_m$ is the $m$'th root of $p_n(x)$ (as usual, taken in decreasing order). In this way the polynomial $\widetilde{q}_k(x)$ corresponding to our scheduled GD is
\begin{equation}\label{eq:wqkx}
    \widetilde{q}_k(x) = p_{n_k,m_k}(x)\prod_{l=0}^{l_k-1} p_{2^l}(x).
\end{equation}
We can now prove the main result.

\paragraph{Proof of Theorem \ref{th:ubgdnonconst}.} As already mentioned, the polynomials $p_n$ satisfy the hypothesis of Lemma \ref{ther:reduced_polynomials} and so, by this lemma, 
$|p_{n,m}(x)|\leq 1$ on $[0,1]$. We apply this bound to $\widetilde{q}_k(x)$:
\begin{equation}
    |\widetilde{q}_k(x)| = |p_{n_k,m_k}(x)| |p_{n_k/2}(x)| \prod_{l=0}^{l_k-2} |p_{2^l}(x)| \leq |p_{n_k/2}(x)|
\end{equation}
Using Corollary \ref{corol:HB_upper_bound}, we then get
\begin{equation}
    L(\mathbf w_k) \leq C_\zeta R(a,r,\zeta)  \left(\frac{n_k}{2}\right)^{-2\zeta} \Big(1+O\big(\frac{1}{n_k}\big)\Big) \leq C_\zeta R(a,r,\zeta) 4^{2\zeta} k^{-2\zeta} \Big(1+O\big(\frac{1}{k}\big)\Big),
\end{equation}
where in the last inequality we used $\frac{k}{2} < n_k \leq k$.
This completes the proof of Theorem \ref{th:ubgdnonconst}.

\subsection{Conjugate Gradients: discrete spectrum}\label{sec:ubcgdiscrete}

\paragraph{Proof of Theorem \ref{th:cgdiscr}. } Consider the degree-$n$ residual polynomial $q_n$ of the form
\begin{equation}
    q_n(\lambda)=\prod_{s=1}^n(1-\lambda/a_s)=\Big(\prod_{s=1}^{\lfloor n/2\rfloor}(1-\lambda/\lambda_s)\Big) r_n(x)
\end{equation}
where $\lambda_1\ge\ldots\ge\lambda_{\lfloor n/2\rfloor}$ are the $\lfloor n/2\rfloor$ largest eigenvalues (atoms of the measure $\rho$) and $r_n$ is some degree-$(n-\lfloor n/2\rfloor)$ residual polynomial. Then, 
\begin{align}
    \int_0^{1} q_n^2(\lambda)\rho(d\lambda)
    ={}&\int_{0}^{\lambda_{\lfloor n/2\rfloor}} q_n^2(\lambda)\rho(d\lambda)\\
    \le{}&\int_{0}^{\lambda_{\lfloor n/2\rfloor}} r_n^2(\lambda)\rho(d\lambda)\\
    \le{}& \int _{0}^{\Lambda (n/2)^{-\nu}} r_n^2(\lambda)\rho(d\lambda)\\
    ={}& \int _{0}^{1} r_n^2(\Lambda (n/2)^{-\nu} t)\rho(d \Lambda (n/2)^{-\nu}t)\\
    ={}& \Lambda^\zeta (n/2)^{-\nu\zeta}\int _{0}^{1} r_n^2(\Lambda (n/2)^{-\nu} t)\rho_n(dt),\label{eq:intrhok}
\end{align}
where the measure $\rho_n$ is defined for Borel subsets $X\subset \mathbb R$ by rescaling
\begin{equation}
    \rho_n(X) = \Lambda^{-\zeta} (n/2)^{\nu\zeta}\rho(\Lambda (n/2)^{-\nu}X).
\end{equation}
The measure $\rho_n$ satisfies the same  power law bound \eqref{eq:rhozeta1} as $\rho$:
\begin{align}
    \rho_n((0,\lambda])={}& \Lambda^{-\zeta} (n/2)^{\nu\zeta}\rho((0,\Lambda(n/2)^{-\nu}\lambda])\\
    \le{}& \Lambda^{-\zeta} (n/2)^{\nu\zeta} (\Lambda(n/2)^{-\nu}\lambda)^\zeta\\
    ={}& \lambda^\zeta.
\end{align}
It follows that we can apply Corollary \ref{corol:HB_upper_bound} and find $r_n$ such that 
\begin{align}
    \int _{0}^{1} r_n^2(\Lambda (n/2)^{-\nu} t)\rho_n(dt)\le 2QC_\zeta R(a,r,\zeta) (n/2)^{-2\zeta} \big(1+o(1)\big).
\end{align}
Combining with \eqref{eq:intrhok}, this gives the desired bound \eqref{eq:cgdisclwn}:
\begin{equation}
    L(\mathbf w_n)\le \frac{1}{2}\int_0^{1} q_n^2(\lambda)\rho(d\lambda)\le QC_\zeta R(a,r,\zeta) \Lambda^{-\zeta} (n/2)^{-(\nu+2)\zeta} \big(1+o(1))\big).
\end{equation}

\section{Non-constant learning rates: lower bounds}

\subsection{Non-adaptive schedules}\label{sec:lbpredefined}
\paragraph{Proof of Theorem \ref{ther:discrete_spectra_predefined_schedule_LB}. }
Consider the power law distribution $\rho_\zeta((0,\lambda])=\lambda^\zeta$ with $\lambda_{\max}=1$. Let us define discrete distributions $(\rho_{\zeta,r})_{r\in[0,1]}$ subject to the spectral conditions \eqref{eq:lndisc1}, \eqref{eq:lndisc2} of the theorem and such that
\begin{equation}\label{eq:rzetadr}
    \rho_\zeta = \int_{0}^1 \rho_{\zeta, r}dr.
\end{equation}
To this end, we set
\begin{equation}
    \rho_{\zeta,r}=\sum_{k=1}^\infty \rho([(k+1)^{-\nu}, k^{-\nu}])\delta_{a_{k,r}}
\end{equation}
with some $a_{k,r}\in [(k+1)^{-\nu},k^{-\nu}].$ It is clear that thus defined $\rho_{\zeta,r}$ satisfies
Eqs. \eqref{eq:lndisc1},
\eqref{eq:lndisc2}, and one can also satisfy Eq. \eqref{eq:rzetadr} by suitably adjusting $a_{k,r}$. 

We will construct the distribution $\rho$ corresponding to the desired $A$ and $\mathbf b$ by joining a sequence of segments of the distributions $\rho_{\zeta,r}:$
\begin{equation}\label{eq:rhohk}
    \rho = \sum_{k=1}^\infty \rho_{\zeta, r_k}|_{[h_k,h_{k-1})},\quad h_0=1.
\end{equation}
It is easy to see that if $h_k\to 0+$ sufficiently fast, say $h_k\le h_{k-1}/2$ for all $k$, then such $\rho$ also satisfies the required conditions \eqref{eq:lndisc1}, \eqref{eq:lndisc2}.

Consider the first step of the construction of $\rho$. Arguing as in Section \ref{sec:poly},  the loss $L(\mathbf w_n)$ of a general multistep method \eqref{eq:multistep} can be written as 
\begin{equation}
    L(\mathbf w_n)=\frac{1}{2}\int_0^1  q_n^2(\lambda)\rho(d\lambda),
\end{equation}
where $q_n$ is some residual polynomial of degree $n$.
We know from the exact solution of the minimization problem \begin{equation}
    \min_{\widetilde q_{n_1}:\deg \widetilde q_{n_1}=n_1, \widetilde q_{n_1}(0)=1}\frac{1}{2}\int_0^1 \widetilde q^2_{n_1}\rho_\zeta(d\lambda)
\end{equation} by a rescaled Jacobi polynomial (see Theorem \ref{th:cgexactpower}) that
\begin{equation}\label{eq:lbjac}
    \frac{1}{2}\int_0^1 q^2_{n_1}\rho_\zeta(d\lambda)> Cn_1^{-2\zeta},
\end{equation}
where $q_{n_1}$ is the residual polynomial corresponding to the given optimization algorithm and $C$ is an absolute constant. 
Choose $n_1$ sufficiently large so that
\begin{equation}
    \frac{1}{2}\int_0^1 q^2_{n_1}\rho_\zeta(d\lambda)> n_1^{-2\zeta-\epsilon}.
\end{equation}
It follows from the decomposition \eqref{eq:rzetadr} that there exists $r_1$ such that this inequality remains valid if we replace $\rho_\zeta$ by $\rho_{\zeta,r_1}:$
\begin{equation}
    \frac{1}{2}\int_0^1 q^2_{n_1}\rho_{\zeta,r_1}(d\lambda)> n_1^{-2\zeta-\epsilon}.
\end{equation}
We can then choose $h_1$ sufficiently small so that 
\begin{equation}
    \frac{1}{2}\int_{h_1}^1 q^2_{n_1}\rho_{\zeta,r_1}(d\lambda)> n_1^{-2\zeta-\epsilon}.
\end{equation}

Consider now the second step of the construction of $\rho$. Using the homogeneity of the distribution $\rho_\zeta,$ the lower bound \eqref{eq:lbjac} extends to the segment $[0,h_1]$ with the additional factor $h_1^\zeta:$
\begin{equation}
    \frac{1}{2}\int_0^{h_1} q^2_{n_2}\rho_{\zeta}(d\lambda)> Ch_1^\zeta n_2^{-2\zeta}.
\end{equation}
Arguing as before, we then choose a sufficiently large $n_2$, a suitable $r_2$, and a sufficiently small $h_2$ such that
\begin{equation}
    \frac{1}{2}\int_{h_2}^{h_1} q^2_{n_2}\rho_{\zeta,r_2}(d\lambda)> n_2^{-2\zeta-\epsilon}.
\end{equation}
Continuing this process, we obtain the full desired expansion \eqref{eq:rhohk}.

\subsection{CG with discrete spectrum}\label{sec:lbcgdiscrete}

\subsubsection{Proof of Proposition \ref{lm:resolv2} for $0<\zeta<1$}\label{sec:lmresolv}
In this section we prove Proposition \ref{lm:resolv2} for $0<\zeta<0$, i.e. we prove only Statement 1 with $m=0$. The remaining cases will be considered in Section \ref{sec:lmresolv2}.

Denote $\mathbf x=(\widetilde A+\epsilon)^{-1}\mathbf e_1$. In coordinates, the equation $(\widetilde A+\epsilon)\mathbf x=\mathbf e_1$ is a system of finite difference equations  
\begin{align}
    -(n-1)^{-\nu}(\tfrac{n}{n-1})^{g} x_{n-1}&\\
    +((n-1)^{-\nu}(\tfrac{n}{n-1})^{2g}+n^{-\nu})x_n&\\
    -n^{-\nu}(\tfrac{n+1}{n})^g x_{n+1} &{}= -\epsilon x_n,\quad n = 2,3,\ldots\\
    x_1-2^g x_{2} &{}= 1-\epsilon x_1,
\end{align}
where we introduced the constant
\begin{equation}\label{eq:gdef}
    g=\frac{1-(2+\nu)\zeta}{2}.
\end{equation}
Let us make the substitution 
\begin{equation}\label{eq:ynxn}
    y_n=n^gx_n.
\end{equation}
Then the finite difference equations become
\begin{align}
    -(n-1)^{-(\nu+2g)}n^gy_{n-1}&\\
    +((n-1)^{-(\nu+2g)}n^{2g}+n^{-\nu})n^{-g}y_n&\\
    -n^{-(\nu+g)}y_{n+1} &{}= -\epsilon n^{-g}y_n,\quad n = 2,3,\ldots\\
    y_1-y_{2} &{}= 1-\epsilon y_1.\label{eq:bcy1g}
\end{align}
We further introduce the variable $h$ by
\begin{equation}\label{eq:heps}
    h = \epsilon^{1/(2+\nu)}
\end{equation}
and the variable $\theta_n$ by
\begin{equation}\label{eq:deftheta}
    1-h\theta_n = \frac{y_n}{y_{n+1}}.
\end{equation}
By multiplying the difference equation by $n^{\nu+g}/y_n,$ we can then rewrite it as
\begin{align}
    -(\tfrac{n-1}{n})^{-(\nu+2g)}(1-h\theta_{n-1})&\\
    +(\tfrac{n-1}{n})^{-(\nu+2g)}+1&\\
    -(1-h\theta_n)^{-1} &{}= -h^2 (hn)^{\nu},\quad n = 2,3,\ldots.
\end{align}
Introducing the variable $s$ by 
\begin{equation}
    s = nh,
\end{equation}
we then get
\begin{equation}\label{eq:thetaniter2}
    \theta_{n-1}=\Big(\frac{\theta_n}{1-h\theta_n}-hs^\nu\Big)\Big(\frac{s-h}{s}\Big)^{\nu+2g},\quad n=2,3,\ldots
\end{equation}
This system of finite difference equations has a one-parameter family of solutions that can be specified by one value $\theta_{n_0}$ at a particular $n_0$. We will now identify a special solution $(\theta^*_n)$ for which $\mathbf x\in l^2$. We expect the components $x_n$ of this $\mathbf x$ to have the same sign and decay to 0 sufficiently fast as $n\to+\infty$. By Eq. \eqref{eq:deftheta}, these conditions will be satisfied if we ensure that $\theta^*_n<0$ for all $n$ and $\theta^*_n\to -\infty$ sufficiently fast as $n\to \infty$ (note that this need not be the case for a generic solution $(\theta_n)$ since it may diverge at a finite $n$ or start increasing at some $n$). Importantly, we will establish growth bounds for the solution $(\theta_n^*)$ that hold uniformly in $h$.  

\begin{lem}\label{lm:thetasab}
Let constants $a,b$ be such that $a>\nu$ and $0<b<\nu/2$. Then there exists a unique solution $(\theta^*_n)$ of Eq. \eqref{eq:thetaniter2} such that we have
\begin{equation}\label{eq:abtheta}
    -s^a\le\theta^*_n\le-s^b\text{ provided }h<h_0\text{ and }s=nh>s_0
\end{equation}
with some constants $h_0,s_0>0$.
\end{lem}
\begin{proof}
Let $G_{s,h}:\mathbb R\to\mathbb R$ denote the transformation in the iteration law \eqref{eq:thetaniter2}:
\begin{equation}
    \theta_{n-1}=G_{nh,h}(\theta_n).
\end{equation}
Consider the intervals
\begin{equation}
    I_s = [-s^a, -s^b].
\end{equation}
We show now that under our iteration law the intervals $I_s$ are ordered by inclusion.
\begin{lem}\label{lm:Gnh} There exist constants $h_0,s_0>0$ such that for all $h<h_0$ and $s=hn>s_0$ we have
\begin{equation}\label{eq:ghihn}
    G_{nh, h}(I_{hn})\subset I_{h(n-1)}.
\end{equation}
\end{lem}
\begin{proof} By monotonicity of $G_{h,s}$, Eq. \eqref{eq:ghihn} will be established if we show
\begin{align}
    -(s-h)^b \ge{}& G_{s,h}(-s^b),\label{eq:shb}\\
    -(s-h)^a \le{}& G_{s,h}(-s^a).\label{eq:sha}
\end{align}
\emph{Fulfilling condition \eqref{eq:shb}.} This inequality is equivalent to
\begin{equation}
    -(1-h/s)^{b-(\nu+2g)}\ge -\frac{1}{1+hs^b}-hs^{\nu-b}.
\end{equation}
Since we assume that $h$ is sufficiently small and $s$ sufficiently large, we can write $-(1-h/s)^{b-(\nu+2g)}\ge -1-Ch/s$ with some absolute constant C. Therefore, it is sufficient to establish 
\begin{equation}
    -Ch/s\ge \frac{hs^b}{1+hs^b}-hs^{\nu-b}.
\end{equation}
Dividing by $h$ and bounding $1+hs^b\ge 1$, this in turn reduces to
\begin{equation}
    -C/s\ge s^b-s^{\nu-b}.
\end{equation}
Clearly, this inequality holds for sufficiently large $s$ if  $b<\nu/2$.

\emph{Fulfilling condition \eqref{eq:sha}.} By a similar argument, it suffices to fulfill 
\begin{equation}
    C/s\le \frac{s^a}{1+h_0s^a}-s^{\nu-a}.
\end{equation}
This holds for all sufficiently large $s$ if we choose any $a>\nu$ and $h_0$ small enough.
\end{proof}
Lemma \ref{lm:Gnh} yields a nested sequence of compact intervals
\begin{equation}
    I_{hn_0} \supset G_{h(n_0+1), h}(I_{h(n_0+1)})\supset G_{h(n_0+1), h}(G_{h(n_0+2), h}(I_{h(n_0+2)}))\supset\ldots, 
\end{equation}
where $n_0=\lceil s_0/h\rceil$. This sequence has a non-empty intersection $I$. Then, a sequence $\theta_n^*$ such that $\theta_{n_0}^*\in I$ satisfies the desired bounds \eqref{eq:abtheta}. 

We argue now that such a sequence $\theta_n^*$ is unique. It is easy to see  that if a solution $\theta_n^*$ satisfies the upper bound in \eqref{eq:abtheta}, then the respective sequence $x_n$ belongs to $l^2$. Different sequences $\theta_n^*$ would correspond to different $l^2$ sequences $x_n$. However, the equation $(\widetilde A+\epsilon)\mathbf x=\mathbf e_1$ has a unique $l^2$ solution $\mathbf x$.
\end{proof}
We study now the behavior of $\theta_n^*$ at small $n$.
It is convenient to introduce the new variables $\omega_n$ by
\begin{equation}\label{eq:thetaomega}
    \theta_n= s^{\nu+2g}\omega_n.
\end{equation}
Then the difference equation \eqref{eq:thetaniter2} becomes
\begin{equation}\label{eq:omeganit2}
    \omega_{n-1}=\frac{\omega_n}{1-hs^{\nu+2g}\omega_n}-hs^{-2g},\quad n=2,3,\ldots
\end{equation}
Let $\omega_n^*$ be the sequence $\omega_n$ corresponding to the sequence $\theta_n^*$ found in Lemma \ref{lm:thetasab}, and $s_0$ be as in this lemma.
\begin{lem}\label{lm:omegacd}
Let $0<\zeta<1$. Then there exist constants $c < d < 0$ such that 
\begin{equation}
    c\le \omega_n^*\le d \text{ provided } h<h_0\text{ and } s=nh<s_0 
    \end{equation}
with some constant $h_0>0$.
\end{lem}
\begin{proof}
\emph{Lower bound.} By Lemma \ref{lm:thetasab} we have $\omega_{n}^*<0$ for $n\ge n_0=\lceil s_0/h \rceil$, and Eq. \eqref{eq:omeganit2} then implies that $\omega_n^*<0$ for all $n$; moreover,
\begin{equation}
    \omega_{n-1}^*\ge\omega_n^*-hs^{-2g},\quad s=nh,\quad n=2,3,\ldots,n_0.
\end{equation}
Note that by the definition of $g$ in Eq. \eqref{eq:gdef} and the inequality $\zeta>0$ we have 
\begin{equation}
    2g < 1.
\end{equation}
It follows that for any $n=1,2,\ldots,n_0$
\begin{align}
    \omega_n^*\ge{}&\omega_{n_0}^*-\int_{nh}^{n_0h}s^{-2g}ds+O(h+h^{1-2g})\\
    \ge{}& -s_0^{-(\nu+2g)}s_0^a-(1-2g)^{-1}s_0^{1-2g}+O(1)\ge c,\quad (h\to 0)
\end{align}
for a suitable constant $c.$

\emph{Upper bound.} On the other hand, Eq. \eqref{eq:omeganit2} implies that
\begin{equation}
    \omega_{n-1}^*\le \frac{\omega_n^*}{1-hs^{\nu+2g}\omega_n^*},\quad n=2,3,\ldots
\end{equation}
This is equivalent to
\begin{equation}
    (\omega_{n-1}^*)^{-1}\ge (\omega_{n}^*)^{-1} -hs^{\nu+2g},\quad n=2,3,\ldots
\end{equation}
Note that by the definition \eqref{eq:gdef} of $g$ we have
\begin{equation}\label{eq:nu2g1}
    \nu+2g+1=(2+\nu)(1-\zeta),
\end{equation}
so that the inequality $\zeta<1$ yields 
\begin{equation}\label{eq:nu2ggrt1}
    \nu+2g > -1.
\end{equation}
It follows that for any $n=1,2,\ldots,n_0$
\begin{align}
    (\omega_n^*)^{-1}\ge{}&(\omega_{n_0}^*)^{-1}-\int_{nh}^{n_0h}s^{\nu+2g}ds+O(h+h^{\nu+2g+1})\\
    \ge{}& -s_0^{\nu+2g}s_0^{-b}-(\nu+2g+1)^{-1}s_0^{\nu+2g+1}+O(1)\ge d^{-1},\quad (h\to 0)
\end{align}
for a suitable constant $d<0,$ which implies the desired bound.\hfill$\Box$
\end{proof}

Lemmas \ref{lm:thetasab} and \ref{lm:omegacd} allow us to control the initial element $x_1$ of the sequence $\mathbf x$. From Eqs. \eqref{eq:ynxn}, \eqref{eq:bcy1g}, and \eqref{eq:deftheta} we have
\begin{equation}\label{eq:x1theta}
    x_1=y_1=\Big(\epsilon-\frac{h\theta_1^*}{1-h\theta_1^*}\Big)^{-1}.
\end{equation}
By Eqs. \eqref{eq:heps}, \eqref{eq:thetaomega}, \eqref{eq:nu2g1}, and Lemma \ref{lm:omegacd} we have 
\begin{equation}
    h\theta_1^*=h^{\nu+2g+1}\omega_1^*=\epsilon^{1-\zeta}\omega^*_1=O(\epsilon^{1-\zeta}),\quad (\epsilon\to 0).
\end{equation}
Since $\omega_1^*<d<0$ for all $\epsilon$, it follows that if $\zeta<1$, then
\begin{equation}\label{eq:x1oeps}
    x_1=O(\epsilon^{\zeta-1}),\quad (\epsilon\to 0).
\end{equation}
This is the desired bound, since $x_1=\langle \mathbf e_1, (\widetilde{\mathcal A}+\epsilon)^{-1}\mathbf e_1\rangle.$

\subsubsection{Proof of Proposition \ref{lm:resolv2} for $\zeta>1$}\label{sec:lmresolv2}
We retain the notation introduced in the previous section. Throughout this section, we write $a_n=O(b_n)$  meaning that $|a_n|\le Cb_n$ for all $n$ with some constant $C>0$ that might depend on $\nu$ and $\zeta$ but not $n$ or $\epsilon.$ 

We start with a technical lemma that describes the special  solution $\theta_n^*$ for $\zeta>1$ (thus complementing Lemma \ref{lm:omegacd} that covers $\zeta<1$). 
\begin{lem}\label{lm:zeta1theta}
If $\zeta>1$, then for sufficiently small $\epsilon$ the special sequence $\theta^*_n$ satisfies
\begin{equation}\label{eq:thetabound2}
    \theta^*_{n}\le \frac{(\nu+2)(1-\zeta)}{nh}.
\end{equation}
\end{lem}
\begin{proof} By Lemma \ref{lm:thetasab}, if $h$ is small enough then for sufficiently large $n$ we have $\theta_n^*\le -s^{\nu/2}$ and hence bound \eqref{eq:thetabound2} is satisfied if $n$ is large enough. We prove now that if it is satisfied for some $n\ge 2$, then it is also satisfied for 
$n-1$. Consider Eq. \eqref{eq:thetaniter2} for  $\theta^*_{n-1}$: 
\begin{equation}\label{eq:thetaniter3}
    \theta^*_{n-1}=\Big(\frac{\theta^*_n}{1-h\theta^*_n}-hs^\nu\Big)\Big(\frac{s-h}{s}\Big)^{\nu+2g}.
\end{equation}
Recall that $\nu+2g+1=(\nu+2)(1-\zeta).$  Denote $a=(\nu+2)(1-\zeta)$. Using the fact that the function $x\mapsto x/(1-hx)$ is increasing on $(-\infty,0)$ and the assumption $\theta_n^*\le a/s$, we get
\begin{align}
    \theta^*_{n-1}\le{}&\frac{\theta^*_n}{1-h\theta^*_n}\Big(\frac{s-h}{s}\Big)^{\nu+2g}\\
    \le{}&\frac{a/s}{1-(h/s)a}\Big(\frac{s-h}{s}\Big)^{\nu+2g}\\
    ={}&\frac{1}{s-h}\frac{a}{1-(h/s)a}\Big(\frac{s-h}{s}\Big)^{a}\\
    \le{}&\frac{a}{s-h},
\end{align}
where in the last step we used the inequality
\begin{equation}
    \frac{(1-x)^a}{1-xa}\ge 1,\quad 0<x<1, a<0
\end{equation}
with $x=h/s$.  
\end{proof}
Our proof of Proposition \ref{lm:resolv2} is based on the following extended version of this proposition that contains bounds on the growth of the involved sequences. 

\begin{prop}\label{lm:resolvseq}Let $n_0=\lfloor s_0/ h\rfloor$ with the constant $s_0$ appearing in Lemma \ref{lm:thetasab}.
\begin{enumerate}
    \item Assuming $2m<\zeta$ for some integer $m\ge 1$, the vectors $\widetilde A^{-m}\mathbf e_1$ and $\widetilde A^{-m}(\widetilde A+\epsilon)^{-1}\mathbf e_1$ exist as elements of $l^2$ and
    \begin{align}
        (\widetilde A^{-m}\mathbf e_1)_n = {}& O(n^{\frac{-1-(2+\nu)(\zeta-2m)}{2}}), \label{eq:ame1}\\
        (\widetilde A^{-m}(\widetilde A+\epsilon)^{-1}\mathbf e_1)_n={}&\begin{cases}
        O(n^{\frac{-1-(2+\nu)(\zeta-2m-2)}{2}}),& n\le n_0\\
        O(\epsilon^{-1}n^{\frac{-1-(2+\nu)(\zeta-2m)}{2}}),&n> n_0
        \end{cases}\label{eq:ame2}
    \end{align}
    \item Assuming $2m+1<\zeta$ for some integer $m\ge 0$, the vectors $J^{-1}\widetilde A^{-m}\mathbf e_1$ and $J^{-1}\widetilde A^{-m}(\widetilde A+\epsilon)^{-1}\mathbf e_1$ exist as elements of $l^2$ and
    \begin{align}
        (J^{-1}\widetilde A^{-m}\mathbf e_1)_n = {}& O(n^{\frac{-1-(2+\nu)(\zeta-2m-1)}{2}}), \label{eq:j1wa1}\\
        (J^{-1}\widetilde A^{-m}(\widetilde A+\epsilon)^{-1}\mathbf e_1)_n={}&\begin{cases}
        O(n^{\frac{-1-(2+\nu)(\zeta-2m-3)}{2}}),& n\le n_0\\
        O(\epsilon^{-1}n^{\frac{-1-(2+\nu)(\zeta-2m-1)}{2}}),&n> n_0
        \end{cases}\label{eq:j1wa2}
    \end{align}
\end{enumerate}
\end{prop}
Let us first show that this proposition implies desired Proposition \ref{lm:resolv2} in all cases except $0<\zeta<1$ (covered in the previous section). Let $2m<\zeta<2m+1$ for some integer $m\ge 1$, then, using Eqs. \eqref{eq:ame1}, \eqref{eq:ame2}, \begin{align}
    \langle \widetilde A^{-m}\mathbf e_1,\widetilde A^{-m}(\widetilde A+\epsilon)^{-1}\mathbf e_1\rangle ={}&\sum_{n=1}^\infty (\widetilde A^{-m}\mathbf e_1)_n (\widetilde A^{-m}(\widetilde A+\epsilon)^{-1}\mathbf e_1)_n\\
    ={}& \sum_{n=1}^{n_0}+\sum_{n=n_0+1}^\infty \\
    ={}& \sum_{n=1}^{n_0} O(n^{\frac{-1-(2+\nu)(\zeta-2m)}{2}}) O(n^{\frac{-1-(2+\nu)(\zeta-2m-2)}{2}})\\
    &+\sum_{n=n_0+1}^\infty O(n^{\frac{-1-(2+\nu)(\zeta-2m)}{2}})O(\epsilon^{-1}n^{\frac{-1-(2+\nu)(\zeta-2m)}{2}})\\
    ={}&\sum_{n=1}^{n_0} O(n^{-1-(2+\nu)(\zeta-2m-1)})+\sum_{n=n_0+1}^\infty O(\epsilon^{-1}n^{-1-(2+\nu)(\zeta-2m)})\\
    ={}&O(n_0^{-(2+\nu)(\zeta-2m-1)})\label{eq:n02nu}\\
    ={}&O(\epsilon^{\zeta-2m-1}),
\end{align}
which is the desired bound \eqref{eq:resolv21}. Note that here we used both inequalities $2m<\zeta<2m+1$ and the identity $\epsilon=h^{2+\nu}$ to get Eq. \eqref{eq:n02nu}. 

By a similar reasoning, if $2m+1<\zeta<2m+2$ with some $m\ge 0$, then Eqs. \eqref{eq:j1wa1}, \eqref{eq:j1wa2} imply desired Eq. \eqref{eq:resolv22} of Proposition \ref{lm:resolv2}. We have  thus fully proved Proposition \ref{lm:resolv2} assuming Proposition \ref{lm:resolvseq}, and it remains to prove the latter.

\medskip
\begin{proof} We prove Proposition \ref{lm:resolvseq} by induction. The base of induction is Statement 2 with $m=0$ (corresponding to $\zeta>1).$ In the induction step, we either derive Statement 1 for $m$ from Statement 2 for $m-1$, or derive Statement 2 for $m$ from Statement 1 with the same $m$. 
\paragraph{Base of induction: Statement 2 for $m=0$.} 
 Given any $\mathbf u\in l^2$, denote $\mathbf w= J^{-1}\mathbf u.$  If $\mathbf w\in l^2$, its components satisfy the equations
 \begin{align}
    w_1 ={}& u_1,\\
    n^{-\frac{\nu}{2}}w_n-(\tfrac{n}{n-1})^{\frac{1-(2+\nu)\zeta}{2}}(n-1)^{-\frac{\nu}{2}}w_{n-1}={}&u_n, \quad n= 2,3,\ldots\label{eq:nwn}
\end{align}
The system can be solved iteratively, starting from $w_1$ and computing $w_n$ from $w_{n-1}$ using Eq.\eqref{eq:nwn}:
\begin{align}
    w_n={}&(\tfrac{n}{n-1})^{\frac{1+\nu-(2+\nu)\zeta)}{2}}w_{n-1}+n^{\frac{\nu}{2}}u_n\\
    ={}&(\tfrac{n}{n-1})^{\frac{1+\nu-(2+\nu)\zeta)}{2}}\Big((\tfrac{n-1}{n-2})^{\frac{1+\nu-(2+\nu)\zeta)}{2}}w_{n-2}+(n-1)^{\frac{\nu}{2}}u_{n-1}\Big)+n^{\frac{\nu}{2}}u_n,\\
    ={}&n^{\frac{-1-(2+\nu)(\zeta-1)}{2}}\sum_{k=1}^n k^{\frac{(2+\nu)\zeta-1}{2}} u_k,\quad
    \quad n= 2,3,\ldots\label{eq:wnj}
\end{align}
In the special case $\mathbf u=\mathbf e_1$ we get the explicit solution
\begin{align}
    w_n =n^{\frac{-1-(2+\nu)(\zeta-1)}{2}},
\end{align}
proving desired Eq. \eqref{eq:j1wa1} for $m=0$. It is also clear that this $\mathbf w\in l^2$ as long as $\zeta>1.$ 

Now let $\mathbf u=(\widetilde A+\epsilon)^{-1}\mathbf e_1$. Let us first bound the components $u_n$, using results of Section \ref{sec:lmresolv}  with $\mathbf x=\mathbf u$ and Lemma \ref{lm:zeta1theta}. First we observe that $u_1$ is uniformly bounded for all suficiently small $\epsilon$: by Eq. \eqref{eq:x1theta} and Lemma \ref{lm:zeta1theta}, as long as $\zeta>1$,
\begin{align}
    |u_1|={}&(\epsilon-\tfrac{h\theta_1^*}{1-h\theta_1^*})^{-1}\\
    ={}&(\epsilon+1-\tfrac{1}{1-h\theta_1^*})^{-1}\\
    \le{}&(1-\tfrac{1}{1+(\nu+2)(\zeta-1)})^{-1}<\infty.
\end{align}
Next we obtain a bound on $u_n$ for $n\le n_0.$ Using the definition of $\theta^*$ and Lemma \ref{lm:zeta1theta},
\begin{align}
    |x_n|={}& n^{-g}|y_n|\\
    ={}&n^{-g}|y_1|\prod_{k=1}^{n-1}(1-h\theta^*_k)^{-1} \\
    \le{}&n^{\frac{-1+(2+\nu)\zeta}{2}}|x_1|\prod_{k=1}^{n-1}\Big(1-\tfrac{(\nu+2)(1-\zeta)}{k}\Big)^{-1}\\
    ={}& n^{\frac{-1+(2+\nu)\zeta}{2}} O(n^{(\nu+2)(1-\zeta)})\\
    ={}& O(n^{\frac{-1-(\nu+2)(\zeta-2)}{2}}).
\end{align}
Now using Eq. \eqref{eq:wnj}, we get for $n\le n_0$
\begin{align}
    |((\widetilde A+\epsilon)^{-1}\mathbf e_1)_n| ={}& n^{\frac{-1-(2+\nu)(\zeta-1)}{2}}\sum_{k=1}^n k^{\frac{(2+\nu)\zeta-1}{2}} O(k^{\frac{-1-(\nu+2)(\zeta-2)}{2}})\\
    ={}&n^{\frac{-1-(2+\nu)(\zeta-1)}{2}}\sum_{k=1}^n O(k^{-1+\nu+2})\\
    ={}&n^{\frac{-1-(2+\nu)(\zeta-1)}{2}} O(n^{\nu+2})\\
    ={}&O(n^{\frac{-1-(2+\nu)(\zeta-3)}{2}}),
\end{align}
which is the desired Eq. \eqref{eq:j1wa1}.

Now consider $n>n_0.$ By Lemma \ref{lm:thetasab}, we can assume w.l.o.g. (if necessary, increasing $s_0$) that $\theta_n^*<-1$ for $n\ge n_0$. Then, with $\mathbf u = (\widetilde A+\epsilon)^{-1}\mathbf e_1$ and $h$ small enough,
\begin{align}
    |u_n|={}& n^{-g}|u_{n_0}|n_0^{g}\prod_{k=n_0}^{n-1}(1-h\theta^*_k)^{-1} \\
    ={}&O\Big(n^{-g}n_0^{(\nu+2)(1-\zeta)}(1+h)^{n_0-n}\Big) \\
    ={}&O\Big(h^{(\nu+2)(\zeta-1)+g}(nh)^{-g}(1+h)^{-n}\Big) \\
    ={}&O\Big(h^{(\nu+2)(\zeta-1)+g}e^{-nh/2}\Big).
\end{align}
It follows that for $n>n_0$
\begin{align}
    |(J^{-1}(\widetilde A+\epsilon)^{-1}\mathbf e_1)_n|\le{}& n^{\frac{-1-(2+\nu)(\zeta-1)}{2}} \sum_{k=1}^\infty  k^{\frac{(2+\nu)\zeta-1}{2}} |u_k|\\
    ={}&  n^{\frac{-1-(2+\nu)(\zeta-1)}{2}} O\Big(n_0^{\nu+2}+ h^{(\nu+2)(\zeta-1)+g}\sum_{k=n_0}^\infty  k^{\frac{(2+\nu)\zeta-1}{2}} e^{-kh/2}\Big)\\
    ={}&  n^{\frac{-1-(2+\nu)(\zeta-1)}{2}} O\Big(h^{-(\nu+2)}+ h^{(\nu+2)(\zeta-1)+2g}\sum_{k=n_0}^\infty  (kh)^{\frac{(2+\nu)\zeta-1}{2}} e^{-kh/2}\Big)\\
    ={}&  n^{\frac{-1-(2+\nu)(\zeta-1)}{2}} O\Big(h^{-(\nu+2)}+ h^{(\nu+2)(\zeta-1)+2g-1}\Big)\\
    ={}&  n^{\frac{-1-(2+\nu)(\zeta-1)}{2}}O(\epsilon^{-1}),
\end{align}
which is the desired bound \eqref{eq:j1wa2}. 

\paragraph{Induction step: Statement 1 for $m$ from Statement 2 for $m-1$.} Note that $(\widetilde A)^{-1}=(J^{\dagger})^{-1}J^{-1}$ and so we can represent
\begin{align}
    \widetilde A^{-m}\mathbf e_1 ={}& (J^{\dagger})^{-1}(J^{-1}\widetilde A^{-(m-1)}\mathbf e_1),\\
    \widetilde A^{-m}(\widetilde A+\epsilon)^{-1}\mathbf e_1={}&(J^{\dagger})^{-1}(J^{-1}\widetilde A^{-(m-1)}(\widetilde A+\epsilon)^{-1}\mathbf e_1).
\end{align}
Let us examine the operator $(J^{\dagger})^{-1}$. Given any $\mathbf u\in l^2$, denote $\mathbf w= (J^{\dagger})^{-1}\mathbf u.$  If $\mathbf w\in l^2$, its components satisfy the equations
\begin{equation}\label{eq:nwu2wn}
    n^{-\frac{\nu}{2}}(w_n-(\tfrac{n+1}{n})^{\frac{1-(2+\nu)\zeta}{2}}w_{n+1}) = u_n,\quad n=1,2,\ldots
\end{equation}
These equations can be solved iteratively, with $w_n$ expressed via $w_{n+1}:$
\begin{align}
    w_n ={}& n^{\frac{-1+(2+\nu)\zeta}{2}} \Big(n^{\frac{1+\nu-(2+\nu)\zeta}{2}}u_n+(n+1)^{\frac{1-(2+\nu)\zeta}{2}}w_{n+1}\Big)\\
    ={}&n^{\frac{-1+(2+\nu)\zeta}{2}}\Big(\sum_{m=n}^{q-1} k^{\frac{-1-(2+\nu)(\zeta-1)}{2}}u_k+q^{\frac{1-(2+\nu)\zeta}{2}}w_{q}\Big)
\end{align}
for any $q> n$. It is convenient to take the limit $q\to \infty$. If we assume that 
\begin{align}\label{eq:mnu1}
    \sum_{k=1}^{\infty} k^{\frac{-1-(2+\nu)(\zeta-1)}{2}}|u_k|<\infty
\end{align}
and
\begin{equation}\label{eq:mnu2}
    w_q = o(q^{\frac{-1+(2+\nu)\zeta}{2}}),
\end{equation}
then we can take this limit, obtaining
\begin{align}\label{eq:wnn1}
    w_n =n^{\frac{-1+(2+\nu)\zeta}{2}}\sum_{k=n}^{\infty} k^{\frac{-1-(2+\nu)(\zeta-1)}{2}}u_k.
\end{align}
In fact, if we just assume condition \eqref{eq:mnu1} and define $w_n$ by Eq. \eqref{eq:wnn1}, these $w_n$ clearly satisfy equations \eqref{eq:nwu2wn} and condition \eqref{eq:mnu2}. Accordingly, it suffices to only check condition \eqref{eq:mnu1}.

We now apply this expansion to $\mathbf u=J^{-1}\widetilde A^{-(m-1)}\mathbf e_1$ and $\mathbf u=J^{-1}\widetilde A^{-(m-1)}(\widetilde A+\epsilon)^{-1}\mathbf e_1.$  Observe first that in both cases condition \eqref{eq:mnu1} is fulfilled thanks to induction hypotheses \eqref{eq:j1wa1}, \eqref{eq:j1wa2} for $m-1$, since they imply 
\begin{align}
    k^{\frac{-1-(2+\nu)(\zeta-1)}{2}}|u_k| ={}& k^{\frac{-1-(2+\nu)(\zeta-1)}{2}} O(k^{\frac{-1-(2+\nu)(\zeta-2m+1)}{2}})\\ ={}& O(k^{-1-(2+\nu)(\zeta-m)})
\end{align}
and, by assumption, $\zeta>2m\ge m$. Taking $\mathbf u=J^{-1}\widetilde A^{-(m-1)}\mathbf e_1$, we obtain desired Eq. \eqref{eq:ame1}:
\begin{align}
    (\widetilde A^{-m}\mathbf e_1)_n ={}& n^{\frac{-1+(2+\nu)\zeta}{2}}\sum_{k=n}^{\infty} k^{\frac{-1-(2+\nu)(\zeta-1)}{2}}(J^{-1}\widetilde A^{-(m-1)}\mathbf e_1)_k\\
    ={}& n^{\frac{-1+(2+\nu)\zeta}{2}}\sum_{k=n}^{\infty} O(k^{-1-(2+\nu)(\zeta-m)})\\
    ={}& n^{\frac{-1+(2+\nu)\zeta}{2}} O(n^{-(2+\nu)(\zeta-m)})\\
    ={}& O(n^{\frac{-1-(2+\nu)(\zeta-2m)}{2}}).
\end{align}
Taking $\mathbf u=J^{-1}\widetilde A^{-(m-1)}(\widetilde A+\epsilon)^{-1}\mathbf e_1$, in the case $n>n_0$ we obtain desired Eq. \eqref{eq:ame2} by a completely similar argument. In the case $n\le n_0$ we obtain the desired bound by
\begin{align}
    (\widetilde A^{-m}(\widetilde A+\epsilon)^{-1}\mathbf e_1)_n ={}& n^{\frac{-1+(2+\nu)\zeta}{2}}\sum_{k=n}^{\infty} k^{\frac{-1-(2+\nu)(\zeta-1)}{2}}(J^{-1}\widetilde A^{-(m-1)}(\widetilde A+\epsilon)^{-1}\mathbf e_1)_k\\
    ={}& n^{\frac{-1+(2+\nu)\zeta}{2}}\Big(\sum_{k=n}^{n_0} k^{\frac{-1-(2+\nu)(\zeta-1)}{2}} O(k^{\frac{-1-(2+\nu)(\zeta-2m-1)}{2}})\\
    &+\sum_{k=n_0+1}^\infty k^{\frac{-1-(2+\nu)(\zeta-1)}{2}} O(\epsilon^{-1}k^{\frac{-1-(2+\nu)(\zeta-2m+1)}{2}})\Big)\\
    ={}& n^{\frac{-1+(2+\nu)\zeta}{2}}\Big(\sum_{k=n}^{n_0} O(k^{-1-(2+\nu)(\zeta-m-1)}) +O(\epsilon^{-1}n_0^{-(2+\nu)(\zeta-m)})\Big)\\
    ={}& n^{\frac{-1+(2+\nu)\zeta}{2}} \Big(O(n^{-(2+\nu)(\zeta-m-1)}-n_0^{-(2+\nu)(\zeta-m-1)})\\
    &+O(n_0^{-(2+\nu)(\zeta-m-1)})\Big)\\
    ={}& n^{\frac{-1+(2+\nu)\zeta}{2}}O(n^{-(2+\nu)(\zeta-m-1)})\\
    ={}&O(n^{\frac{-1-(2+\nu)(\zeta-2m-2)}{2}}),
\end{align}
where we used the fact that $m\ge 1$ and $\zeta>2m\ge m+1$. 

\paragraph{Induction step: Statement 2 for $m$ from Statement 1 for the same $m$.} Applying again Eq. \eqref{eq:wnj} with $\mathbf u=\widetilde {\mathcal A}^{-m}\mathbf e_1$, we get for $\zeta>2m+1$ with $m\ge 1$
\begin{align}
        |(J^{-1}\widetilde A^{-m}\mathbf e_1)_n| \le{}&n^{\frac{-1-(2+\nu)(\zeta-1)}{2}}\sum_{k=1}^n k^{\frac{(2+\nu)\zeta-1}{2}} |u_k|\\
={}&n^{\frac{-1-(2+\nu)(\zeta-1)}{2}}\sum_{k=1}^n k^{\frac{(2+\nu)\zeta-1}{2}} O(k^{\frac{-1-(2+\nu)(\zeta-2m)}{2}})\\  
={}&n^{\frac{-1-(2+\nu)(\zeta-1)}{2}}\sum_{k=1}^n  O(k^{-1+(2+\nu)m})\\ 
={}&n^{\frac{-1-(2+\nu)(\zeta-1)}{2}} O(n^{(2+\nu)m})\\ 
={}&O(n^{\frac{-1-(2+\nu)(\zeta-2m-1)}{2}}), 
\end{align}
which is the desired bound \eqref{eq:j1wa1}. The $n\le n_0$ case of bound \eqref{eq:j1wa2} is obtained similarly. In the case $n>n_0$ we have
\begin{align}
        |(J^{-1}\widetilde A^{-m}\mathbf e_1)_n| \le{}&n^{\frac{-1-(2+\nu)(\zeta-1)}{2}}\Big(\sum_{k=1}^{n_0}k^{\frac{(2+\nu)\zeta-1}{2}} |u_k|+\sum_{k=n_0+1}^n k^{\frac{(2+\nu)\zeta-1}{2}} |u_k|\Big)\\
={}&n^{\frac{-1-(2+\nu)(\zeta-1)}{2}}\Big(\sum_{k=1}^{n_0}k^{\frac{(2+\nu)\zeta-1}{2}} O(k^{\frac{-1-(2+\nu)(\zeta-2m-2)}{2}})\\
&+ \sum_{k=n_0+1}^n k^{\frac{(2+\nu)\zeta-1}{2}} O(\epsilon^{-1}k^{\frac{-1-(2+\nu)(\zeta-2m)}{2}})\Big)\\  
={}&n^{\frac{-1-(2+\nu)(\zeta-1)}{2}}\Big(O(n_0^{(2+\nu)(m+1)}) + O(\epsilon^{-1}\sum_{k=n_0+1}^n k^{-1+(2+\nu)m})\Big)\\ 
={}&n^{\frac{-1-(2+\nu)(\zeta-1)}{2}}\Big(O(\epsilon^{-1}n_0^{(2+\nu)m}) + O(\epsilon^{-1}(n^{(2+\nu)m}-n_0^{(2+\nu)m}))\Big)\\ 
={}&n^{\frac{-1-(2+\nu)(\zeta-1)}{2}} O(\epsilon^{-1}n^{(2+\nu)m})\\ 
={}&O(\epsilon^{-1}n^{\frac{-1-(2+\nu)(\zeta-2m-1)}{2}}),
\end{align}
which is the desired bound \eqref{eq:j1wa2} for $n>n_0.$ This completes the proof of the proposition. 
\end{proof}

\section{Experiments}

\subsection{Details of experiments}\label{sec:exp_details}

\paragraph{Algorithms.} Let us describe details of each of the eight algorithms present in our experiments (see the legend of Figure \ref{fig:Losses}). 

For ``constant rate GD'' and ``constant rate HB'' we used parameters $\alpha=1$ and $\beta=0.9$. 

The algorithm ``scheduled HB'' uses schedule \eqref{eq:Jacobi_ansatz_parameters} for $\alpha_n,\beta_n$ with parameters $a=\zeta, b=0, r=1$. The ``asymptotic scheduled HB'' uses asymptotic $n\to\infty$ version of Jacobi schedule given by the rightmost part of \eqref{eq:Jacobi_ansatz_parameters}, and additionally set $\alpha_n=1$ (it is rather unnecessary artifact of our experimentation). The difference between limiting values of $\alpha_n$ in "scheduled HB" ($\alpha_n \rightarrow 2$) and "asymptotic scheduled HB" ($\alpha_n=1$) explains the slight advantage of the former in Figure \ref{fig:Losses}. 

The algorithm "scheduled GD" uses the schedule based on the roots of the same residual polynomials we used for "scheduled HB", as described in the proof of theorem \ref{th:ubgdnonconst} in section \ref{sec:ubgdpredefined}. As this schedule continuously "fills" roots of Jacobi polynomials of degrees $2^l$, we see respective stair-like structure in figure \ref{fig:Losses} and spikes in figure \ref{fig:MNIST_net}.    

The adaptive algorithms "Steepest Descent" and "basic CG" in our experiments are given by formulas \eqref{eq:SD_quadratic_choice} and \eqref{eq:CG_quadratic_choice_1},\eqref{eq:CG_quadratic_choice_2} since we apply them only to quadratic problems. 

The "numerically stable CG" algorithm is meant to fix the problems of "basic CG" as we expect the convergence rate for CG to be $\sim n^{-(2+\nu)\zeta}$. In particular,  numerical errors accumulate during the run of the algorithm, leading to non-exact placement of the roots of respective residual polynomial $p_n(\lambda)$ at spectral points $\lambda_k$. We resolve this issue by introducing some kind of checking procedure on each step. As CG is known to produce a system of orthogonal steps $\Delta{\mathbf{w}}_n = \mathbf{w}_{n+1}-\mathbf{w}_{n}$ for quadratic problems, we directly check this orthogonality on each step. Specifically, before making new step  $\Delta \mathbf{w}_{n+1}$, we first eliminate all its components along previously made steps $\Delta \mathbf{w}_l, l\leq n$ (which we store during the run of the algorithm). Then, after $\Delta \mathbf{w}_{n+1}$ is made orthogonal to all $\Delta \mathbf{w}_l, l\leq n$, we correct the magnitude of the step to fully eliminate the component of $\mathbf{f}_{n+1}-\mathbf{f}_*$ in this direction. The described procedure is equivalent to formulas \eqref{eq:CG_quadratic_choice_1}, \eqref{eq:CG_quadratic_choice_2} in exact arithmetic, but is required for actual implementation of CG to reach convergence rate $\sim n^{-(2+\nu)\zeta}$, as can be seen from experiments in Figure \ref{fig:Losses}. 

\paragraph{The MNIST-based quadratic problem.} In figures \ref{fig:MNIST_loss_and_spectral_distributions} and \ref{fig:Losses} we took the first $M=30000$ MNIST images from the usual train subset, and flattened them into $d=784$--dimensional vectors $\{\mathbf{x}_i\}_{i=1}^M = \mathcal{D}$. Then we normalize each vector using the dataset mean $\mathbf{m}=\frac{1}{M}\sum_i \mathbf{x}_i$ and variance $r^2 = \frac{1}{M}\sum_i \|\mathbf{x_i} -\mathbf{m}\|^2$ by $\mathbf{x}_i \rightarrow (\mathbf{x}_i-\mathbf{m}) / r$. Then the scalar targets $y_i$ were obtained simply as numerical values $\{0, 1, \ldots, 8, 9\}$ of the digits corresponding to the images $\mathbf{x}_i$. Then, instead of formulating the quadratic problem in parameter space $\mathbf{w}$, where we would need to specify matrix $A$ and vector $\mathbf{b}$, we consider the problem in output space where we need matrix $\widetilde{A}$ and vector $\mathbf{f}_*$. Components of the latter are given simply by our targets $y_i$, and the matrix $\widetilde{A}$ is obtained by evaluating a kernel $K(\mathbf{x},\mathbf{x}')$ on our data with $\widetilde{A}_{ij}=K(\mathbf{x}_i,\mathbf{x}_j)$. For the kernel we take the NTK of infinitely wide shallow ReLU network given by (see e.g. \cite{NEURIPS2019_0d1a9651})
\begin{equation}
    K(\mathbf{x},\mathbf{x}') = \frac{\|\mathbf{x}\|\|\mathbf{x}'\|\big(\sin \varphi + 2 \cos \varphi (\pi-\varphi) \big)}{2\pi}, \quad \cos\varphi = \frac{\mathbf{x}^T\mathbf{x}'}{\|\mathbf{x}\|\|\mathbf{x}'\|}.
\end{equation}

For figure \ref{fig:MNIST_UBs} we repeat the same procedure but on a full training set of MNIST ($M=5000$), which was possible due to the availability of additional computational resources at the later times of our work on this paper. Also, in figure \ref{fig:MNIST_UBs} we changed the kernel from NTK to sigmoid kernel $K(\mathbf{x},\mathbf{x}')=\tanh(\mathbf{x}^T\mathbf{x}'+1)$, which seem to better illustrate the described phenomenon.

\paragraph{The neural network experiment.} We consider the standard MNIST classification problem with one-hot encoding of the 10 classes. We consider a simple shallow ReLU network of width $N=1000$ with the NTK parametrization. Its function $f(\mathbf{x})$ can be written as
\begin{equation}
    f(\mathbf{x}) = \frac{1}{\sqrt{N}} \sum\limits_{l=1}^N \mathbf{c}_l \operatorname{ReLU}(\mathbf{w}_l^T \mathbf{x} + b_l)
\end{equation}
where $\mathbf{w}_l \in \mathbb{R}^{784}, \; b_k \in \mathbb{R}, \; \mathbf{c}_l \in \mathbb{R}^{10}$ are the parameters of the neuron $l$.

Then we train this network on the full MNIST dataset with standard train-test split. Importantly, we don't use mini-batches during training steps, but process the whole train dataset of size $M=50000$ during optimization. Thus, considering full-batch GD allows us to stay close to our main setting, with the only difference being non-linearity of the model.

\subsection{Finding the end of the loss power law region} \label{sec:powerlaw_end} 
Let us formulate a general principle allowing to estimate the transition point $n_{\text{th}}$ on the loss curve where the power-law region ends. The end of the power-law region in the loss is due to the end of the power law region in the spectral measure asymptotic. In particular, we assume that the power-law asymptotic $\rho[(\lambda_1,\lambda_2)]\sim \lambda_2^\zeta-\lambda_1^\zeta$ holds in the region $\lambda_1,\lambda_2 \gtrsim \lambda_{\text{low}}$ with $\lambda_{\text{low}}$ being (an estimated) end of this power-law region. For synthetic data from figure \ref{fig:Losses} (a) this would be simply the lowest eigenvalue $\lambda_{\text{low}} = M^{-\nu}$; for the MNIST-based quadratic problem from figure \ref{fig:Losses} (b-d) we visually set $\lambda_{\text{low}} = 5 \times 10^{-5}$. Finally, $n_{\text{th}}$ is simply a step when for a chosen optimization algorithm the region $[0,\lambda_{\text{low}}]$ can no longer be ignored.  

Next, define (approximately) a point $\lambda(n)$ as the point where the residual polynomial $p_n(\lambda)$ of considered optimization algorithm starts to significantly deviate from its value at the origin $p_n(\lambda=0)=1$, and then is expected to rapidly converge to zero $p_n(\lambda)\rightarrow 0$ as $\lambda \gg \lambda(n)$. Then the loss of the algorithm at step $n$ can be estimated as $L_n \sim \rho([0,\lambda(n)])$. Suppose that only the $\rho([\lambda_{\text{low}}, \lambda(n)])$ part of this loss is defined by asymptotic spectral power law, while the $\rho([0, \lambda_{\text{low}}])$ part is unknown. Then the fraction $h$ of "controlled" loss on step $n$ can be estimated as
\begin{equation}
h(n) = \frac{\lambda(n)^\zeta - \lambda_{\text{low}}^\zeta}{\lambda(n)^\zeta} = 1 - \left(\frac{\lambda(n)}{\lambda_{\text{low}}}\right)^{-\zeta}     
\end{equation}
As this fraction reaches some predefined tolerance threshold $h_0$ (e.g. $h_0=0.5$) we can say that the loss no longer follows its power-law and therefore we are at threshold step $n_\text{th}$. Formally, $n_\text{th}$ is defined by the equation
\begin{equation}\label{eq:n_th_eq}
    h_0 = 1-\left(\frac{\lambda(n_{\text{th}})}{\lambda_{\text{low}}}\right)^{-\zeta}   
\end{equation}

To actually apply this principle we need to know $\lambda(n)$. Let us find it for the algorithms considered in this work. For algorithms with constant rates $\alpha, \beta$ the residual polynomial at small $\lambda$ has large $n$ asymptotic $p_n(\lambda) \sim \exp(-\frac{n \alpha \lambda}{1-\beta})$. Indeed, for $\beta=0$ we have $p_n(\lambda)=(1-\alpha\lambda)^n\approx \exp(-\alpha n)$, while for the case with momentum, one needs to use representation \eqref{eq:const_lr_HB_poly} together with \eqref{eq:HB_Chebyshev_explicit_form} near $\lambda=0$ ($z=\tfrac{1+\beta}{2\sqrt{\beta}}$). Thus we define $\lambda(n) = \frac{1-\beta}{\alpha n}$. Next, for the algorithms with predefined schedules based on Jacobi polynomials we recall asymptotic \eqref{eq:Jacobi_square_asym}, which says that the polynomial start to deviate from $1$ at $\lambda n^2 \sim 1$, therefore $\lambda(n) = n^{-2}$. Finally, for (stable) Conjugate Gradient method we may assume that $p_n(\lambda)$ is simply $q_{n/2}^{(a,b)}(\lambda / \lambda_{n/2})$ as in the proof of \ref{th:cgdiscr}. Then we again apply Jacobi polynomials asymptotic to find that $\lambda(n) = n^{-\nu-2}$. To summarize, we have established
\begin{equation}\label{eq:lambda_n_vals}
    \lambda(n) = 
    \begin{cases}
    \frac{1-\beta}{\alpha n}, \qquad& \text{for constant rate algorithms}\\
    \frac{1}{n^2}, \qquad& \text{for algorithms based on Jacobi polynomials}\\
    n^{-\nu-2}, \qquad& \text{for (numerically stable) Conjugate Gradients}
    \end{cases}
\end{equation}
Solving \eqref{eq:n_th_eq} with \eqref{eq:lambda_n_vals} gives

\begin{equation}
    n_\text{th} = 
    \begin{cases}
    \vspace{2mm}
    (1-h_0)^{\frac{1}{\zeta}}\frac{1-\beta}{\alpha \lambda_{\text{low}}}, \qquad& \text{for constant rate algorithms}\\
    \vspace{2mm}
    (1-h_0)^{\frac{1}{\zeta}}\frac{1}{\sqrt{\lambda_{\text{low}}}}, \qquad& \text{for algorithms based on Jacobi polynomials}\\
    \vspace{2mm}
    (1-h_0)^{\frac{1}{\zeta}}\lambda_{\text{low}}^{-\frac{1}{\nu+2}}, \qquad& \text{for (numerically stable) Conjugate Gradients}
    \end{cases}
\end{equation}
This result agrees with Table \ref{tab:bounds} and also admits the following interpretation: the critical values $n_{\text{th}}$ approximately correspond to the step numbers at which the order of the loss magnitude approximately matches the value $\lambda_{\text{low}}^\zeta$ associated with the measure $\rho((0,\lambda_{\text{low}}])$ under the power-law spectral assumption.

\addcontentsline{toc}{section}{References}

\vskip 0.2in
\bibliography{jmlr1}

\begin{thebibliography}{68}
\providecommand{\natexlab}[1]{#1}
\providecommand{\url}[1]{\texttt{#1}}
\expandafter\ifx\csname urlstyle\endcsname\relax
  \providecommand{\doi}[1]{doi: #1}\else
  \providecommand{\doi}{doi: \begingroup \urlstyle{rm}\Url}\fi

\bibitem[Akaike(1959)]{akaike1959successive}
Hirotugu Akaike.
\newblock On a successive transformation of probability distribution and its
  application to the analysis of the optimum gradient method.
\newblock \emph{Annals of the Institute of Statistical Mathematics},
  11\penalty0 (1):\penalty0 1--16, 1959.

\bibitem[Atanasov et~al.(2021)Atanasov, Bordelon, and
  Pehlevan]{atanasov2021neural}
Alexander Atanasov, Blake Bordelon, and Cengiz Pehlevan.
\newblock Neural networks as kernel learners: The silent alignment effect.
\newblock 2021.
\newblock \doi{10.48550/ARXIV.2111.00034}.
\newblock URL \url{https://arxiv.org/abs/2111.00034}.

\bibitem[Bahri et~al.(2021)Bahri, Dyer, Kaplan, Lee, and
  Sharma]{bahri2021explaining}
Yasaman Bahri, Ethan Dyer, Jared Kaplan, Jaehoon Lee, and Utkarsh Sharma.
\newblock Explaining neural scaling laws.
\newblock \emph{arXiv preprint arXiv:2102.06701}, 2021.

\bibitem[Basri et~al.(2020)Basri, Galun, Geifman, Jacobs, Kasten, and
  Kritchman]{basri2020frequency}
Ronen Basri, Meirav Galun, Amnon Geifman, David Jacobs, Yoni Kasten, and Shira
  Kritchman.
\newblock Frequency bias in neural networks for input of non-uniform density,
  2020.
\newblock URL \url{https://arxiv.org/abs/2003.04560}.

\bibitem[Berthier et~al.(2020{\natexlab{a}})Berthier, Bach, and
  Gaillard]{Berthier2020_gossip}
Rapha\"{e}l Berthier, Francis Bach, and Pierre Gaillard.
\newblock Accelerated gossip in networks of given dimension using jacobi
  polynomial iterations.
\newblock \emph{SIAM Journal on Mathematics of Data Science}, 2\penalty0
  (1):\penalty0 24--47, 2020{\natexlab{a}}.
\newblock \doi{10.1137/19M1244822}.
\newblock URL \url{https://doi.org/10.1137/19M1244822}.

\bibitem[Berthier et~al.(2020{\natexlab{b}})Berthier, Bach, and
  Gaillard]{berthier2020tight}
Rapha{\"e}l Berthier, Francis Bach, and Pierre Gaillard.
\newblock Tight nonparametric convergence rates for stochastic gradient descent
  under the noiseless linear model.
\newblock \emph{arXiv preprint arXiv:2006.08212}, 2020{\natexlab{b}}.

\bibitem[Bietti(2021)]{bietti2021approximation}
Alberto Bietti.
\newblock Approximation and learning with deep convolutional models: a kernel
  perspective, 2021.
\newblock URL \url{https://arxiv.org/abs/2102.10032}.

\bibitem[Birman and Solomjak(1970)]{Birman_1970}
M~{\v{S}} Birman and M~Z Solomjak.
\newblock Asymptotic behavior of the spectrum of weakly polar integral
  operators.
\newblock \emph{Mathematics of the {USSR}-Izvestiya}, 4\penalty0 (5):\penalty0
  1151--1168, oct 1970.
\newblock \doi{10.1070/im1970v004n05abeh000948}.
\newblock URL \url{https://doi.org/10.1070/im1970v004n05abeh000948}.

\bibitem[Birman and Solomjak(2012)]{Birman1987SpectralTO}
M.S. Birman and M.Z. Solomjak.
\newblock \emph{{Spectral Theory of Self-Adjoint Operators in Hilbert Space}}.
\newblock Mathematics and its Applications. Springer Netherlands, 2012.
\newblock ISBN 9789400945869.
\newblock URL \url{https://books.google.ru/books?id=unPrCAAAQBAJ}.

\bibitem[Bj{\"o}rck et~al.(1998)Bj{\"o}rck, Elfving, and
  Strakos]{bjorck1998stability}
{\AA}ke Bj{\"o}rck, Tommy Elfving, and Zdenek Strakos.
\newblock Stability of conjugate gradient and lanczos methods for linear least
  squares problems.
\newblock \emph{SIAM Journal on Matrix Analysis and Applications}, 19\penalty0
  (3):\penalty0 720--736, 1998.

\bibitem[Bordelon and Pehlevan(2021)]{bordelon2021learning}
Blake Bordelon and Cengiz Pehlevan.
\newblock Learning curves for sgd on structured features, 2021.
\newblock URL \url{https://arxiv.org/abs/2106.02713}.

\bibitem[Brakhage(1987)]{brakhage1987ill}
Helmut Brakhage.
\newblock On ill-posed problems and the method of conjugate gradients.
\newblock In \emph{Inverse and ill-posed Problems}, pages 165--175. Elsevier,
  1987.

\bibitem[Canatar et~al.(2021)Canatar, Bordelon, and
  Pehlevan]{canatar2021spectral}
Abdulkadir Canatar, Blake Bordelon, and Cengiz Pehlevan.
\newblock Spectral bias and task-model alignment explain generalization in
  kernel regression and infinitely wide neural networks.
\newblock \emph{arXiv preprint arXiv:2006.13198}, 2021.

\bibitem[Caponnetto and De~Vito(2007)]{caponnetto2007optimal}
Andrea Caponnetto and Ernesto De~Vito.
\newblock Optimal rates for the regularized least-squares algorithm.
\newblock \emph{Foundations of Computational Mathematics}, 7\penalty0
  (3):\penalty0 331--368, 2007.

\bibitem[Chizat et~al.(2019)Chizat, Oyallon, and Bach]{chizat2019lazy}
L{\'e}na{\"\i}c Chizat, Edouard Oyallon, and Francis Bach.
\newblock On lazy training in differentiable programming.
\newblock \emph{Advances in Neural Information Processing Systems},
  32:\penalty0 2937--2947, 2019.

\bibitem[Cho and Saul(2009)]{NIPS2009_kernel}
Youngmin Cho and Lawrence Saul.
\newblock Kernel methods for deep learning.
\newblock In Y.~Bengio, D.~Schuurmans, J.~Lafferty, C.~Williams, and
  A.~Culotta, editors, \emph{Advances in Neural Information Processing
  Systems}, volume~22, pages 342--350. Curran Associates, Inc., 2009.
\newblock URL
  \url{https://proceedings.neurips.cc/paper/2009/file/5751ec3e9a4feab575962e78e006250d-Paper.pdf}.

\bibitem[Cui et~al.(2021)Cui, Loureiro, Krzakala, and
  Zdeborová]{cui2021generalization}
Hugo Cui, Bruno Loureiro, Florent Krzakala, and Lenka Zdeborová.
\newblock Generalization error rates in kernel regression: The crossover from
  the noiseless to noisy regime, 2021.
\newblock URL \url{https://arxiv.org/abs/2105.15004}.

\bibitem[Daniel(1971)]{daniel1971approximate}
J.W. Daniel.
\newblock \emph{The Approximate Minimization of Functionals}.
\newblock Prentice-Hall series in automatic computation. Prentice-Hall, 1971.
\newblock ISBN 9780130438775.
\newblock URL \url{https://books.google.ru/books?id=kEUZAQAAIAAJ}.

\bibitem[Dou and Liang(2021)]{dou2021training}
Xialiang Dou and Tengyuan Liang.
\newblock Training neural networks as learning data-adaptive kernels: Provable
  representation and approximation benefits.
\newblock \emph{Journal of the American Statistical Association}, 116\penalty0
  (535):\penalty0 1507--1520, 2021.

\bibitem[Feller(1991)]{feller1991introduction}
William Feller.
\newblock \emph{An introduction to probability theory and its applications,
  Volume 2}, volume~81.
\newblock John Wiley \& Sons, 1991.

\bibitem[Ferreira and Menegatto(2009)]{ferreira2009eigenvalues}
JC~Ferreira and VA~Menegatto.
\newblock Eigenvalues of integral operators defined by smooth positive definite
  kernels.
\newblock \emph{Integral Equations and Operator Theory}, 64\penalty0
  (1):\penalty0 61--81, 2009.

\bibitem[Fischer(2011)]{fischer2011polynomial}
Bernd Fischer.
\newblock \emph{Polynomial based iteration methods for symmetric linear
  systems}.
\newblock SIAM, 2011.

\bibitem[Flammarion and Bach(2015)]{flammarion2015averaging}
Nicolas Flammarion and Francis Bach.
\newblock From averaging to acceleration, there is only a step-size.
\newblock In \emph{Conference on Learning Theory}, pages 658--695. PMLR, 2015.

\bibitem[Fletcher and Reeves(1964)]{fletcher1964quadratically}
R.~Fletcher and C.~M. Reeves.
\newblock {Function minimization by conjugate gradients}.
\newblock \emph{The Computer Journal}, 7\penalty0 (2):\penalty0 149--154, 01
  1964.
\newblock ISSN 0010-4620.
\newblock \doi{10.1093/comjnl/7.2.149}.
\newblock URL \url{https://doi.org/10.1093/comjnl/7.2.149}.

\bibitem[{Fridman}(1962)]{fridman1952steepest}
V.~M. {Fridman}.
\newblock {On the convergence of methods of steepest descent type}.
\newblock \emph{{Usp. Mat. Nauk}}, 17\penalty0 (3(105)):\penalty0 201--204,
  1962.
\newblock ISSN 0042-1316.

\bibitem[Gilyazov and Gol'dman(2013)]{gilyazov2013regularization}
Sergei~Farshatovich Gilyazov and Natali{\^a}~L'vovna Gol'dman.
\newblock \emph{Regularization of ill-posed problems by iteration methods},
  volume 499.
\newblock Springer Science \& Business Media, 2013.

\bibitem[Hanke(1991)]{hanke1991accelerated}
Martin Hanke.
\newblock Accelerated landweber iterations for the solution of ill-posed
  equations.
\newblock \emph{Numerische mathematik}, 60\penalty0 (1):\penalty0 341--373,
  1991.

\bibitem[Hanke(1996)]{hanke1996asymptotics}
Martin Hanke.
\newblock Asymptotics of orthogonal polynomials and the numerical solution of
  ill-posed problems.
\newblock \emph{Numerical Algorithms}, 11\penalty0 (1):\penalty0 203--213,
  1996.

\bibitem[Hestenes and Stiefel(1952)]{Hestenes&Stiefel:1952}
M.~R. Hestenes and E.~Stiefel.
\newblock Methods of conjugate gradients for solving linear systems.
\newblock \emph{Journal of research of the National Bureau of Standards},
  49:\penalty0 409--436, 1952.

\bibitem[Hestenes(2012)]{hestenes2012conjugate}
M.R. Hestenes.
\newblock \emph{Conjugate Direction Methods in Optimization}.
\newblock Stochastic Modelling and Applied Probability. Springer New York,
  2012.
\newblock ISBN 9781461260486.
\newblock URL \url{https://books.google.ru/books?id=nc3cBwAAQBAJ}.

\bibitem[Inc.()]{Mathematica}
Wolfram~Research{,} Inc.
\newblock Mathematica, {V}ersion 13.2.
\newblock URL \url{https://www.wolfram.com/mathematica}.
\newblock Champaign, IL, 2022.

\bibitem[Jacot et~al.(2018)Jacot, Gabriel, and Hongler]{jacot2018neural}
Arthur Jacot, Franck Gabriel, and Cl{\'e}ment Hongler.
\newblock Neural tangent kernel: Convergence and generalization in neural
  networks.
\newblock \emph{arXiv preprint arXiv:1806.07572}, 2018.

\bibitem[Jin et~al.(2021)Jin, Banerjee, and Mont{\'u}far]{jin2021learning}
Hui Jin, Pradeep~Kr Banerjee, and Guido Mont{\'u}far.
\newblock Learning curves for gaussian process regression with power-law priors
  and targets.
\newblock \emph{arXiv preprint arXiv:2110.12231}, 2021.

\bibitem[Kammerer and Nashed(1971)]{kammerer1971steepest}
William~J Kammerer and M~Zuhair Nashed.
\newblock Steepest descent for singular linear operators with nonclosed range.
\newblock \emph{Applicable Analysis}, 1\penalty0 (2):\penalty0 143--159, 1971.

\bibitem[Kammerer and Nashed(1972)]{kammerer1972convergence}
William~J Kammerer and M~Zuhair Nashed.
\newblock On the convergence of the conjugate gradient method for singular
  linear operator equations.
\newblock \emph{SIAM Journal on Numerical Analysis}, 9\penalty0 (1):\penalty0
  165--181, 1972.

\bibitem[Kantorovich and Akilov(1964)]{kantorovich1964functional}
Leonid Kantorovich and Gleb Akilov.
\newblock \emph{Functional analysis in normed spaces}.
\newblock Number~46. Pergamon Press;[distributed in the Western Hemisphere by
  Macmillan, New York], 1964.

\bibitem[Karamata(1930)]{karamata1930certains}
Jovan Karamata.
\newblock {Sur certains 'Tauberian theorems' de M.M. Hardy et Littlewood}.
\newblock 1930.

\bibitem[Kingma and Ba(2014)]{kingma2014adam}
Diederik~P Kingma and Jimmy Ba.
\newblock Adam: A method for stochastic optimization.
\newblock \emph{arXiv preprint arXiv:1412.6980}, 2014.

\bibitem[Kopitkov and Indelman(2020)]{kopitkov2020neural}
Dmitry Kopitkov and Vadim Indelman.
\newblock Neural spectrum alignment: Empirical study.
\newblock In \emph{International Conference on Artificial Neural Networks},
  pages 168--179. Springer, 2020.

\bibitem[Krasnoselskii et~al.(1972)Krasnoselskii, Vainikko, Zabreiko,
  Rutitskii, and Stetsenko]{krasnoselskii52approximate}
MA~Krasnoselskii, GM~Vainikko, PP~Zabreiko, Ya~B Rutitskii, and V~Ya Stetsenko.
\newblock Approximate solutions of operator equations, noordhoff, groningen,
  1972.
\newblock \emph{MR}, 52:\penalty0 6515, 1972.

\bibitem[K{\"u}hn(1987)]{kuhn1987eigenvalues}
Thomas K{\"u}hn.
\newblock Eigenvalues of integral operators with smooth positive definite
  kernels.
\newblock \emph{Archiv der Mathematik}, 49\penalty0 (6):\penalty0 525--534,
  1987.

\bibitem[Lacotte and Pilanci(2020)]{lacotte2020optimal}
Jonathan Lacotte and Mert Pilanci.
\newblock Optimal randomized first-order methods for least-squares problems.
\newblock In \emph{International Conference on Machine Learning}, pages
  5587--5597. PMLR, 2020.

\bibitem[LeCun et~al.(2010)LeCun, Cortes, and Burges]{lecun2010mnist}
Yann LeCun, Corinna Cortes, and CJ~Burges.
\newblock Mnist handwritten digit database.
\newblock \emph{ATT Labs [Online]. Available:
  http://yann.lecun.com/exdb/mnist}, 2, 2010.

\bibitem[Lee et~al.(2019)Lee, Xiao, Schoenholz, Bahri, Novak, Sohl-Dickstein,
  and Pennington]{NEURIPS2019_0d1a9651}
Jaehoon Lee, Lechao Xiao, Samuel Schoenholz, Yasaman Bahri, Roman Novak, Jascha
  Sohl-Dickstein, and Jeffrey Pennington.
\newblock Wide neural networks of any depth evolve as linear models under
  gradient descent.
\newblock In H.~Wallach, H.~Larochelle, A.~Beygelzimer, E.~Fox, and R.~Garnett,
  editors, \emph{Advances in Neural Information Processing Systems}, volume~32.
  Curran Associates, Inc., 2019.
\newblock URL
  \url{https://proceedings.neurips.cc/paper/2019/file/0d1a9651497a38d8b1c3871c84528bd4-Paper.pdf}.

\bibitem[Lee et~al.(2020)Lee, Schoenholz, Pennington, Adlam, Xiao, Novak, and
  Sohl-Dickstein]{lee2020finite}
Jaehoon Lee, Samuel~S Schoenholz, Jeffrey Pennington, Ben Adlam, Lechao Xiao,
  Roman Novak, and Jascha Sohl-Dickstein.
\newblock Finite versus infinite neural networks: an empirical study.
\newblock \emph{arXiv preprint arXiv:2007.15801}, 2020.

\bibitem[Meurant and Strako{\v{s}}(2006)]{meurant2006lanczos}
G{\'e}rard Meurant and Zden{\v{e}}k Strako{\v{s}}.
\newblock The lanczos and conjugate gradient algorithms in finite precision
  arithmetic.
\newblock \emph{Acta Numerica}, 15:\penalty0 471--542, 2006.

\bibitem[Neal(2012)]{neal2012bayesian}
Radford~M Neal.
\newblock \emph{Bayesian learning for neural networks}, volume 118.
\newblock Springer Science \& Business Media, 2012.

\bibitem[Nemirovskiy and Polyak(1984{\natexlab{a}})]{nemirovskiy1984iterative1}
Arkadi~S Nemirovskiy and Boris~T Polyak.
\newblock {Iterative methods for solving linear ill-posed problems under
  precise information. I.}
\newblock \emph{{Izv. Akad. Nauk SSSR. Tekhn. Kibernet.}}, \penalty0 (2),
  1984{\natexlab{a}}.
\newblock [In Russian].

\bibitem[Nemirovskiy and Polyak(1984{\natexlab{b}})]{nemirovsky1984iterative2}
Arkadi~S Nemirovskiy and Boris~T Polyak.
\newblock {Iterative methods of solving linear ill-posed problems with precise
  information. II."}.
\newblock \emph{{Izv. Akad. Nauk SSSR. Tekhn. Kibernet.}}, \penalty0 (3),
  1984{\natexlab{b}}.
\newblock [In Russian].

\bibitem[Nesterov(2003)]{nesterov2003introductory}
Yurii Nesterov.
\newblock \emph{Introductory lectures on convex optimization: A basic course},
  volume~87.
\newblock Springer Science \& Business Media, 2003.

\bibitem[Nesterov(1983)]{nesterov1983method}
Yurii~E Nesterov.
\newblock A method for solving the convex programming problem with convergence
  rate $o(1/k^2)$.
\newblock In \emph{{Dokl. Akad. Nauk SSSR}}, volume 269, pages 543--547, 1983.

\bibitem[Nitanda and Suzuki(2021)]{nitanda2021optimal}
Atsushi Nitanda and Taiji Suzuki.
\newblock Optimal rates for averaged stochastic gradient descent under neural
  tangent kernel regime.
\newblock In \emph{International Conference on Learning Representations}, 2021.
\newblock URL \url{https://openreview.net/forum?id=PULSD5qI2N1}.

\bibitem[Pedregosa and Scieur(2020)]{pedregosa2020acceleration}
Fabian Pedregosa and Damien Scieur.
\newblock Acceleration through spectral density estimation.
\newblock In \emph{International Conference on Machine Learning}, pages
  7553--7562. PMLR, 2020.

\bibitem[Polyak(1964)]{polyak1964some}
Boris~T Polyak.
\newblock Some methods of speeding up the convergence of iteration methods.
\newblock \emph{Ussr computational mathematics and mathematical physics},
  4\penalty0 (5):\penalty0 1--17, 1964.

\bibitem[Polyak(1987)]{Polyak87}
Boris~T. Polyak.
\newblock \emph{Introduction to Optimization}.
\newblock Optimization Software, New York, 1987.

\bibitem[Polyak and Juditsky(1992)]{polyak1992acceleration}
Boris~T Polyak and Anatoli~B Juditsky.
\newblock Acceleration of stochastic approximation by averaging.
\newblock \emph{SIAM journal on control and optimization}, 30\penalty0
  (4):\penalty0 838--855, 1992.

\bibitem[Pronzato et~al.(2001)Pronzato, Wynn, and
  Zhigljavsky]{pronzato2001renormalised}
Luc Pronzato, Henry~P Wynn, and Anatoly~A Zhigljavsky.
\newblock Renormalised steepest descent in hilbert space converges to a
  two-point attractor.
\newblock \emph{Acta Applicandae Mathematica}, 67\penalty0 (1):\penalty0 1--18,
  2001.

\bibitem[Qian(1999)]{qian1999momentum}
Ning Qian.
\newblock On the momentum term in gradient descent learning algorithms.
\newblock \emph{Neural networks}, 12\penalty0 (1):\penalty0 145--151, 1999.

\bibitem[Ritter et~al.(1995)Ritter, Wasilkowski, and
  Wo{\'z}niakowski]{ritter1995multivariate}
Klaus Ritter, Grzegorz~W Wasilkowski, and Henryk Wo{\'z}niakowski.
\newblock Multivariate integration and approximation for random fields
  satisfying sacks-ylvisaker conditions.
\newblock \emph{The Annals of Applied Probability}, pages 518--540, 1995.

\bibitem[Steinwart et~al.(2009)Steinwart, Hush, Scovel,
  et~al.]{steinwart2009optimal}
Ingo Steinwart, Don~R Hush, Clint Scovel, et~al.
\newblock Optimal rates for regularized least squares regression.
\newblock In \emph{COLT}, pages 79--93, 2009.

\bibitem[Szego(1939)]{Szvego1939}
Gabor Szego.
\newblock \emph{Orthogonal polynomials}.
\newblock American Mathematical Society Providence, 4th ed. edition, 1939.
\newblock ISBN 0821810235.

\bibitem[Szeg{\"o}(1959)]{szego1959orthogonal}
Gabor Szeg{\"o}.
\newblock \emph{Orthogonal Polynomials}.
\newblock Number v. 23 in American Mathematical Society colloquium
  publications. American Mathematical Society, 1959.
\newblock ISBN 9780821889527.

\bibitem[Varre et~al.(2021)Varre, Pillaud-Vivien, and
  Flammarion]{varre2021last}
Aditya Varre, Loucas Pillaud-Vivien, and Nicolas Flammarion.
\newblock Last iterate convergence of sgd for least-squares in the
  interpolation regime.
\newblock \emph{arXiv preprint arXiv:2102.03183}, 2021.

\bibitem[Velikanov and Yarotsky(2021)]{velikanov2021explicit}
Maksim Velikanov and Dmitry Yarotsky.
\newblock Explicit loss asymptotics in the gradient descent training of neural
  networks.
\newblock \emph{Advances in Neural Information Processing Systems}, 34, 2021.

\bibitem[Velikanov et~al.(2022)Velikanov, Kuznedelev, and
  Yarotsky]{velikanov2022view}
Maksim Velikanov, Denis Kuznedelev, and Dmitry Yarotsky.
\newblock {A view of mini-batch SGD via generating functions: conditions of
  convergence, phase transitions, benefit from negative momenta}.
\newblock In \emph{The Eleventh International Conference on Learning
  Representations}, 2022.
\newblock URL \url{https://openreview.net/forum?id=bzaPGEllsjE}.

\bibitem[Widom(1963)]{widom1963asymptotic}
Harold Widom.
\newblock Asymptotic behavior of the eigenvalues of certain integral equations.
\newblock \emph{Transactions of the American Mathematical Society},
  109\penalty0 (2):\penalty0 278--295, 1963.

\bibitem[Williams and Rasmussen(2006)]{williams2006gaussian}
Christopher~K Williams and Carl~Edward Rasmussen.
\newblock \emph{Gaussian processes for machine learning}, volume~2.
\newblock MIT press Cambridge, MA, 2006.

\bibitem[Zou et~al.(2021)Zou, Wu, Braverman, Gu, and Kakade]{zou2021benign}
Difan Zou, Jingfeng Wu, Vladimir Braverman, Quanquan Gu, and Sham~M Kakade.
\newblock Benign overfitting of constant-stepsize sgd for linear regression.
\newblock \emph{arXiv preprint arXiv:2103.12692}, 2021.

\end{thebibliography}

\end{document}